\documentclass{amsart}
  \usepackage{amscd,amssymb,epsfig}
  \usepackage{graphicx}
   \usepackage{epic,eepic}

                     \numberwithin{equation}{subsection}

                     \newtheorem{propo}{Proposition}[subsection]
                     \newtheorem{corol}[propo]{Corollary}
                     \newtheorem{theor}[propo]{Theorem}
                     \newtheorem{lemma}[propo]{Lemma}
                     \theoremstyle{definition}

                     \theoremstyle{remark}

                     \newcommand{\ZZ}{\mathbb{Z}}

\newcommand{\A}{\mathcal{A}}

 \newcommand{\naw}{nanoword}
\newcommand{\naws}{nanowords}

\newcommand{\N}{\mathcal{N}}

\newcommand{\K}{\mathcal{K}}

\newcommand{\card}{\operatorname{card}}

                     \newcommand{\id}{\operatorname{id}}

\newcommand{\modu}{\operatorname{mod}}

\begin{document}
      \title{Topology of   words}
                     \author[Vladimir Turaev]{Vladimir Turaev}
                     \address{%
              IRMA, Universit\'e Louis  Pasteur - C.N.R.S., \newline
\indent  7 rue Ren\'e Descartes \newline
                     \indent F-67084 Strasbourg \newline
                     \indent France \newline
\indent e-mail: turaev@math.u-strasbg.fr }
                     \begin{abstract}  We introduce a topological approach to words.   Words are approximated by Gauss words and then  studied up to   natural modifications inspired by  homotopy transformations   of curves on the plane.
                     \end{abstract}
                     \maketitle

   \section {Introduction}
 Words are finite sequences of letters in a given  alphabet.
   Every word has its own personality  and
should be treated with the same respect  and attention  as say, a polyhedron or a   manifold.  In this paper we attempt to  study words as topological objects. A word in an  alphabet $\alpha$ can be viewed as a way of interaction or interlacement of  the letters of $\alpha$.   For example, two letters $a, b\in \alpha$ are interlaced in the word $abab$ and are not interlaced in the word $aabb$.    From this perspective, a  word  can be compared with  a link  of circles   in Euclidean
 3-space, the letters  being the counterparts of the circles.  
 
\begin{figure}
\centerline{\includegraphics[width=5cm]{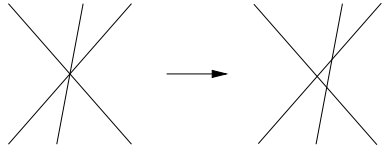}}
\caption{Desingularization}\label{figure1}
\end{figure}

Another  geometric viewpoint is suggested  by the Gauss-Rosenstiehl \cite{ro} correspondence between certain words and plane curves. This viewpoint consists in treating an arbitrary word in the alphabet $\alpha$  as a \lq\lq curve"   passing through the points of $\alpha$ in the prescribed order.  To make this work, one needs to embed $\alpha$ in a bigger space, say  a surface, or better to label certain points   of a surface with letters of  $\alpha$. These  ideas, albeit  imprecise at this stage, suggest   further directions of thought. First of all,   the letters appearing in the word with multiplicity $\geq 3$  are      \lq\lq singular" self-crossings of the curve. They  may  be  \lq\lq  desingularized"  to obtain only double self-crossings, see Figure \ref{figure1}.   
Secondly,   one may allow different points  of the ambient surface  to be labeled with the same letter  which leads us to  so-called \'etale words generalizing the ordinary words.  Thirdly, the notion  of  homotopy for   curves   leads us to a   notion of homotopy  for words.

The key  new concepts introduced in this paper are those of  \'etale words and nanowords. An  \'etale
word  over an alphabet $\alpha$  is a word in an alphabet $\A$ endowed with a projection       $\A\to \alpha$. 
The image of a letter $A\in \A$ under this projection  is denoted $\vert A\vert$. 
For instance, $ABCABC$   is an \'etale
word  over $\alpha$ provided we specify   $
\vert A\vert, \vert B\vert , \vert C\vert\in \alpha$.    Every word in the alphabet $\alpha$ becomes an \'etale word  over $\alpha$ by using the identity   $\id:\alpha\to \alpha$ as the projection.  
An \'etale word   in which every      letter      appears   twice or not at all is  called  a   nanoword.   For instance,   $ABCABC$    is a nanoword. Every  \'etale
word can be  
 approximated  by a nanoword  via a  desingularization which replaces each letter of multiplicity $m\geq 3$ by $m(m-1)/2$  letters of multiplicity 2.   This allows us to focus on nanowords; all definitions and results
concerning them extend to arbitrary \'etale words (and, in particular, to  arbitrary words) via   desingularization.

\begin{figure}
\centerline{\includegraphics[width=8cm]{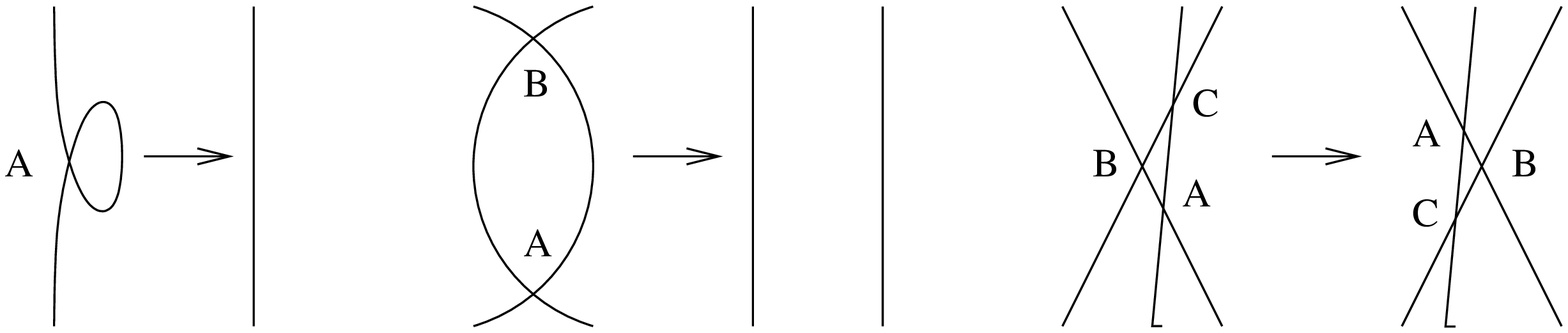}}
\caption{Local deformations of plane curves}\label{figure2}
\end{figure}

We now fix an additional piece of data:   an  involution
$\tau :\alpha\to \alpha$ (it may be the identity).    Given $\tau$, we introduce   an
equivalence relation of   homotopy   on the set of  \'etale words  over  $\alpha$.     It is implied  that homotopic \'etale words   give rise to the same 
interlacement of    letters of $\alpha$  with respect to $\tau$.        
 The relation of homotopy is generated by three  transformations  or moves  on nanowords. The first move consists in deleting two consecutive  entries of the
same   letter.  The second move has the form $xABy BA z\mapsto xyz$ where $x,y,z$ are words and $A,B$ are   letters such that $\tau (\vert A \vert)=\vert B\vert$.  The third move has the form  $xAByACzBCt\mapsto xBAyCAzCBt$ where   $x,y,z,t$ are words and $A,B,C$ are   letters such that $
\vert A\vert= \vert B\vert = \vert C\vert$.  These moves are   suggested by 
the  standard local deformations of plane curves, see Figure \ref{figure2}.

  We shall study properties and characteristics of words   and \'etale words   preserved under homotopy.  A number of methods
developed in   topology  find an echo in this setting.     The homological intersection theory of curves on surfaces suggests a family of  homotopy invariants of nanowords  (the self-linking function, the  linking  form,   the linking pairing).    The technique of Alexander matrices in knot theory leads us  to modules and polynomial invariants of nanowords including a rather powerful invariant  $\lambda$.  
The technique of  colorings of knot diagrams suggests a notion   of a coloring for nanowords.  The theory of knot quandles   finds its analogue in the form of   $\alpha$-keis.  The theory of virtual strings leads to so-called characteristic sequences of words.  Other  methods of low-dimensional topology may possibly apply in this setting.

As an application of our algebraic methods, we give a homotopy classification of  nanowords of length $\leq 6$.  All nanowords of length 2 are contractible, i.e., homotopic to an empty nanoword. The classification of nanowords of length $4$ is as follows:    the   nanowords of the form $AABB$ or $ABBA$  are  contractible (by the first move); a nanoword of the form $ABAB$  is contractible if and only if $\vert A\vert =\tau (\vert B \vert)$; two non-contractible nanowords $ABAB$ and $A'B'A'B'$ are homotopic if and only if $\vert A\vert=\vert A'\vert$  and $\vert B\vert=\vert B'\vert$ (Theorem \ref{classid}). The homotopy classification of    nanowords of length $6$ is more involved. Excluding those   homotopic to nanowords of  length $\leq 4$ by the first move, we obtain  5 families  of nanowords of the form $$ABCABC, \,ABCACB, \,ABCBAC, \,ABCBCA, \,ABACBC.$$ We fully describe when such nanowords are homotopic to nanowords of length $\leq 4$ and when they   are homotopic to each other  (Theorem \ref{1classid}).   

We also give a homotopy classification of   words  of length $\leq 5$ in the alphabet $\alpha$ (Theorems \ref{erf4444sid} and  \ref{erfclassid}).  For  $\tau=\id:\alpha\to \alpha$ this classification is especially simple:  a     word  of length $\leq 5$, in which every letter appears at least twice or not at all,   is non-contractible if and only if it  has one of the following six forms: 
 $abab,   abaab, baaba,   
 aabab,  ba baa,  ababa$ with distinct $a,b\in \alpha$. Two  words from this list are  homotopic   if and only if they coincide letter-wise.

  In the case where the alphabet $ \alpha$  consists of two elements   permuted by $\tau$, the  notion of a nanoword over $\alpha$ is equivalent to the notion of an open virtual string introduced in \cite{tu1}.  The results of the present paper generalize a number of results of \cite{tu1} and of the subsequent paper by   D. Silver and S. Williams \cite{sw} on     open virtual strings.  In particular, the invariant $\lambda$ and the characteristic sequences generalize the invariants  of open strings 
  introduced by Silver and   Williams.

To make this paper   accessible to readers not interested in topology, it is written    in  purely algebraic terms. Relations with  topology will be discussed elsewhere.  

The papers consists of 4 parts. The first part (Sections 2 -- 5) is  devoted to the basics. We define   \'etale words, nanowords, desingularization of \'etale words, homotopy and coverings of nanowords, and discuss a group-theoretic approach  to nanowords. In the second part (Sections 6 -- \ref{cl555v}) we  construct    linking invariants of nanowords  and  give a homotopy classification of nanowords of length 4 and of words of length $\leq 5$.  
The third part  (Sections \ref{col56822}  --  \ref{dqmqmqmqmqm2}) begins with a discussion of colorings of nanowords and proceeds to    modules and polynomials of nanowords.  Then  we give a homotopy classification of nanowords of length 6.   In the fourth  part  (Sections \ref{4dd822} -- \ref{dfdfdfdfd}) we introduce      $\alpha$-keis  and characteristic sequences of nanowords. They are used to distinguish certain nanowords of length 6 from nanowords of length 4 and to accomplish thus the homotopy classification of  nanowords of length $\leq 6$.

 Throughout the paper the symbol $\alpha$ denotes a   set endowed with an involution $\tau:\alpha\to \alpha$.


 \section{\'Etale words and {\naw}s}

\subsection{Words}\label{word}    An {\it alphabet}  is a    set and   
{\it 
letters} are its elements.     A {\it   word of length} $n\geq 1$  in an alphabet  $\alpha$ is  a mapping 
$w: \hat n\to \alpha  $   where $\hat
n=\{1,2,...,n\}$.    Such a 
word
 is encoded by the  sequence $w(1) w(2) \cdots  w(n)$.     For example, the sequence
$aba   $ in the alphabet $\alpha=\{a,b\}$ encodes  the word  $ \hat 3 \to \alpha $ sending $1, 2, 3 $
to
$a, b, a $ respectively.   By definition, there is a   unique  {\it empty word} $\emptyset$  of length 0.

 Writing down consecutively the letters of two  words $w$ and   $v$ we 
obtain
their concatenation $wv$. For instance,   the concatenation of  $w=abb$ 
and
$v=aa$ is the word $wv=abbaa$. Writing the letters of a word $w$ in the opposite order  we obtain the {\it opposite word} $w^-$. For instance, if $w=abb$,
then $w^-=bba$.

 The {\it multiplicity} of a letter $a\in \alpha$ in  a word $w: \hat n\to \alpha  $ is the number   $m=m_w(a)=\card\,\{i\in \hat n\,\vert\, w(i)=a \}$. We say   that   
  $a $ appears $m$ times in  $w $. For example,  $m_{aab}(a)=2$ and $m_{aab}(b)=1$.
  
  A {\it monoliteral} word of length $m\geq 2$ is the word
 $aa\cdots a$ formed by $m$ copies of  the same letter $a\in \alpha$.   This word is denoted $a^m$.
     
A mapping $f$ from a set $\alpha_1$ to  a set $\alpha_2$ induces a mapping $f_\#$ from the set of words in the alphabet 
$\alpha_1$ to  the set of words in the alphabet  $\alpha_2$. It is obtained by applying $f$ letterwise.  

 \subsection{\'Etale words}\label{nano} The  class of    words in the alphabet $\alpha$ is too narrow for our purposes.
We  introduce  here  a wider class of  \'etale words over     $\alpha$. 

An {\it   $\alpha$-alphabet} is a set $\A$ endowed with  a  mapping  $ \A\to \alpha$ 
called {\it projection}. The image  of a letter $A\in
\A$   under the projection is denoted $\vert A\vert$.   Any subset $\A'$ of an  
$\alpha$-alphabet $\A$ becomes an   $\alpha$-alphabet   by
restricting the projection  $\A\to \alpha$   to $\A'$.

A {\it morphism}  of   $\alpha$-alphabets $\A_1$,  $\A_2$ is a set-theoretic
mapping $f:\A_1\to
\A_2$ such that    $\vert A\vert=\vert f(A)\vert$ for all $A\in \A_1$. If $f$ is bijective, then this
morphism is an {\it  isomorphism}.    

An {\it \'etale word} over   
$\alpha$ is a pair (an $\alpha$-alphabet
$\A$,   a  word  
 in  the alphabet $\A$).  Two \'etale words    ($\A_1$, 
$w_1$) and  ($\A_2$, 
$w_2$)  over $\alpha$ are {\it   isomorphic} if there is an isomorphism    $f:\A_1\to \A_2$ 
  such that   $w_2=f_\# (w_1)$.  The relation of isomorphism for  \'etale words is denoted  $\approx$.

We define a product of \'etale word  ($\A_1$, 
$w_1$) and  ($\A_2$, 
$w_2$)  over $\alpha$  as follows. Replacing if necessary  ($\A_1$, 
$w_1$) with an isomorphic \'etale word we can assume that $\A_1\cap \A_2=\emptyset$.  Then the product in question is the  \'etale word  $(\A_1 \cup
\A_2, w_1 w_2)$ over $\alpha$.  It is well defined up to isomorphism. Multiplication of  \'etale words    is associative and has a unit represented by an empty
\'etale word in an empty $\alpha$-alphabet.

For each \'etale word $(\A,w)$, we   have the {\it opposite} \'etale word    $(\A, w^-)$ where $w^-$ is the word opposite to $w$.
 For each $\alpha$-alphabet $(\A, p:\A\to \alpha)$,    the {\it inverse} 
$\alpha$-alphabet    $\overline \A$  is the same set $\A$ with projection
$\tau p:\A\to \alpha$ where    $\tau:\alpha\to \alpha$ is the fixed involution.   The {\it inverse} of an  \'etale word $(\A,w)$ over $\alpha$ is the \'etale word $\overline w= (\overline \A, w)$. An   \'etale word $w$   is {\it symmetric} (resp.\ {\it skew-symmetric}) if it is isomorphic to $ w^-$ (resp.\ to $\overline w^-$). For instance,   $ v   v^-$ is symmetric and $v \overline v^-$  is skew-symmetric  for any  \'etale word $v$. 
 
  Each word $w$ in the alphabet $\alpha$ gives rise to an  \'etale
word
$(\A=\alpha,w)$ over $\alpha$ where the projection $\A\to \alpha$ is the identity.  In this way,   \'etale words over $\alpha$ generalize 
  words in the alphabet $\alpha$. 

Warning:  concatenation of words in the alphabet $\alpha$ differs from multiplication of the corresponding \'etale words. 
For instance, for the words $abb$ and $aa$ in the alphabet $\alpha=\{a,b\}$,  the corresponding \'etale words     are $(\A_1=\{A,B\}, ABB)$ with $\vert A\vert
=a,
\vert B
\vert =b$ and $(\A_2=\{C\}, CC)$ with $\vert C\vert
= a$. Their  product  is the \'etale word $(\{A,B,C\}, ABBCC)$ with  $\vert A\vert = \vert C\vert
=a,
\vert B
\vert 
= b$. On the other hand, the \'etale word corresponding to $abbaa$ is 
$(\{A,B\}, ABBAA)$ with  $\vert A\vert 
=a,
\vert B
\vert 
= b$.

     \subsection{Gauss words and nanowords}\label{gawo} A  word $w$ in a finite alphabet $\A$ is a {\it Gauss word} if    every letter of
$\A$  appears in $w$ exactly twice.   For instance, 
$ABAB
$ is a Gauss word in the alphabet $\{A,B\}$ while $ABA $ and  $AA$ are not Gauss words in this alphabet.  
If   $\A$    consists of $N$
letters, then there are  $  (2N) !$ Gauss words in the alphabet $\A$. Note that concatenation of two or more Gauss words
 in a non-empty alphabet is never a Gauss
word.

  A {\it nanoword} is an  \'etale Gauss word.  More precisely, 
a {\it {\naw} over}   
$\alpha$ is a pair (a finite $\alpha$-alphabet
$\A$,   a Gauss word  
 in  the alphabet $\A$).  The {\it length} of this  nanoword is the length  of the Gauss word in question, i.e.,      $2 \card (\A)$.

The set of
nanowords over
$\alpha$ is denoted  
$\N(\alpha)$. Two nanowords     are {\it   isomorphic} if they are isomorphic as
\'etale words.  Observe that the  \'etale word  opposite or inverse to a nanoword is  a nanoword. 
The product of two nanowords is a nanoword (defined up to isomorphism).   An   empty
\'etale word in an empty $\alpha$-alphabet is a nanoword called the {\it empty nanoword}.

Instead of writing $(\A,w)$ for a nanoword over $\alpha$, we shall often   write simply $w$. The alphabet $\A$ can be uniquely recovered from   $w$
as the set of all letters appearing in $w$. However, the projection $\A\to \alpha$ should be always specified. 

Each Gauss word $w$ in a finite alphabet $\alpha$ gives rise to a nanoword $(\A=\alpha,w)$ over $\alpha$ where the projection
$\A\to
\alpha$ is the identity.  In this way, nanowords  over $\alpha$ generalize  Gauss words  in the alphabet $\alpha$.

 \subsection{Example}\label{rexacrr}   Let $\alpha=\{a,b\}$ with involution permuting $a$ and $b$. Let  $\A=\{A_1, A_2, B_1\}$ with $\vert A_1\vert
=\vert A_2\vert =a, \vert B_1\vert =b$. Then
$(\A,   A_1 B_1 A_1 A_2 A_2 B_1 )$ is a nanoword over $\alpha$.  The opposite nanoword is  $(\A,   B_1 A_2 A_2 A_1 B_1  A_1  
 )$ and the inverse   is  $(\overline \A,   A_1 B_1 A_1 A_2 A_2 B_1 )$ where  $\overline \A=\{A_1, A_2, B_1\}$ with $\vert A_1\vert =\vert A_2\vert
=b,
\vert B_1\vert =a$.

\subsection{Desingularization}\label{eteahomot}  For an \'etale word $(\A, w )$ over $\alpha$, we define a nanoword
$(\A^d, w^d)$ over $\alpha$ called the {\it desingularization} of $w$.    The alphabet
$\A^d$ consists of the triples
$(A,i,j)$ where
$A\in \A$  and   $1\leq i <j \leq m_w(A)$ where $m_w(A)$ is the multiplicity of $A$ in $w$.  For   brevity, we   write $A_{i,j}$ for  $(A,i,j)$.  We make $\A^d$ into an $\alpha$-alphabet
by $\vert A_{i,j}\vert =\vert A\vert \in \alpha$ for all $A, i,j$.  The word
$w^d$ is obtained  from $w$ by  first
deleting   all $A\in \A$ with $m_w(A)=1$.  Then for each    $A\in \A$ with $m_w(A)\geq 2$ and each $i=1,2,..., m_w(A)$,  we replace the $i$-th entry of $A$ in
$w$ by 
$$A_{1,i} A_{2,i}\cdots A_{i-1,i} A_{i, i+1} A_{i, i+2} \cdots A_{i, m_w(A)}.$$ 
The resulting word $w^d$ in the alphabet $\A^d$  is a Gauss word so that   $(\A^d,w^d)$ is a nanoword of length   $ \sum_{A\in \A} m_w(A) (m_w(A)-1) $.
For example,  if $w=AABABC$, then
 $\A^d=\{A_{1,2},A_{1,3},  A_{2,3},   B_{1,2} \}$ 
and
 $w^d=A_{1,2} A_{1,3}   A_{1,2} A_{2,3}   B_{1,2} A_{1,3} A_{2,3}   B_{1,2} $. 
 
 The desingularization may considerably increase the length   of a word   but this is the price to pay for obtaining a Gauss word.
 If $w$ is a nanoword,  then   $  w^d \approx  w $.
If $w=\emptyset$, then   $ w^d=\emptyset$.
It follows from the definitions that   $(w_1w_2)^d= w_1^d w_2^d$ for any \'etale words $w_1, w_2$,   $\overline w^d ={\overline w}^d$,   and  
 $(w^-)^d \approx (w^d)^-$ for any  $w$. (The latter isomorphism    sends $A_{i,j}$ to $A_{m_w(A)+1-j, m_w(A)+1-i}$ for all $A,i,j$).

\section{Homotopy}\label{456FGH}

\subsection{Homotopy moves}\label{honawo} To define   homotopy
  of nanowords we fix a  {\it homotopy data} consisting of a
set $\alpha$ with involution $\tau:\alpha\to \alpha$  and a  set
$S\subset \alpha^3= \alpha\times \alpha \times \alpha$. Given this
  data, we define three transformations of {\naws} over
$\alpha$ called {\it $S$-homotopy moves} or, when $S$ is fixed, {\it homotopy moves}.

(1). The first move  applies to any nanoword   of the form    $(\A,   x 
A
A  y)$ where $A\in \A $ and $x,y$ are words  in the alphabet 
$\A'=\A-\{A\}$. It
transforms $(\A,   x A
A  y)$  into
the nanoword
$(\A',  xy)$ where  the structure of an $\alpha$-alphabet in $\A' $ is 
obtained
by restricting the  one  in  $\A$.  Note that
$ x  y$ is   a Gauss word in the alphabet $\A'$.

 The inverse move
$(\A',   xy)\mapsto (\A'\cup \{A\} ,   x A  A  y)
$    adds a  new   letter
$A $ to   $\A' $  with   arbitrary  $\vert A\vert \in \alpha$ and   
replaces
the  Gauss word
$xy$ in  the $\alpha$-alphabet $\A'$ with $xAAy$.

(2). The second move  applies to a  nanoword  of the form $(\A ,   xA  
B  y BAz
 )$   where $A,B\in  \A$ with
$\vert B\vert =\tau ( \vert A\vert)$ and $x,y,z$ are words in the 
alphabet 
$\A'=\A-\{A,B\}$.  This  nanoword is transformed   into
$(\A' ,   xyz)$ where  the structure of an $\alpha$-alphabet in $\A' $ 
is
obtained by restricting
the  one in  $\A$.

The inverse move  $(\A', xyz)\mapsto  (\A'\cup \{A,B\},   xA
B  y BAz)$      adds   two  new letters
$A , B$ to  $\A' $  with   arbitrary $\vert A\vert \in \alpha$ and       
$\vert
B\vert =\tau ( \vert A\vert)$ and    replaces the Gauss
word
$xyz$ in  the alphabet $\A'$ with   $xA
B  y BAz$.

(3)   The third move  applies to a  nanoword  of the form $(\A, 
xAByACzBCt)$  
where
 $A,B,C \in \A$ are distinct letters with
$(\vert A\vert ,\vert B\vert , \vert C\vert)\in S $ and $x,y,z,t$
are words in the alphabet $\A-\{A,B,C\} $ such that $xyzt$ is a
Gauss word in this alphabet. The
  move  transforms  $(\A, xAByACzBCt)$  into  $ (\A, xBAyCAzCBt)$. The 
inverse
move
transforms  $ (\A, xBAyCAzCBt)$ into $(\A, xAByACzBCt)$. It
applies if   $(\vert A\vert ,\vert B\vert , \vert
C\vert)\in S $.

\subsection{Homotopy  of {\naw}s}\label{homot}  {\it Homotopy} is the
equivalence relation in the class of nanowords
generated by   isomorphism and the homotopy moves. For a  homotopy
data $(\alpha,   S)$,    two nanowords over $\alpha$ are {\it
  $S$-homotopic} if they can be
obtained from each other by a finite sequence of    $S$-homotopy moves
(1) -- (3), the inverse moves, and isomorphisms.  The relation of
$S$-homotopy is denoted $\simeq_{S}$. A nanoword $S$-homotopic to an
empty nanoword $\emptyset$ is  said to be {\it $S$-contractible}.
All these notions crucially depend on the choice of $\tau :\alpha
\to \alpha$ and   $S$.

The set of $S$-homotopy classes of nanowords over  $\alpha$ is
denoted $\N^S_\bullet(\alpha) =\N(\alpha)/\simeq_S$.
Multiplication of nanowords   is compatible with $S$-homotopy: if
nanowords $w_1, w_2$ are $S$-homotopic to nanowords $w'_1, w'_2$,
respectively, then $w_1w_2$ is $S$-homotopic to $ w'_1 w'_2$.
Multiplication of nanowords   makes $\N^S_\bullet(\alpha)$ into a
monoid with unit represented by $\emptyset$. If $S$ is invariant
under the involution $\tau \times \tau \times \tau$ on $\alpha^3$,
then nanowords inverse   to $S$-homotopic nanowords are themselves
$S$-homotopic. If $S$ is invariant under the involution
$(a,b,c)\mapsto (c,b,a)$ on $\alpha^3$, then nanowords opposite
to $S$-homotopic nanowords are $S$-homotopic.

A nanoword $w$   is {\it $S$-homotopically symmetric} (resp.\ {\it
$S$-homotopically skew-symmetric})  if   $w \simeq_S  w^-$ (resp.\
if   $w \simeq_S \overline w^-$). Clearly,   (skew-)  symmetric
nanowords  are $S$-homotopically (skew-) symmetric.
$S$-contractible nanowords are   $S$-homoto\-pically symmetric
and   skew-symmetric.  If nanowords $w  , w'$ are
$S$-homotopic, then $w' w^-$ is $S$-homotopically symmetric and
$w'\, \overline w^-$ is $S$-homotopically skew-symmetric.

For a  nanoword $w$,  let $\vert\vert w\vert \vert_S$ be half of
the minimal length of a nanoword homotopic to $w$.  This defines a
$\ZZ$-valued \lq\lq norm" $\vert \vert \cdot \vert \vert_S$ on
$\N^S_\bullet(\alpha)$ whose value on an $S$-homotopy class of
nanowords is half of the minimal length of a nanoword in this
class. Clearly, $\vert \vert w \vert \vert_S=0$ if and only if $w$
is $S$-contractible. Since all nanowords of length 2 are
$S$-contractible, the function $\vert \vert \cdot \vert \vert_S$
does not take  the value 1.   It is obvious that  $\vert \vert w_1
w_2 \vert \vert_S\leq \vert \vert w_1 \vert \vert_S + \vert \vert
w_2 \vert \vert_S$ for any  $w_1,w_2$.

\begin{lemma}\label{2edcolbrad}   Let $A,B,C $ be distinct letters in  
an
$\alpha$-alphabet $ \A$ and let $x,y,z,t$ be words in the alphabet
$\A-\{A,B,C\} $ such
that
$xyzt$ is a Gauss word in this alphabet. Then
$$ (\A, xAByCAzBCt) \simeq_S (\A, xBAyACzCBt) \, \,\,  { if}\,\,\, 
(\vert
A\vert , \tau (\vert B\vert) ,
\vert C\vert )\in S,  \leqno (i)$$
$$ (\A, xAByCAzCBt) \simeq_S  (\A, xBAyACzBCt)\, \,\,  { if}\,\,\, 
(\tau(\vert
A\vert ), \tau(\vert B\vert
), \vert C\vert)\in S, \leqno (ii)$$
$$ (\A, xAByACzCBt) \simeq_S  (\A, xBAyCAzBCt)\, \,\,  {  if}\,\,\, 
(\vert
A\vert , \tau(\vert B\vert ), \tau
(\vert C\vert))\in S. \leqno (iii)$$ \end{lemma}
                     \begin{proof}     Pick two symbols  $D,E $   not 
belonging
to $\A$.

We   verify (i). Set    $\vert D\vert =  \vert B\vert\in \alpha$
and $ \vert E\vert = \tau (\vert B\vert) \in \alpha$. Applying an
inverse (2)-move inserting $\cdots DE \cdots  ED \cdots $ and an
inverse (3)-move we obtain
$$(\A, xAByCAzBCt) \simeq_S   (\A\cup
\{D,E\}, xD\underline {EA}B y\underline {CA}zB\underline {CE}D t)   $$
$$
 \simeq_S   (\A\cup \{D,E\}, xDA  {EB}yACz {BE}CDt)$$
 where  for   visual convenience   the letters modified by the
inverse  (3)-move  are underlined. Here we use that  $(\vert
A\vert ,  \vert E\vert  , \vert C\vert )\in S$. Applying to the
resulting nanoword a (2)-move    removing $B, E$ we obtain the
nanoword $(\A\cup \{D\}-\{B\}, xDA yACz CDt)$. It is isomorphic to
$ (\A, xBAyACzCBt)$.

 We   verify (ii).  Set    $\vert D\vert =  \tau (\vert A\vert)$
 and $ \vert E\vert =  \vert A\vert$.    Applying        an inverse
(2)-move  inserting $\cdots DE \cdots  ED \cdots $ and an  inverse
to the homotopy (i)    we obtain
$$(\A, xAByCAzCBt) \simeq_S   (\A\cup \{D,E\}, xA\underline 
{BD}EyE\underline
{DC}Az\underline {CB}t)  $$
$$
 \simeq_S   (\A\cup \{D,E\}, x\underline {AD}BEyEC\underline {DA}zBCt)  
\simeq_S  (\A\cup \{E\}-\{A\}, xBE yECz BCt)  .$$
The latter nanoword is isomorphic to $ (\A, xBAyACzBCt)$.

To prove (iii), set    $\vert D\vert =   \vert C\vert$
 and $ \vert E\vert =  \tau (\vert C\vert)$. Then
$$xAByACzCBt\simeq_S x\underline {AB}yD \underline {E}  \underline {A}
CzC\underline {BE}Dt\simeq_S   xBAyDA\underline {EC}z\underline {CE} 
BDt$$
$$
\simeq_S
 xBAyDA z  BDt  \approx
xBAyCAzBCt.$$
\end{proof}

\begin{lemma}\label{1edcolbrad}  Suppose that $S\cap (\alpha \times b 
\times b)
\neq \emptyset$ for all $b\in  \alpha$. Let $(\A, xAByABz)$ be a 
nanoword over
$\alpha$ where  $A,B\in
\A$ with $\vert B\vert =\tau ( \vert A\vert)$ and $x,y,z$ are
words in $\A-\{A,B\}$.   Then $ (\A, xAByABz) \simeq_S  ( \A
-\{A,B\}, xyz)$. \end{lemma}
    \begin{proof}  Set $b=\vert B\vert\in \alpha$. By assumption, there 
is
$e\in \alpha$ such that $(e,b,b)\in S$.
    Pick a symbol $E$ not belonging to $\A$ and set    $   \vert
E\vert = \tau(e)$.     Then
$$(\A, xAByABz) \simeq_S   (\A\cup \{E\}, x\underline {AE}\, \underline
{EB}y\underline {AB}z)   \simeq_S   (\A \cup \{E\}, xE\underline
{AB}Ey\underline {BA}z)$$
$$  \simeq_S   ( \A  \cup \{E\} -\{A,B\}  , xEE y z)  \simeq_S   (
\A -\{A,B\}, x yz) .$$ Here the second homotopy is provided by
item  (ii) of the previous lemma  where   $A,B,C$ are replaced
with $E,A,B$, respectively. We use   that $(\tau (\vert E\vert),
\tau(\vert A\vert), \vert B\vert)=(e,b,b)\in S$.
\end{proof}

 \subsection{Homotopy of \'etale words}\label{homsmsmsomot} We say that  
\'etale words $w_1, w_2$ over $\alpha$ are
  {\it $S$-homotopic}  and write $w_1
\simeq_S w_2$ if   $w_1^d\simeq_S w_2^d$.    Clearly,  isomorphic
\'etale words are $S$-homotopic. $S$-homotopy is  an equivalence
relation in the class of \'etale words. For nanowords, this
relation of $S$-homotopy coincides with the one  defined above.
Every $S$-homotopy invariant  $I$ of nanowords extends to  an
$S$-homotopy invariant $I^d$ of   \'etale words
 by
$I^d (w)=I(w^d)$.   In particular, the  {\it $S$-norm}  of an
\'etale word $w$ is defined by $\vert \vert w \vert \vert_S=\vert
\vert w^d \vert \vert_S$.

 An \'etale word $w$ is {\it $S$-contractible} (resp.\ {\it 
$S$-homotopically
(skew-) symmetric}) if
the nanoword $w^d$ is $S$-contractible (resp.\ $S$-homotopically
(skew-) symmetric). Since    words in the alphabet $\alpha$ may be
treated as \'etale words,  we can apply to them all these
definitions.

 Although there are various interesting choices for $S$, we focus  here  on 
  one
simplest $S$.
 In the remaining part of the paper we     assume that $S$ is the 
diagonal of
$\alpha^3$ that is
 $S=\{(a,a,a)\}_{a\in \alpha}$. In other words  we allow the third 
homotopy
move $(\A, xAByACzBCt)\mapsto  (\A, xBAyCAzCBt)$
 and the inverse move   if and only if $\vert A\vert=\vert B\vert 
=\vert
C\vert$.
 Under this convention, we shall omit the prefix $S$- and speak simply 
of
homotopy rather then $S$-homotopy.
 We shall also omit the  index $S$ and
 write $\simeq, \vert \vert . \vert \vert, \N_\bullet(\alpha) $
 for $\simeq_S, \vert \vert . \vert \vert_S,  \N^S_\bullet(\alpha)
 $. Note that $S\cap (\alpha \times b \times b) =(b,b,b)$ for   $b\in
 \alpha$ so that we can apply Lemma \ref{1edcolbrad}.

 The principal aim of the homotopy theory of nanowords is  a 
computation of
$\N_\bullet(\alpha)$.
  As a step in this direction, we shall    construct several
homotopy invariants of nanowords.   They will often allow us to
distinguish nanowords up to homotopy and to decide whether a given
nanoword  is homotopically (skew-) symmetric or homotopic to its
inverse.   A  related  problem  is a  computation or at   least an
estimate  for    $\vert \vert \cdot \vert \vert$.

The following lemma shows that for words in a given alphabet, the
relation of homotopy does not change under a passage to a bigger
alphabet.

\begin{lemma}\label{1f559317olbrad}   Let $ \beta$ be a 
$\tau$-invariant 
subset of $\alpha$.
  If two \'etale words over $\beta$ are homotopic  in the class of 
\'etale
words over $\alpha$, then they are homotopic in the class of \'etale 
words over
$\beta$. \end{lemma}
    \begin{proof} Given
an \'etale word $(\A, w)$ over $\alpha$ we define its {\it pull-back} 
to
$\beta$  to be the  \'etale word over
$\beta$ obtained by deleting from both $\A$ and $w$ all letters $A\in 
\A$ with 
$\vert A\vert \in
\alpha-
\beta$.   The pull-back commutes with desingularization and  transforms
isomorphic (resp.\   homotopic) \'etale words  over $\alpha$  into 
isomorphic
(resp.\
homotopic) \'etale words over $\beta$.

Observe that any equivariant mapping of  sets with involution 
$f:\alpha'\to
\alpha$  induces a mapping
$f_\#:\N(\alpha')\to
\N(\alpha)$  by
$f_\# (\A, p:\A\to \alpha',w)= (\A, fp:\A\to \alpha,w)$.  Clearly,
$f_\#$ transforms isomorphic (resp.\   homotopic) nanowords   into 
isomorphic
(resp.\
homotopic) nanowords.  Denote the induced mapping  $ \N_{\bullet} 
(\alpha') 
\to  \N_{\bullet} (\alpha) $ by $f_\bullet$.  The composition of  
$f_\bullet
$ with the pull-back to $f(\alpha')\subset \alpha$ is the mapping $
\N_{\bullet} (\alpha')  \to  \N_{\bullet} (f(\alpha)) $ induced by $f$. 
If $f$
is
injective, then the latter
mapping is    an isomorphism and therefore $f_\bullet$ is
injective.  Applying this to the inclusion $\beta \hookrightarrow 
\alpha$, we
obtain the claim of the lemma.
  \end{proof}
  
  \subsection{Examples}\label{exithhhg}   1. The nanoword  $(\A=\{A\},    
A
A   )$ is contractible  for any choice of $\vert A\vert \in \alpha $.  
The
nanowords
$  ABBA  $ and $AABB$  (in the alphabet $\{A,B\}$) are  contractible 
for any  
$\vert A\vert, \vert B\vert \in \alpha $.
The nanoword
$  ABAB$ is  contractible  provided  $\vert A\vert =\tau (\vert B
\vert)$.     The nanoword  $ ABCDABCD $  (in the alphabet 
$\{A,B,C,D\}$)
 is  contractible  provided  either $\vert A\vert =\tau (\vert B
\vert)$ and  $\vert C\vert =\tau (\vert D\vert)$ or $\vert B\vert =\tau 
(\vert
C\vert)$ and  $\vert A\vert =\tau (\vert D\vert)$.

2. A more sophisticated example:  the nanoword $ABCDEFBGDHFAGCHE$
is contractible provided $\vert A\vert = \vert B\vert = \vert F\vert 
=\tau
(\vert C
\vert) =\tau (\vert D
\vert)  =\tau (\vert G
\vert)$ (the elements $\vert E\vert, \vert H\vert \in \alpha$ may be
arbitrary). Indeed, applying   Lemma \ref{2edcolbrad}(iii),
Lemma \ref{1edcolbrad}, and the first homotopy move  we obtain
$$A\underline {BC} DEF\underline {BG}DHFA\underline {GC}HE \simeq 
\underline
{AC}BDEFGBDHF\underline {AC}GHE$$
$$\simeq BDE\underline {FG}BDH\underline
{FG}HE \simeq \underline {BD}E \underline{BD} {HH}E \simeq EHHE \simeq 
EE
\simeq \emptyset.$$
Note that there is no way to apply to the nanoword
$ABCDEFBGDHFAGCHE$ the first or second homotopy moves; to contract this
nanoword one needs the third homotopy move.

3.    Let $a\in \alpha$ be a letter such that $\tau(a)= a$.  We show 
that the
monoliteral  words $a^3, a^4, a^5$  are
contractible. The  \'etale word corresponding to $a^3=aaa$ is 
$(\A=\alpha,
aaa)$ and the
desingularization yields the nanoword
 $   -\!12\!-\!13\!-\!12\!-\!23\!-\!13\!-\!23 $ where  we write $-ij$ 
for  $
a_{i,j} $. We have $$ \underline{ -\!12\!-\!13} \, 
\underline{-\!12\!-\!23}  \,
 \underline{-\!13\!-\!23}    \simeq \!-\!13 \underline{-\!12\!-\!23}  
\,
\underline{-\!12\!-\!23} \!-\!13\simeq \!-\!13\!-\!13 \simeq 
\emptyset.$$
 Here is a contracting homotopy for the desingularization of 
$a^4=aaaa$:   $$
\underline{-\! 12 \!-\!13} \!-\!14\underline{-\!12\!-\!23}
\!-\!24\underline{-\!13\!-\!23} \!-\!34\!-\!14\!-\!24\!-\!34$$
$$\simeq
\!-\!13 \underline{-\!12\!-\!14}  \!-\!23\underline{-\!12\!-\!24}
\!-\!23\!-\!13\!-\!34 \underline{-\!14\!-\!24}  \!-\!34$$
$$\simeq
\!-\!13\!-\!14 \underline{-\!12\!-\!23} \!-\!24 
\underline{-\!12\!-\!23} 
\!-\!13\!-\!34\!-\!24\!-\!14\!-\!34$$
$$ \simeq \!-\!13 \underline{-\!14\!-\!24} \!-\!13\!-\!34
\underline{-\!24\!-\!14}  \!-\!34
\simeq \!-\!13\!-\!13\!-\!34\!-\!34 \simeq \emptyset.$$
Here is a contracting homotopy for the desingularization of 
$a^5=aaaaa$:
$$\underline{-\! 12 \!-\!13} \!-\!14  \!-\!15 \underline{-\! 12 
\!-\!23}
\!-\!24  \!-\!25
\underline{-\! 13 \!-\!23} \!-\!34  \!-\!35 \!-\! 14 \!-\!24  \!-\!34  
\!-\!45
\!-\! 15 \!-\!25  \!-\!35  \!-\!45$$
$$\simeq \!-\!13 \underline{-\!12\!-\!14}  \!-\!15
\!-\!23\underline{-\!12\!-\!24} \!-\!25
\!-\!23\!-\!13\!-\!34 \!-\!35 \underline{-\!14\!-\!24} \!-\! 34\!-\! 45 
\!-\!
15 \!-\!25  \!-\!35  \!-\!45$$
$$\simeq \!-\!13\! -\!14 \underline{-\!12\!  \!-\!15}
 \!-\!23 \!-\!24 \underline{-\!12\!-\!25}
\!-\!23\!-\!13\!-\!34 \!-\!35 \!-\!24 \!-\!14  \!-\! 34\!-\! 45 
\underline{-\!
15 \!-\!25}  \!-\!35  \!-\!45$$
$$\simeq \!-\!13\! -\!14  \!-\!15 \underline{-\!12\!
 \!-\!23}  \!-\!24 \!-\!25 \underline{-\!12
\!-\!23}\!-\!13\!-\!34 \!-\!35 \!-\!24 \!-\!14  \!-\! 34\!-\! 45  
\!-\!25 \!-\!
15  \!-\!35  \!-\!45$$
$$\simeq \underline{-\!13\! -\!14}  \!-\!15    \!-\!24 \!-\!25 
\underline{-\!13\!-\!34} \!-\!35 \!-\!24 \underline{-\!14  \!-\! 
34}\!-\! 45 
\!-\!25 \!-\! 15  \!-\!35  \!-\!45$$
$$\simeq \! -\!14 \underline{ -\!13  \!-\!15 }   \!-\!24 \!-\!25 
\!-\!34
\underline{-\!13  \!-\!35} \!-\!24 \!-\!34  \!-\! 14\!-\! 45  \!-\!25
\underline{-\! 15  \!-\!35}  \!-\!45$$
$$\simeq \! -\!14 \!-\!15 \underline{-\!13     \!-\!24} \!-\!25 \!-\!34 
\!-\!35 \underline{-\!13     \!-\!24} \!-\!34  \!-\! 14\!-\! 45  
\!-\!25 
\!-\!35 \!-\! 15   \!-\!45$$
$$\simeq \underline{-\!14 \!-\!15}   \!-\!25 \!-\!34  \!-\!35   \!-\!34 
\underline{-\! 14\!-\! 45}  \!-\!25  \!-\!35 \underline{-\! 15   
\!-\!45}$$
$$\simeq   \!-\!15 \underline{-\!14   \!-\!25} \!-\!34  \!-\!35   
\!-\!34 
\!-\! 45  \underline{-\!14   \!-\!25}  \!-\!35 \!-\!45 \!-\! 15   $$
$$\simeq   \!-\!15   \underline{-\!34  \!-\!35} \,  \underline{-\!34  
\!-\! 45}
 \,   \underline{-\!35 \!-\!45 } \!-\! 15    $$
$$\simeq   \underline{-\!15    \!-\!35}  \!-\!34    \!-\! 45 \!-\!34    
\!-\!45
 \underline{-\!35    \!-\! 15}    \simeq     \!-\!34    \!-\! 45 
\!-\!34   
\!-\!45  \simeq  \emptyset.  $$
Uisng a topological argument one can show that     $a^m$ is
contractible  for all $m\geq 2$. (One can realize the
desingularization of $a^m$ as the Gauss word of a generic closed
curve in $\mathbf R^2$ and use that all curves in $\mathbf R^2$
are contractible).

4. Let $a\in \alpha$ be a letter such that $\tau(a)\neq a$. Set 
$b=\tau(a)$.
Then the word $w=ababa$  is
contractible. Indeed, the
desingularization yields the nanoword
$$(\A^d=\{a_{1,2},  a_{1, 3}, a_{ 2,3}, b_{1,2}\}, w^d=a_{1,2}  a_{1, 
3}
b_{1,2} a_{1,2}  a_{2, 3} b_{1,2} a_{1,3}  a_{2, 3})$$
where $\vert a_{i,j}\vert =a$ for all $i,j$ and  $\vert b_{1,2}\vert 
=b$. We
have
$$w^d=\underline{a_{1,2} a_{1, 3}} \, \underline{b_{1,2} a_{1,2}}  
a_{2, 3}
\underline{b_{1,2} a_{1,3}}  a_{2, 3} \simeq
 a_{1,3} \underline{a_{1, 2} a_{1,2}}  b_{1,2}   a_{2, 3}  {a_{1,3} 
b_{1,2}} 
a_{2, 3}$$
$$\simeq  \underline{ a_{1,3}    b_{1,2}}   a_{2, 3}  
\underline{a_{1,3}
b_{1,2}}  a_{2, 3} \simeq    a_{2, 3}     a_{2, 3}\simeq \emptyset.$$

5. Let $a\in \alpha$ be a letter such that $\tau(a)\neq  a$. Set 
$b=\tau(a)$.
Then the  words  $ aabab, babaa$, and $baaab$ are pairwise homotopic  
and
homotopic to the nanoword $w=AA'AA'$ with $\vert A\vert =\vert A' \vert 
=a$:
$$(aabab)^d= a_{1,2}  a_{1,3} a_{1,2}   \underline{ a_{2,3} b_{1,2}}  
a_{1,3} 
\underline{ a_{2,3}  b_{1,2}}
 \simeq a_{1,2}    a_{1,3} a_{1,2} a_{1,3}  \approx w,$$
$$
(babaa)^d= \underline{ b_{1,2} a_{1,2}}  a_{1,3}   \underline{ b_{1,2} 
a_{1,2}}
 a_{2,3} a_{1,3} a_{2,3}\simeq  a_{1,3}    a_{2,3} a_{1,3} a_{2,3} 
\approx w
.$$
$$( baaab)^d= b_{1,2}  \underline {a_{1,2}  a_{1,3}}\, \underline 
{a_{1,2}     
a_{2,3} } \, \underline {  a_{1,3}   a_{2,3}  } b_{1,2}
 \simeq  \underline {b_{1,2}  a_{1,3}}  a_{1,2}  a_{2,3}   { a_{1,2} }  
a_{2,3}
  \underline {  a_{1,3}  b_{1,2}}
\simeq
 w.$$
The word
  $aabab$ is not symmetric but  is homotopically symmetric since 
$(aabab)^-=babaa\simeq aabab$.

\subsection{Remarks}\label{remdfpvhhhg}  1.    In analogy  with the 
theory of
framed knots in topology, one can define a {\it framed} homotopy move  
on
nanowords  $(\A, xAAyBBz)\mapsto (\A, xyz)$ which
applies   when
$\vert  A\vert =\tau (\vert B\vert)$.  This move together with homotopy 
moves
of types 2 and 3 generates  an equivalence relation on nanowords called 
{\it
framed homotopy}. It is stronger than homotopy. We shall not study 
framed
homotopy here.

2. One can generalize    the second homotopy move
$(\A, xA B y BAz) \mapsto (\A-\{A, B\},   xyz)    $ allowing it
whenever the pair $(\vert A\vert ,\vert B\vert )$ lies in a fixed
subset of $\alpha^2$. We study here only    the case where this
subset   is the graph of   $\tau$.

3. A  {\it homotopy automorphism} of an \'etale word $w$ over $\alpha$ is a
$\tau$-equivariant permutation   $f:\alpha\to \alpha$ such that
$f_{\#}  (w)\simeq w$. The classification theorems below allow to
compute the group of homotopy automorphisms for   words of length $\leq 5$
and   nanowords of length $\leq 6$.

4. It is possible that all \'etale words over  a 1-element alphabet
are   contractible; at least I do not know obstructions to this.

\section{Homomorphism $\gamma$ }\label{ccc56822}

We introduce an elementary  homotopy invariant of nanowords  over
$\alpha$. The idea is to   associate with a nanoword an element of  a group generated by elements of $\alpha$.

\subsection{Group $\Pi$ and homomorphism $\gamma$}\label{groupga} Recall that an {\it orbit} of  the involution 
$\tau:\alpha\to
\alpha
$   is a   subset    of 
$ 
\alpha$ consisting  either of one element  preserved by $\tau$ or of two elements permuted by $\tau$; in the latter  case the orbit is   {\it free}.   Let $\Pi$ be the group with generators $\{z_a\}_{a\in \alpha}$ and defining relations $z_a z_{\tau
(a)}=1$ for   $a\in \alpha$.       It is clear that $\Pi$ is a free product of   $k$ infinite cyclic
groups and   $l$ cyclic
groups of  order 2 where $k$ is the number of 
free  orbits of $\tau$ and $l$ is the number of fixed points of $\tau$.  The formula $z_a \mapsto z_{\tau(a)}$ defines an involutive automorphism of  $\Pi$ denoted $\tau_\ast$.   

  Consider a nanoword $ (\A, w:\hat n \to \A)$ over $\alpha$. For $i=1,...,n$,  set  $\gamma_i 
  =z_{\vert w(i) \vert}  $ if  $i$ numerates  the first entry of  $w(i)$ in   $w$, that  is  if $w(i)\neq w(j)$ for $j<i$. Otherwise, set  $\gamma_i =z_{\tau(\vert w(i) \vert)}
  =(z_{\vert w(i) \vert})^{-1}   $. Set $\gamma (w)=\gamma_1 \cdots  \gamma_n \in \Pi$.  For example, if $w= ABAB $ with $\vert A\vert =a\in \alpha, \vert B\vert =b\in \alpha$, then $\gamma(w)= z_a z_b z_a^{-1} z_b^{-1} $.   It is
easy to check that $
\gamma (w)
$    is invariant under homotopy moves on $w$. This provides an efficient and  easily computable homotopy invariant of  nanowords. The mapping $w\mapsto \gamma (w)
$  defines a monoid  homomorphism $\gamma:\N_\bullet(\alpha)\to \Pi$ such that $\gamma (\overline w)=\tau_\ast(\gamma (w))$ and $\gamma (w^-)=(\gamma (w))^{-1}$ for any   $w$.

\begin{lemma}\label{25edcolbrad}   $\gamma(\N_\bullet(\alpha))= [\Pi, \Pi]$. \end{lemma}
    \begin{proof}     It is obvious that the composition of $\gamma$ with the projection $\Pi \to  \Pi/ [\Pi, \Pi]$ is a trivial homomorphism. Hence    $\gamma(\N_\bullet(\alpha)) \subset    [\Pi, \Pi]$.  To prove the opposite inclusion, pick a set $\alpha_0\subset \alpha$ meeting each orbit of $\tau$ in one element. The group $\Pi$ can be presented by the generators $\{z_a\}_{a\in \alpha_0}$ and defining relations $z_a z_a=1$ for all $a\in \alpha_0$ such that $\tau(a)=a$.   Each element of $ [\Pi, \Pi]$ can be 
presented by a word
$W$ in the alphabet  $\{z_a^{\pm 1}\}_{a\in \alpha_0}$ in which every
$z_a$   appears with total power $0$. Partition   the entries of    $z_a^{\pm 1}$ in $W$ into pairs  $(z_a, z_a^{-1})$ or $(z_a^{-1}, z_a)$ in an arbitrary way and let $s_a$ be the resulting set of  
pairs. For a pair $A\in s_a$, set $\vert A\vert =a$ if the first entry of this pair   is  $z_a$
and set   $\vert A\vert =\tau(a) $ if the first entry of this pair   is  $z_a^{-1}$. This makes  the set 
$\A=\amalg_a s_a$ into an 
$\alpha$-alphabet. Replacing in $ W$ each entry $z_a^{\pm 1}$ with the only   $A\in
s_a$  containing this entry  we obtain a nanoword $(\A,   w)$ such that $\gamma ( w)=W$. For example if $W=z_a^2 z_b^{-1} z_a^{-1} z_b z_a^{-1} $, then  $(\A=\{A_1, A_2, B\}, w=A_1A_2 BA_1 B A_2  )$ with $\vert A_1\vert =\vert A_2\vert =a, \vert B\vert =\tau(b)$ is  a nanoword such that $\gamma(w)=W$. 
\end{proof}

The  commutator group  $[\Pi, \Pi]$ is a free group   for any $\tau$. Its rank is infinite if (in the notation above) $k\geq 2$  or $k=1$ and $l\geq 1$.  If $k=1$ and $l=0$, then   $[\Pi, \Pi]=1$. If $k=0$, then    $[\Pi, \Pi]$  is a free group of rank $2^{l-1} (l-2)+1$.  One way to see  it goes by realizing $\Pi$ as the
fundamental group of the connected sum $X$ of $l$ copies of $RP^3$ and observing that the maximal abelian covering  of $X$ has the same  fundamental group as  a
connected graph with $2^{l-1} l$ vertices and 
$2^l (l-1)$ edges. 

If $\tau$ has
at least two orbits, then     Lemma
\ref{25edcolbrad} implies that $\N_\bullet(\alpha) $ is an infinite monoid.
If $\tau$ has at least three orbits, then     Lemma \ref{25edcolbrad} implies that  $\N_\bullet(\alpha)$ is  
non-abelian.

\subsection{Homomorphism $\gamma'$}\label{applega}   It is  sometimes  sufficient to consider  a weaker invariant $\gamma':\N_\bullet(\alpha) \to \Pi'$ where $\Pi'$ is obtained from $\Pi$  by  adding  relations $z_a z_a=1$ for all $a$ and $\gamma'$ is the composition of $\gamma$ with the projection $\Pi\to \Pi'$.    To give an  example,  pick an integer
$m\geq 1$ and elements
$a,b\in
\alpha$ lying in different orbits of
$\tau$.  Then $z_az_b \in \Pi' $ is an element of infinite order  and $(z_az_b)^{-1}=z_bz_a$ in $\Pi'$. Consider the nanoword
$w_m=A_1B_1A_2B_2\cdots A_mB_mA_1B_1A_2B_2\cdots A_mB_m$ with  $m\geq 1$ and 
$\vert A_i\vert=a, \vert B_i\vert =b$ for all $i$. It is clear that $\gamma' (w_m)= (z_az_b)^{2m}   $. Therefore the nanowords $\{w_m\}_m$ are non-contractible and
mutually non-homotopic. The opposite nanowords $w_m^-=B_mA_m\cdots  B_1A_1 B_mA_m\cdots  B_1A_1 $ are also mutually non-homotopic
and homotopically distinct  from $\{w_m\}_m$. Indeed $  \gamma' (w_m^-)=(\gamma' (w_m))^{-1} =(z_b z_a)^{2m}= (z_az_b)^{-2m} $ for any $m \geq 1$.  

  We   use $\gamma'$ to associate with any nanoword $w$ a   function $\mu_w:\hat \alpha
\times
\hat \alpha \to \ZZ$ where $\hat \alpha= \alpha/\tau$ is the set of orbits of $\tau$.   The orbit $\{a,\tau(a)\}$ of $a\in \alpha$ will be denoted $\hat a$. Pick
$a,b\in
\alpha$. If $\hat a=\hat b$, then   $\mu_w(\hat a, \hat b)=0$. Suppose that $a$ and $b$   lie in different orbits.  Consider the group  $\Pi_0=(x, y\,\vert
\, x^2=y^2=1)$ and 
 the group homomorphism
$\rho :\Pi'
\to
\Pi_0$ such that  $\rho(z_a)=\rho( z_{\tau (a)})=x$,  $\rho(z_b)=\rho( z_{\tau (b)})=y$, and $ \rho(z_c)=1$ for all other $c\in \alpha$. 
Observe that   $[\Pi_0, \Pi_0]$ is an infinite cyclic group generated by $xyxy= (yxyx)^{-1}$. Since $\rho \gamma' (w)\in [\Pi_0, \Pi_0]$, there is a
unique 
$m\in \ZZ$ such that $\rho \gamma' (w) =(xyxy)^m$. Set $\mu_w(\hat a, \hat b)=m$.  It follows from the definitions that the mapping  $ \mu_w: \hat \alpha
\times
\hat \alpha \to
\ZZ$ is skew-symmetric in the sense that  $\mu(\hat a, \hat b)= -\mu(\hat b, \hat a)$ for all $a,b\in \alpha$.   Since $\gamma'(w)$ is a homotopy invariant of $w$,
so is
$\mu_w $.  Note   that
$\mu_{w_1w_2}=\mu_{w_1}+
\mu_{w_2}$ for any nanowords
$w_1,w_2$.  For  $w=ABAB$ with $ \vert A\vert=a\in \alpha ,  \vert B\vert=b \in \alpha$ lying in different orbits of $\tau$, the function  $\mu_w$ takes values 1
and
$-1$ on the pairs
$\hat a, \hat  b$ and $ \hat b, \hat a  $, respectively, and the value 0 on all other pairs.  This implies that every skew-symmetric mapping 
$\hat \alpha \times \hat \alpha \to \ZZ$ is the $\mu$-function of a   nanoword over $\alpha$.

\subsection{Homomorphism $\tilde \gamma$}\label{refinga} A  more careful approach yields the following  refined version of $\gamma$.   
Let $\tilde \Pi$ be the group with generators $\{\tilde z_a\}_{a\in \alpha}$ and defining relations 
  $ \tilde z_a \tilde z_{\tau(a)} \tilde z_b=\tilde z_b \tilde z_a \tilde z_{\tau
(a)} $ for   all $a,b\in \alpha$.       The formula $\tilde z_a\mapsto  z_a$ defines a  projection $\tilde \Pi \to \Pi$.
Replacing $z_a$ with $\tilde z_a$ in the definition of $\gamma$, we obtain  a  lift of $\gamma$ to a monoid  homomorphism $\tilde \gamma:\N_\bullet(\alpha)\to
\tilde
\Pi$.   Clearly, $\tilde \gamma(\N_\bullet(\alpha)) \subset    [\tilde \Pi, \tilde \Pi]$.

  Despite  many nice features, the homomorphisms $\gamma$, $\gamma'$, and $\tilde \gamma$ often fail  to distinguish homotopy classes of nanowords.   For   example, as we shall see
in Sect.\ \ref{nanexacacacacota1},  the nanoword
$  w= ABAB $  with $\vert A\vert=\vert B\vert  \neq \tau (\vert A\vert) $ 
is non-contractible    but obviously    $\tilde \gamma (w)=1$.  
  
 \section{Coverings of nanowords}\label{orticcdfd}

We define  for each nanoword   a family of nanowords  called its coverings.   

  \subsection{Group $\pi$ and interlacement of letters}\label{covfunct}      Let
$\pi=\pi(\alpha, \tau)$ be the multiplicative abelian group with  generators
$\{a\}_{a\in
\alpha}$ and defining relations  $ a \, \tau (a)=1 $ for all ${a\in
\alpha}$. Clearly,      $\pi=\Pi/ [\Pi, \Pi]$ where $\Pi$ is  the group
  considered in Sect.\ \ref{groupga}.
The group  
$\pi$ is a direct product  of
 copies of $\ZZ/2\ZZ$ numerated by $\{ a\in \alpha \,\vert\, \tau (a)=a\}$ and copies of $\ZZ$ numerated by free orbits of $\tau$.   
 
Consider a nanoword  $(\A , w)$ over   $ \alpha $. We  say that two letters
$A, B\in \A$ are {\it   $w$-interlaced} if   
  $$w=\cdots A \cdots B\cdots A 
\cdots B  \cdots  \,\,\,\,\,\, {\rm or}Ê\,\,\,\,\,\,w=\cdots  B\cdots A 
\cdots B 
\cdots A  \cdots.$$ Set in the first case  $n_w(A,B)= 1\in \ZZ$ and in the second case $n_w (A,B)=-1\in \ZZ$. For all other pairs of letters $A,B\in \A$, we set
$n_w(A,B)=0\in \ZZ$ and say that $A,B$ are  {\it   $w$-unlaced}.      Note that   $n_w(A,B)=-n_w(B,A)$ and   $n_w(A,A)=0$.  For  any
$A\in \A$, set  
$$ [A]_w=
\prod_{ B \in \A }   \vert B\vert^{n_w (A,B)} \in \pi . $$

\subsection{Coverings}\label{invaricovfunct}    By an {\it $\alpha$-family of subgroups} of $\pi$, we  mean a family $H=\{H_a\subset \pi\}_{a\in \alpha}$  where
$H_a$ is a subgroup of $\pi$ such that $H_a=H_{\tau(a)}$ for all $a$.  For a   nanoword $(\A,w)$  over $\alpha$,      consider  the  nanoword  $(\A^H,w^H)$ over
$\alpha$  obtained from 
$(\A,w)$ by deleting   from both $\A$ and $w$ all letters $A$ such that $[A]_w \notin H_{\vert A\vert}$. It is called the {\it $H$-covering} of $(\A,w)$.

\begin{lemma}\label{covl} If  two nanowords  are homotopic, then their $H$-coverings
 are homotopic.
\end{lemma} 
                     \begin{proof}  It is clear that $H$-coverings of isomorphic nanowords are isomorphic.  It  remains to prove that if a nanoword $v$ is
obtained from a nanoword 
$w$ by a homotopy move then  $v^H$   is obtained from   $w^H$ by a homotopy move.  Consider the   move
$w=(\A,  x A A  y) 
\mapsto  (\A- \{A\},  x  y)=v$.  Then $n_w(A,B)=n_w(B,A)=0$ for all $B\in \A$.  Hence   $[B]_w=[B]_v$  for any  $B\in \A-\{A\}$
and
$[A]_w=1\in \pi$.    Therefore  $A$ survives in $w^H$ and  the nanoword $v^H$   is obtained from   $w^H$ by deleting $A$, i.e., by the first   homotopy move.  

Consider the   move $w=(\A,  xA 
B  y BAz) \mapsto (\A-\{A, B\},   xyz)  =v  $ where  $\vert A\vert =a\in \pi$ and $\vert B\vert= \tau (a) $.   It follows from the definitions that $n_w(A,B)=0$
and   
$n_w(A,C)= n_w(B,C)$ for all $C\in   \A-\{A, B\}$.  Therefore $[A]_w=[B]_w $.  Observe  also  that
for    
$  C\in   \A-\{A, B\}$,    $$[C]_w \,[C]_v^{-1}= \vert A\vert^{n_w(C,A) }
\vert B\vert^{n_w(C,B)}=  (a\,  \tau (a))^{n_w(C,A) }=1.$$ 
Hence    $  [C]_w= [C]_v $. 
 Therefore in the case   $[A]_w  \notin H_a$, we have $w^H=v^H$ and in the case  $[A]_w  \in H_a$, the nanoword $v^H$   is obtained from  
$w^H$ by the second homotopy move.  

Consider the move $w=(\A, xAByACzBCt)\mapsto  (\A, xBAyCAzCBt)=v$. Set  $a=\vert A\vert =\vert B\vert = \vert C\vert \in \alpha$. 
It is clear that for   $D\in \A-\{A,B,C\}$ and $E\in \A$, we have  $n_w (D,E)=n_v(D,E)$.   Thus 
$  [D]_w=[D]_v 
$.   Therefore 
  $D $   either  disappears in both  $w^H, v^H$ or survives in both $w^H, v^H$. 
  For     $D\in \{A,B,C\}$,  the
sum 
$n_w(D,A)  +n_w(D,B)  +n_w(D,C)  $
 is preserved  under the move.  The product
$$ \vert A\vert^{n_w(D,A)}   \vert B\vert^{n_w(D,B)} \vert
C\vert^{n_w(D,C)}=a^{n_w(D,A)+{n_w(D,B)} + n_w(D,C)}$$ is  also preserved.  This implies that $[D]_w=[D]_v$. 
We claim that $[A]_w[C]_w=[B]_w$. To see this,   set 
$x\ast  y=\prod_E [E]_w\in \pi$ where $E$ runs over the letters which appear both in $x$ and   in $y$.  The expressions $y\ast   z, y\ast  t$, etc.
are defined similarly.  It follows from the definitions  that
$$[A]_w= a^{n_w(A,A)+ n_w(A,B) +n_w(A,C)} (y\ast  z) (y\ast  t) (x\ast  y)^{-1}=a (y\ast  z) (y\ast  t) (x\ast  y)^{-1} ,$$
$$[B]_w= a^{n_w(B,A)+ n_w(B,B) +n_w(B,C)} (y\ast  t)(z\ast  t) (x\ast  y)^{-1}(x\ast  z)^{-1}$$
$$=(y\ast  t)(z\ast  t) (x\ast  y)^{-1}(x\ast  z)^{-1},$$
$$[C]_w= a^{n_w(C,A)+ n_w(C,B) +n_w(C,C)}  (z\ast  t) (x\ast  z)^{-1} (y\ast  z)^{-1}= a^{-1}  (z\ast  t) (x\ast  z)^{-1} (y\ast  z)^{-1}.$$
These computations imply that $[A]_w[C]_w=[B]_w$.  There are three possibilities:   the elements $[A]_w, [B]_w, [C]_w\in \pi$ belong to
$H_a$;  exactly one of them belongs to
$H_a$;   neither  of them belongs to
$H_a$. In the first case     $v^H$   is obtained from   $w^H$ by
the third homotopy move. In the second and third cases   $v^H=w^H$. 
 \end{proof}
 
 Note that  $(w^-)^H=(w^H)^-$, ${\overline w}^H=
  \overline {w^H}$, $(ww')^H=w^H (w')^H$ for any nanowords $w,w'$. 
   
  Any $\tau$-invariant subgroup $H$ of $\pi$ determines an    $\alpha$-family of subgroups  of $\pi$   by $H_a  =H$ for all $a\in \alpha $. This  family is  denoted by the same symbol $H$.
   
The following example shows that   homotopy invariants  of the coverings can be efficiently used to 
distinguish homotopy classes of nanowords.

\subsection{Example}\label{nezzexixixremar}      Set   $\alpha=\{a,b\}$ with   $\tau=\id$.   Consider the nanoword 
$w=A_1B_1B_2A_2A_1A_3B_1A_3A_2B_2$  with  
$\vert A_i\vert=a$ and    
$\vert B_j\vert =b$ for all $i,j$.  We have    $\pi= \pi(\alpha, \tau)=(a,b\,\vert \, a^2=b^2=1, ab
=ba)$,  $[A_1]_w=a\in \pi,
[A_3]_w=b\in \pi$ and
$[A_2]_w=[B_1]_w=[B_2]_w=ab\in \pi$.  Let $H  =\{1,ab\}\subset  \pi$.  Then 
$w^H=B_1B_2A_2 B_1 A_2B_2$. Clearly, $\gamma (w^H)=z_az_b z_az_b\neq 1$ and $ \mu_{w^H} (a,b)=1$. Thus $w$ is non-contractible.   All the invariants of $w$ defined in Sect.\ \ref{ccc56822} are trivial since  $\tilde \gamma(w)=1 $.

\section{Self-linking}\label{wxcv}

Starting  from the notion of linking (or interlacement) of letters in  a word, we construct   a family of polynomial   
invariants of nanowords.

 \subsection{Self-linking  class}\label{27covfunct}    Let $\pi=\pi(\alpha, \tau)$ be the group defined in Sect.\ \ref{covfunct}. Recall that given a nanoword  $(\A , w)$
over  
$
\alpha $, we   associate with every $A\in \A$ an element $ [A]_w$ of $\pi$, see Sect.\ \ref{covfunct}.  
Set $\A_w=\{A\in \A\,\vert \, [A]_w\neq 1\}$. For    a letter  $a\in \alpha$,  define its {\it self-linking class} by
$$[a]_w= \sum_{ A \in \A_w ,  \vert A\vert =a} [A]_w
\in
\ZZ\pi$$
 where  
$\ZZ\pi$ is the group ring of $\pi$ with integer coefficients.
 By definition, the neutral element $1\in \pi$ appears in this expression with coefficient $0$. 
The key property of the self-linking class is provided by the following theorem.

   \begin{theor}\label{dcppmolbrad}   If $\tau(a)=a$, then $[a]_w (\modu 2) \in (\ZZ/ 2\ZZ) \pi $ is a homotopy invariant    of $w$. 
 If $\tau(a)\neq a$, then $[a]_w-[\tau(a)]_w\in \ZZ\pi$ is a homotopy invariant    of $w$.\end{theor}
    \begin{proof}     
It is obvious that isomorphic nanowords give rise to the same  self-linking class  for all   $a\in \alpha$.    Consider the first  homotopy move
$w=(\A,  x A A  y) 
\mapsto  (\A- \{A\},  x  y)=v$.  As  was shown in the proof of Lemma \ref{covl},   $[B]_w=[B]_v$  for any  $B\in \A-\{A\}$
and
$[A]_w=1$.   Therefore  $[a]_w=[a]_v$   for  any $a\in  \alpha$.

Consider the      move $w=(\A,  xA 
B  y BAz) \mapsto (\A-\{A, B\},   xyz)  =v  $ where  $\vert A\vert =\tau (\vert B\vert) $.  
By the proof of  Lemma \ref{covl},  $[A]_w=[B]_w $ and   
$  [C]_w= [C]_v $ for    
$  C\in   \A-\{A, B\}$.
  This implies that if $a\in \alpha-\{  \vert A\vert, \vert B\vert\}$, then  $[a]_w=[a]_v$.
If $a=\vert A\vert=\vert B\vert $, then $A,B$ contribute
$[A]_w+ [B]_w=2[A]_w$ to $[a]_w$ so that $[a]_w=[a]_v (\modu 2 )$.
If $a=\vert A\vert  \neq \tau(a)=\vert B\vert $ then $A,B$ contribute
$[A]_w- [B]_w=0$ to $[a]_w-[\tau(a)]_w$ so that $[a]_w -[\tau(a)]_w=[a]_v -[\tau(a)]_v$.

Consider the     move $w=(\A, xAByACzBCt)\mapsto  (\A, xBAyCAzCBt)=v$.  
By the proof of  Lemma \ref{covl}, $[D]_w=[D]_v$  for any  $D\in
\A$.  Hence   $[a]_w=[a]_v$  for all  $a\in \alpha$.    \end{proof}

\subsection{Self-linking function}\label{vpamxct} Theorem \ref{dcppmolbrad} suggests   the following definitions. For $a\in \alpha$, set $R_a=
\ZZ/2\ZZ$ if $\tau (a)=a$ and 
$\ZZ_a= \ZZ $ if $\tau (a)\neq a$.  Consider  the system of commutative rings $\{R_a\pi \}_{a\in \alpha}$.  
By a {\it section} of this system we mean a mapping which assigns to every $a\in \alpha$ an element of $R_a\pi $. Each 
nanoword
$w$ determines such a section $u^w$    by 
 $ u^w(a)= 
[a]_w (\modu 2) \in R_a\pi$ if 
$\tau (a)=a$ and $
  u^w(a)=  [a]_w-[\tau(a)]_w \in R_a\pi$ if $  \tau (a)\neq a$.   By Theorem
\ref{dcppmolbrad},
$u^w $, also denoted $u(w)$, is a homotopy invariant of
$w$. It is called the {\it self-linking function}   of $w$. If $w$ is contractible, then $u^w=0$.

It is clear that $u^w(\tau (a))=-u^w(a)$ for all  $ a$. Therefore the function $u^w$ is entirely determined by its restriction to any 
set 
$\alpha_0\subset
\alpha
$ meeting every  orbit  of $\tau$ in exactly one element. We   call such $\alpha_0$ an {\it orientation} of $\tau$.
Given  
$\alpha_0$ and $a\in \alpha_0$, we can   view $u^w(a)$ as a Laurent polynomial as follows.   Any  element of  $ \pi$     expands
as $ \prod_{b\in \alpha_0} b^{m_b}   $
with   $m_b\in \ZZ$.  The integer $m_b$ is  uniquely determined   if $\tau(b)\neq b$ and is  determined   modulo 2 if $\tau(b)= b$. The elements
of 
$R_a\pi$ (and in particular $u^w(a)$)  are just   Laurent polynomials with   coefficients in $R_a$ in the commuting variables  
${b\in \alpha_0}$ subject to $b^2=1$ for  $\tau(b)=b$.

  For  a nanoword $w$ obtained as  the 
product of nanowords   
$w_1$  and   
$w_2$,  we have  
   $u^w(a)=u^{w_1} (a)+ u^{w_2}(a)$ for all  $a\in \alpha$.
  This  follows from the definitions and the fact that   letters of $w_1$ and $w_2$ are  unlaced in $w$.
Thus the formula  $w\to u^w$ defines a monoid homomorphism from $\N_\bullet (\alpha)$ to the additive group of
sections of the system   $\{R_a\pi \}_{a\in \alpha}$.  We describe the image  of this homomorphism in the next subsection.

To compute the self-linking function of the opposite nanoword $w^-$, observe that   $[a]_{w^-} =\overline {[a]_w}$  for any $a\in \alpha$ where the   overbar
 denotes the   involution in  $R_a \pi $ sending   elements of $\pi$ to their  inverses.  
Therefore   
$u({w^-})=\overline {u(w)}$. For the
inverse nanoword $\overline w$,  we have  the identity $[a]_{\overline w} =\overline {[\tau (a)]_w}$.
Therefore      
$u(
\overline
 w ) =-  \overline {u(w)}$ and $u(
\overline
 w^- ) =-    {u(w)}$.  In other words $u(w\,\overline
 w^-)=0$ for all $w$.    If $w$ is homotopically skew-symmetric, then $u^w(a)=0$ for all $a\in  \alpha$ such that $\tau(a)\neq a$.

The self-linking function can be used to estimate $\vert  \vert w\vert \vert $ from below. If $w$ is homotopic to a nanoword of length $2m$, then for any $a\in \alpha$, the element $u^w(a)\in R_a\pi$ is an algebraic sum of monomials in  the  generators $\{b\}_{b\in \alpha}\subset \pi$
of degree $\leq m-1$.

\subsection{Characterization}\label{readfryopn1} Consider the following characterization problem:  given an orientation 
$\alpha_0\subset
\alpha
$, when a  system $u=\{u(a)\in R_a\pi
\}_{a\in \alpha_0}$   is the  self-linking function   of a  nanoword   ?   We call such a system a {\it section} of  $\{ R_a\pi
\}_{a\in \alpha_0}$.  
Denote by $\delta_a$ the additive homomorphism $R_a\pi \to R_a$ sending $1\in \pi$ to $1\in R_a$ and sending all other elements of $\pi$ to $0$.        For   
${a,b\in \alpha_0}$, set $R_{a,b}=\ZZ$ if $\tau(a)\neq a$ and $\tau (b)\neq b$ and
$R_{a,b}=\ZZ/2\ZZ$  otherwise. (Clearly $R_{a,b}=R_a\otimes R_b$ and $R_{a,a}=R_a$). 
Let $ \partial_b$  be the additive
homomorphism $R_a\pi\to R_{a,b}$ sending a monomial $\prod_{b\in \alpha_0} b^{m_b}$   to $m_b$. We say that  a section  $u$ of  $\{ R_a\pi
\}_{a\in \alpha_0}$  is {\it
skew-symmetric} if 
 $\delta_a (u(a))= \partial_a (u (a )) =0 $ for any $a \in \alpha_0$ and  $
\partial_a (u (b))+ \partial_b (u (a))=0 $ for any $a,b\in \alpha_0$.

 \begin{theor}\label{bdbdbdad}  A  section    of $\{ R_a\pi
\}_{a\in \alpha_0}$  is the  self-linking function   of a   nanoword    if and only if it is   skew-symmetric.\end{theor}
    \begin{proof}   Pick a nanoword $w$ and set $u=u^w$.  The equality $\delta_a (u(a))=0$ follows    from the definitions. Let us
prove the other identities.  For  
$a,b\in
\alpha_0$, set 
$$a\circ b= \sum_{ A , B \in \A ,  \vert A\vert =a  ,  \vert B\vert =b}  n_w (A,B) \in R_{a,b}.$$
The identities $n_w(A,B)+n_w(B,A)=0$ and  $n_w(A,A)=0$ imply that $a\circ b +b \circ a=a\circ a=0$.  

 Suppose first that $\tau (a)=a, \tau (b)=b$.  Observe that for any $A\in
\A-\A_w$ we have $\partial_b ([A]_w)= 0$. Therefore  
$$\partial_b (u(a))=\partial_b ([a]_w )=\sum_{ A \in \A_w ,  \vert A\vert =a} \partial_b ( [A]_w)
=\sum_{ A \in \A ,  \vert A\vert =a}  \partial_b ( [A]_w)$$
$$
=\sum_{ A \in \A ,  \vert A\vert =a}  \partial_b (\prod_{ B \in \A }   \vert B\vert^{n_w (A,B)})
=\sum_{ A , B \in \A ,  \vert A\vert =a  ,  \vert B\vert =b}  n_w (A,B) =a\circ b.$$
This implies that $ \partial_a (u (a )) = 0=
\partial_a (u (b))+ \partial_b (u (a)) $.

 Suppose  that $\tau (a)= a, \tau (b) \neq b$.    Then as above
$$\partial_b (u(a))=\partial_b ([a]_w  )
=\sum_{ A \in \A ,  \vert A\vert =a}  \partial_b (\prod_{ B \in \A }   \vert B\vert^{n_w (A,B)})$$
$$
=  \sum_{ A,B \in \A ,  \vert A\vert =a, \vert B\vert =b}   {n_w (A,B)}  -\sum_{ A,B \in \A ,  \vert A\vert =a, \vert B\vert =\tau (b)}   {n_w (A,B)} =
a\circ b -a 
\circ \tau (b), $$  
$$\partial_a (u(b))=\partial_a ([b]_w -[\tau(b)]_w)$$
$$
=\sum_{ B \in \A ,  \vert B\vert =b}  \partial_a ( [B]_w) -\sum_{ B \in \A ,  \vert B\vert =\tau (b)}  \partial_a ( [B]_w) =b\circ a -\tau (b) \circ a.$$
Hence
$$\partial_a (u(b))+ \partial_b (u(a))=b\circ a-\tau (b) \circ a+ a\circ b -a 
\circ \tau (b)=0.$$
 
In the case $\tau (a)\neq a, \tau (b)\neq b$, similar  computations  give 
$$\partial_b (u(a))=a\circ b - a \circ \tau (b) -\tau (a) \circ b+ \tau (a) \circ \tau (b)$$
and the equalities $ \partial_a (u (a )) = 0=
\partial_a (u (b))+ \partial_b (u (a)) $  follow  from the properties of $\circ$.  

We now prove the sufficiency of the conditions stated in the theorem.    For $a\in \alpha$, we  call  an element of
$R_a\pi$ 
 {\it linear} if it is a  linear combination of elements of $  \alpha_0\subset \pi$.  (Elements of $\alpha- \alpha_0\subset \pi$ and $1\in R_a\pi$ are not   linear). 
A section  $u$  of $  \{  R_a\pi
\}_{a\in \alpha_0}$   is {\it linear} if  its values on all
elements of $  \alpha_0$ are linear.  A     section  $u$ of $  \{  R_a\pi
\}_{a\in \alpha_0}$   is {\it realizable} if it is the restriction  to $\alpha_0$ of the  self-linking function   of a  nanoword.  Note that if $u$ is realizable then so is $-u$.
By the additivity of the self-linking function, the sum  of realizable systems is realizable.

The idea of the proof is to take any    skew-symmetric system
$u$ and to deduce from it consecutively  certain \lq\lq model" realizable sections  so that at the end   we obtain $0$. Since $0$ is
realized by an  empty nanoword,  
$u$ is realizable.  
A model section is defined as follows.  Pick
$a\in \alpha_0, b_1,...,b_n \in \alpha 
$ (possibly with repetitions) where $n\geq 1$. Consider  the nanoword $w=  B_1B_2\cdots B_n AB_n\cdots B_2B_1A $ with  
$\vert A\vert =a$ and $ \vert B_i\vert =\tau(b_i) $ for all $i$.  It is clear that $[A]_w=\prod_i b_i $ and $[B_i]_w=a$ for all $i$.  
It follows from the definitions that   $u^w(a)= \prod_i b_i $ modulo linear terms (corresponding to $i$ such that $b_i=a$ or $b_i=\tau(a)$). The values of
$u^w$ on all   elements of
$\alpha_0-\{a\}$ are linear.   Adding or subtracting such model sections $u^w$ from the original section $u$ and using that $\delta_a (u(a))=0$ we can obtain a   
linear skew-symmetric section, again  denoted  $u$. For any $a\in \alpha_0$,  we can uniquely expand  $u(a)=\sum_{b\in
\alpha_0} r(a,b) b$ where $r(a,b)\in R_a$. Pick  
$a\in
\alpha_0$ such that
$\tau (a)=a$.      For   $b\in
\alpha_0-\{a\}$, consider  the nanoword $w =   AB AB $ with  
$\vert A\vert =a$ and $ \vert B\vert =b $.  It is clear that $[A]_w=  b$ and $[B]_w=a^{-1}=a\in \pi$.  Therefore $u^w(a)=b, u^w(b)=a$, and   all other
values of
$u^w$ are $0$.  Adding $\pm u^w$ to $u$, we   change  $r(a,b)$ by $\pm 1$ keeping all  $r(a,b')$ for $b'\neq  b$.   
 Proceeding by induction, we transform $u$ into a linear    skew-symmetric section, still denoted  $u$,  such that $r(a,b)= 0$ for all $b\in \alpha-\{a\}$. 
By skew-symmetry,  $r(b,a)=r(a,b)=0$ and $r (a,a)=\partial_a (a) =0  \in \ZZ/2\ZZ$.  Thus $u(a)=0$ and  $a$ does not show up in the expansions of
$u(b)$ for all $b\in \alpha_0$.  Inductively, we can ensure  that this holds for all $a\in \alpha_0$ such that $\tau (a) =a$.
Pick now   $a\in \alpha_0$ with $\tau(a)\neq a$.  Pick $b\in
\alpha_0-\{a  \}$ with $\tau(b)\neq  b$.  Consider  the nanoword $w =   AB AB $ with  
$\vert A\vert =\tau (a) $ and $ \vert B\vert =b $.  It is clear that  $u^w(a)=-[A]_w=-b, u^w(b)=[B]_w=(\tau (a))^{-1}=a$. 
Adding the linear section $\pm u^w$ to $u$, we   change  $ u(a)$ by $\pm b$ keeping all  other values of $u$ except $u(b)$. Proceeding in this way we can ensure
that
$u(a)= r(a,a) a$. By skew-symmetry, $r (a,a)=\partial_a (a) =0  \in \ZZ$.  Then $u(a)=0$ and again by skew-symmetry, $r(b,a)=0$ for all $b\in
\alpha_0-\{a\}$.   Proceeding further in this way we eventually transform $u$ into the  zero section as required.
 \end{proof}

   \subsection{Examples}\label{unkn1}  1.  Pick $a\in \alpha$ with $\tau(a)\neq a$ and consider the nanoword $w=ABAB$ with $\vert A\vert =\vert
B
\vert =a$. Clearly,  $[A]_w=a, [B]_w=a^{-1}$ and the self-linking function $u^w$ takes the value $a+a^{-1} \in R_a\pi=
\ZZ\pi$ on
$a$, the value $-a-a^{-1} $ on
$\tau(a)$, and the value
$0$ on all other elements of $\alpha$.  Since $u^w\neq 0$, the word $w$ is not contractible.  Since  $u^w\neq -u^w$, the word $w$ is not homotopic
to its inverse.  For   $n\geq 2$, the polynomial  $u^w(a)$ is not divisible by $n$ and therefore  $w$ is not homotopic to the $n$-th power of a nanoword.
Note that
$ \overline {u(w)}=u(w)$ which is compatible with the fact that $w$ is  isomorphic  to its opposite.

2.  Fix    $a,b \in  \alpha$ belonging to  different orbits of $\tau$. For finite sequences  
$\varepsilon=(\varepsilon_i)_{i=1}^m, \mu=(\mu_j)_{j=1}^n $ consisting of zeros and ones  we define a nanoword $$ 
 w=w(\varepsilon, \mu) =A_1 A_2 \cdots A_m B_1 B_2\cdots B_n A_m \cdots A_2 A_1 B_n\cdots B_2 B_1 $$
 over
$\alpha$  with  
$\vert A_i\vert=\tau^{\varepsilon_i} (a)$ and    
$\vert B_j\vert =\tau^{\mu_j} (b)$ for all $i,j$.   We   give  a  homotopy classification  of
these nanowords using the self-linking function.  

Case $\tau (a)=a, \tau (b)=b$. Observe that $u^w(a)=0$ if $mn$ is even    and 
$u^w(a)= b
$  if
$mn$ is odd.  Hence   $mn
(\modu 2)$ is a   homotopy invariant of      $ w=w(\varepsilon, \mu)$.  Applying  the   first two homotopy moves it is easy to see that      $mn
(\modu 2)$ is a full homotopy invariant of      $ w(\varepsilon, \mu)$: any two nanowords of this type with the same   $mn
(\modu 2)$ are homotopic.  
 
Case  $\tau (a)\neq a, \tau (b)=b$.   By the second homotopy move,  if  the length  $n$ of $\mu$  is even 
then    $ w=w(\varepsilon, \mu)  $  is contractible. Suppose that $n$ is odd. 
Set $$k =\card \, \{i=1,...,m\,\vert\,   \varepsilon_i =0\}- \card \, \{i=1,...,m\,\vert\,   \varepsilon_i =1\}  \in \ZZ.$$  A direct computation  gives 
 $u^w(a)=kb$ so that   $k$ is a   homotopy invariant.   It is easy to see that  $k$ is a full  homotopy invariant:  two  nanowords   
$ w(\varepsilon, \mu), w(\varepsilon', \mu')  $ with  $\mu, \mu'$ of odd length  are homotopic if and only if they give rise to the same $k$.  In particular,  such a
nanoword    is contractible iff $k=0$.  

Case  $\tau (a)= a, \tau (b)\neq b$. This case is similar to the previous one.  If  the length  $m$ of $\varepsilon$  is even 
then    $ w(\varepsilon, \mu)  $  is contractible. For odd $m$, the classifying invariant is the integer
$$l =\card \, \{j=1,...,m\,\vert\,   \mu_j =0\}- \card \, \{j=1,...,m\,\vert\,   \mu_j =1\}.$$

Case $\tau (a)\neq a, \tau (b)\neq b$.  Let  $k,l$ be the integers determined by $\varepsilon, \mu$  as above.     We have    $u^w (a)=kb^l$ if $l\neq 0$ and
$u^w (a)=0$ if $l=0$. Similarly, 
$u^w (b) =l a^{-k}$ if $k\neq 0$ and $u^w (b)=0$ if $k=0$.  Therefore   $w(\varepsilon, \mu) $  is contractible iff
$k=0$ or $l=0$.  For the nanowords  $w(\varepsilon, \mu) $ with $kl\neq 0$,  the pair $(k,l)$ is a  full homotopy invariant:  two  such nanowords
   are homotopic if and only if they give rise to the same  $k,l$.

3. The following example  illustrates  the   homotopy invariants derived from the self-linking function via  
coverings.      Let
$\alpha=\{a,b\}$ with   involution
$\tau$ permuting
$a $ and $b$.  The corresponding group $\pi$ is an infinite cyclic group generated by $a$ and $b=a^{-1}$.  The self-linking function $u^w$ of any nanoword $w$ over
$\alpha$ is  determined by the Laurent polynomial $u^w(a)\in \ZZ
\pi=
\ZZ[ a^{\pm 1}]$. For  integers
$m,n
\geq 1$, consider the nanoword 
$$w=w_{m,n}=A_1 A_2 \cdots A_m B_1 B_2\cdots B_n A_m \cdots A_2 A_1 B_n\cdots B_2 B_1.$$
   with  
$\vert A_i\vert=a$ and    
$\vert B_j\vert =b$ for all $i,j$.  It is clear that $[A_i]_w=b^n$ and $[B_j]_w=a^{-m}$ for all $i,j$. Therefore $$u^w(a)=[a]_w- [\tau (a)]_w= [a]_w-[b]_w=
mb^n -na^{-m}= ma^{-n}-n a^{-m}.$$
 Pick integers
$p,q,s,t\geq 1$ and consider  the nanoword $v$ over $\alpha$ of length $12(p+ q+s+t)$ obtained as the product (in an arbitrary order) of the nanowords 
$$ w_{p,s}, w_{p,t}, w_{q,s}, w_{q,t}, w_{s+t,p}, w_{s+t, q}, w_{ s, p+q}, w_{t, p+q}, 
 w_{   p+q, s+t}.$$     The additivity of the self-linking function
easily implies that
$u(v)=0$.  Let $r$ be the greatest common divisor of $p$ and $s$. Let $H$ be the  subgroup  of $\pi$   generated by
$a^r$.   Suppose  that
$r\geq 2$ and
$q,t$ are prime to
$r$.   It is easy to check that $v^H\simeq w_{p,s}$   and therefore $ u(v^H)=pa^{-s}-s a^{-p}\neq 0$ for $p\neq s$. Thus $u(v^H)$ detects the non-contractibility of $v$.  

4. Let $a,b,w,H$ be   as in Example \ref{nezzexixixremar}.  The nanoword $v=ww$
  lies  in the kernel of $\tilde \gamma$ and has a zero  self-linking function. (Indeed,   $R_a=R_b=\ZZ/2\ZZ$   and therefore $u(ww)=u(w) + u(w)=0$ for any nanoword $w$.) The non-contractibility of $v$  is detected by  
$\mu_{v^H} $, since  $ \mu_{v^H} (a,b)=2
\mu_{w^H} (a,b)=2$.

 \subsection{Applications}\label{wor2597unkn1} We can use the self-linking function to give a homotopy classification of the monoliteral words $a^m$ where $m\geq 3$ and $a\in \alpha$ with $a\neq \tau (a)$.
 Desingularizing $a^m$,  we obtain  a nanoword  $w=(a^m)^d$ of length $m(m-1)$. 
 A direct computation from the definitions gives
 $[b]_w=0$ for all $b\in \alpha-\{a\}$ and $$[a]_w=
 \sum_{1\leq i<j \leq m, i+j\neq m+1} a^{(m+1-i-j) (j-i)}.$$
 Then $u^w(a)=[a]_w, u^w(\tau (a))=-[a]_w  
$ and  $ u^w(b)=0$ for $b\in \alpha-\{a, \tau(a)\}$. Since $u^w\neq 0$, the word $a^m$ is non-contractible. 
We can estimate its norm from below applying the last remark of Sect.\ \ref{vpamxct} to the monomial determined by $i=1, j=[m/2]+1 $. This gives  
\begin{equation}\label{ine} \vert \vert a^m\vert \vert \geq \left [\frac{m} {2}\right ] \, \left [\frac{m-1} {2}\right ] +1.\end{equation}
(Presumably, $\vert \vert a^m\vert \vert=m (m-1)/2$.)

\begin{theor}\label{lipidyd}  Let $a,b$ be letters in $\alpha$ with $a\neq \tau (a), b\neq \tau (b)$. The words $a^m$ and $ b^n$ with $m,n\geq 3$ are homotopic if and only if $a=b$ and $m=n$. \end{theor}
\begin{proof} For $w= a^m $, we can recover both  $a$ and   $m$ from $u(w)=u({w^d})$. Indeed, $a$ is  the only letter on which the value of $u^w$ is a sum of monomials with positive coefficients. The   sum of these coefficients   is equal to
$m(m-1)/2 -[m/2]$.  This implies the theorem.\end{proof}

Let $a \in \alpha$ be a letter  with $a\neq \tau (a)$.  The self-linking function   allows us to show that the word $ a^m$ with $  m\geq 3$ is not homotopic to a word $v=bcbc$ with distinct $b,c\in \alpha$. If $b=\tau(c)$, then $v$ is contractible while $a^m $ is not. If $b\neq \tau(c)$, then the equalities $[b]_v=c, [c]_v=b^{-1}$ imply that the only  non-zero values of the self-linking function $u^v$ are  $   \pm c, \pm b^{-1} $. The sum  of coefficients  of each of these values is   $  \pm  { 1} \neq m(m-1)/2 -[m/2] $.  Hence     $a^m$ is not homotopic to $v$.

\section{ $\alpha$-forms and $\alpha$-pairings}\label{based22}

The interlacement of letters in a nanoword $(\A,w)$ is   reflected in 
the function
 $n_w:\A\times \A
\to
\{-1,0,1\}$ defined in Sect.\ \ref{covfunct}. A study of $n_w$  and related pairings leads   to 
notions of   
$\alpha$-forms and $\alpha$-pairings. We introduce  here an abstract algebraic theory  of  
$\alpha$-forms and $\alpha$-pairings. Connections with words will be discussed in the next section.

 \subsection{ $\alpha$-forms}\label{extend1}  Let $\pi=\pi(\alpha, \tau)$ be the multiplicative abelian group defined in Sect.\ \ref{covfunct}.  Given a set $S$,   a pairing $l:S\times S\to
\pi$  is {\it skew-symmetric} if   $l(A,B)= l(B,A)^{-1}$ and $
l(A,A)=1  $  for all $A,B\in
S$. A  (skew-symmetric)  {\it 
$\alpha$-form}   is  an  $\alpha$-alhabet $\A$ endowed with two
skew-symmetric pairings $n:\A\times \A \to \{-1,0,1\}$  and  $l:\A\times
\A
\to \pi$. We   view $n $ and $l $  as   matrices 
whose rows and columns are numerated by the letters of $\A$.  The skew-symmetry of  $n$ means  that
$n(A,B)=-n(B,A)
 $
 for all $A,B\in \A$ and in particular  $n(A,A)=0$ for all $A \in \A$. 
Two 
$\alpha$-forms  $(\A, n, l)$ and  $(\A', n', l')$ are {\it isomorphic} if there is an isomorphism of $\alpha$-alphabets $\A\to \A'$
transforming $n, l$ into $n', l'$, respectively.

Given two $\alpha$-forms $(\A_1, n_1, l_1)$ and $(\A_2, n_2, l_2)$, we   define their direct sum 
by $$(\A_1, n_1, l_1)\oplus (\A_2, n_2, l_2)=(\A_1\amalg \A_2, n_1\oplus n_2, l_1\oplus l_2)$$ where $\amalg$ is the disjoint union of $\alpha$-alphabets and $\oplus$ is the   direct sum of matrices. This means that the restriction of $n=n_1\oplus n_2$ to $\A_i$ is $n_i$ for $i=1,2$ and 
$n  (\A_1, \A_2)=n  (\A_2, \A_1)=0$. The restriction of $l=l_1\oplus l_2$ to $\A_i$ is $l_i$ for $i=1,2$ and 
$l  (\A_1, \A_2)=l  (\A_2, \A_1)=1$. Here we   assume  $\A_1$ and $\A_2$ to be disjoint; if it is not the case we replace $(\A_1, n_1, l_1)$ by an isomorphic $\alpha$-form and proceed as above.
The {\it trivial $\alpha$-form}  $(\A=\emptyset, n, l)$ yields a neutral element for the direct sum.

For each $\alpha$-form $f=(\A, n, l)$, we   define  the {\it opposite} $\alpha$-form    $f^-=(\A, n^-,l^-)$ where $n^-=-n$ and $l^-(A,B)=(l(A,B))^{-1}=l(B,A)$ for $A,B\in \A$.
We also define   the {\it inverse} 
$\alpha$-form   $\overline f=(\overline \A,  n ,l^-)$ where $\overline A$ is the inverse $\alpha$-alphabet.

 \subsection{ Relation of homology}\label{edrextend1} We define three transformations  (moves) on  $\alpha$-forms.

(i)$^*$  The move applies to an  $\alpha$-form $(\A, n, l)$   if there is a letter $A\in \A$ such
that $n(A,\A)=0$ and   $l (A,\A)=1$. The move replaces $(\A, n, l)$ with  $(\A'= \A-\{A\},
n\vert_{{\A'}}, l\vert_{{\A'}})$.

(ii)$^*$  The   move applies to an  $\alpha$-form $(\A, n, l)$   if there are letters $A,B\in \A$ such
that $\vert A\vert =\tau (\vert B\vert) $   and $  n(A,C)=n(B,C), l (A,C)=l(B,C) \vert
B\vert^{n(B,C)} $
for all $C\in \A $.
The move replaces 
$(\A, n, l)$ with  $(\A'= \A-\{A, B\},
n\vert_{{\A'}}, l\vert_{{\A'}})$.

Note that for $C=A$ and $C=B$, the conditions  $  n(A,C)=n(B,C), l (A,C)=l(B,C) \vert
B\vert^{n(B,C)} $ are equivalent to $n(A,B)=0, l (A,B)=1$. Note also that the formula $l
(A,C)=l(B,C)
\vert B\vert^{n(B,C)} $ is symmetric in $A,B$. Indeed, this formula 
is equivalent to  $ l (B,C)=l(A,C)
\vert A\vert^{n(A,C)} $ since $\vert A\vert=\tau (\vert B\vert)= \vert B\vert^{-1}\in \pi$  and $ 
n(A,C)=n(B,C)$.

(iii)$^*$  The  move applies to an  $\alpha$-form $(\A, n, l)$   if there are letters $A,B, C\in \A$ such
that $\vert A\vert =  \vert B\vert =  \vert C\vert$ and $ n(A,B)=n(B,C)=1, n(A,C)=0 $.
The move replaces 
$(\A, n, l)$ with  $( \A ,n',l')$ where   $n',l'$ coincide  with $n, l$,
respectively, except on the pairs $(A,B),
(B,A),
(A,C), (C,A), (B,C), (C,B)$. The values of  $n',l'$  on these pairs  are 
determined by
$n'(A,B)=n'(B,C)=0, n'(A,C)=1$ and  $$l'(A,B)=l (A,B)\,\vert C\vert ,\,\,\, l'(A,C)=l (A,C)\,\vert
B\vert^{-1},\,\,\, l'(B,C)=l (B,C)\,\vert A\vert  .$$

Two 
$\alpha$-forms are {\it homologous} if they can be related by a finite sequence of  moves (i)$^*$--(iii)$^*$,  the
inverse moves, and isomorphisms.  This defines an equivalence relation of {\it homology} in the class of $\alpha$-forms.  It is easy to check  that   $\alpha$-forms    opposite (resp.\ inverse) to homologous $\alpha$-forms are themselves homologous.

 \subsection{$\alpha$-pairings}\label{23extend1}       To study $\alpha$-forms,  we introduce a
simpler concept  of an $\alpha$-pairing. The  principal advantage of $\alpha$-pairings  is that
for them the   move (iii)$^*$ becomes trivial.

A  (skew-symmetric)  {\it $\alpha$-pairing}   is  a set $S$ endowed with a distinguished element
$s\in S$,    a
mapping $S-\{s\}\to \alpha$, and   a skew-symmetric
pairing $b:S\times
 S  
\to \pi=\pi(\alpha, \tau)$.  The conditions on $S $ can be rephrased by saying that $S$ is a disjoint union of an $\alpha$-alphabet and a distinguished element $s$. As usual, the image of   $A\in S-\{s\}$ under the
projection  to $\alpha $
 is denoted $\vert A\vert$. Two 
$\alpha$-pairings  $(S, s,b)$ and  $(S', s',b')$ are {\it isomorphic} if there is a bijection $S\to S'$ transforming $s$  to $s'$, inducing an isomorphism of 
$\alpha$-alphabets $S-\{s\}\to S'-\{s'\}$, and 
transforming $b$ into $b'$.

The {\it direct sum} of  two $\alpha$-pairings $(S_1, s_1, b_1)$ and $(S_2, s_2, b_2)$ is the $\alpha$-pairing  $$(S_1, s_1, b_1)\oplus (S_2, s_2, b_2)= ((S_1-\{s_1\} )\amalg (S_2-\{s_2\} )\amalg \{s\} , s, b)$$ where  $b$ is uniquely determined by the following conditions: $b(S_1-\{s_1\}, S_2-\{s_2\})=1$,  for $i=1,2$,  the restriction of $b$  to $S_i-\{s_i\}$ is $b_i$,   and 
$b(A,s)=b_i(A,s_i)$ for all $A\in S_i-\{s_i\}$.  Here we   assume  $S_1$ and $S_2$ to be disjoint; if it is not the case we replace $(S_1, s_1, b_1)$ by an isomorphic $\alpha$-pairing and proceed as above. The {\it trivial $\alpha$-pairing}   $(S=\{s\}, s, b=1)$  yields a neutral element for the direct sum.

For each $\alpha$-pairing $p=(S, s,b)$, we   define  the {\it opposite} $\alpha$-pairing    $p^-=(S, s ,b^-)$ where   $b^-(A,B)=(b(A,B))^{-1}=b(B,A)$ for $A,B\in S$.
We   define   the {\it inverse} 
$\alpha$-pairing   $\overline p=(\overline S, s ,   b^-)$ where $\overline S$ is the   set $S$ with the given mapping $S-\{s\}\to \alpha$  replaced by  its composition with $\tau:\alpha\to \alpha$.

 \subsection{ Relation of homology for $\alpha$-pairings}\label{chcemqpnd1}  
 Given an $\alpha$-pairing $(S, s,b)$, we say that $A\in S-\{s\}$ is an {\it annihilating element}   if $b(A,C)= 1 $ for all $C\in S$. Two elements
$A,B\in
S-\{s\}$  are   {\it twins} if $\vert A\vert =\tau (\vert B\vert)$   and $  b(A,C)=b(B,C) $
for all $C\in S $. Note that if $A,B$ are twins, then $b(A,B)=b(B,B)=1$.   We define two \lq\lq compression" moves $M_1, M_2$  on   $(S, s,b)$.
  The  move 
$M_1$ deletes  an annihilating element $A\in S-\{s\}$ that is replaces $(S,s,b)$ with  $(S'= S-\{A\}, s\in S', 
b\vert_{{S'}} )$.
  The   move $M_2$ deletes a pair of twins     $A,B\in
S-\{s\}$ that is replaces 
$(S,s,b)$ with  $(S'= S-\{A, B\}, s\in S',
b\vert_{{S'}} )$.

 Two $\alpha$-pairings are {\it homologous} if they can be related by a finite sequence of  moves $ M_1, M_2 $, the
inverse moves, and isomorphisms. This defines    an equivalence relation of {\it homology}  in the  class of $\alpha$-pairings.   It is easy to check  that   $\alpha$-pairings    opposite (resp.\ inverse) to homologous $\alpha$-pairings are themselves homologous.

We now classify  
$\alpha$-pairings up to homology. Let us call an $\alpha$-pairing    {\it primitive} if it has no annihilating elements and no pairs of twins. 
 The following lemma   shows that  homology classes of $\alpha$-pairings bijectively correspond to isomorphism classes of primitive   $\alpha$-pairings.

\begin{lemma}\label{l:t51} Every $\alpha$-pairing can be transformed by the moves $M_1, M_2$ into a primitive   $\alpha$-pairing. Two homologous primitive $\alpha$-pairings are isomorphic.
		   \end{lemma}
                     \begin{proof}  The first claim is obvious: 
deleting  
annihilating elements  and   pairs of twins by 
  $M_1, M_2$  we can  compress any 
$\alpha$-pairing   into a  primitive one.  
   To prove the second claim, we need the following assertion:
  
  $(\ast)$ for any $i,j \in \{1,2 \}$, a  move $M_i^{-1}$ followed by   $M_j$ yields the same 
result as    
an isomorphism, or a   move $ M_k $, or     a  move $M_k $ 
followed 
by $M_l^{-1}$ with $k,l\in \{1,2 \}$.

  This assertion will imply the second claim of the lemma. Indeed,    
suppose that two primitive  $\alpha$-pairings $p,p'$ are related by a finite 	
sequence 
of moves $M_1^{\pm 1}, M_2^{\pm 1} $ and 
isomorphisms. 	
 An isomorphism   followed by $ M_i^{\pm 1}$     
can be  
decomposed  as  $ M_i^{\pm 1}$    followed by an isomorphism.  Therefore 
all 
isomorphisms in our sequence   can be   accumulated at the end. The 
claim  
$(\ast)$ implies that $p, p'$
can be  related by a finite sequence of moves consisting of  several 
  moves of type $  M_i$ followed by  several 
  moves of type $ M_i^{-1}$ and   isomorphisms. Since   $p$ is 
primitive  
  we cannot apply to it a move of type $  M_i $. Hence there are no 
such 
moves in our sequence. Similarly, since $
  p'$ (and any isomorphic $\alpha$-pairing) is primitive, it cannot be obtained by 
an 
application of $M_i^{-1}$. Therefore our sequence  consists solely of 
isomorphisms so 
that $p$ is isomorphic to $p'$.

 Let us   prove $(\ast)$. We have to consider  four  cases depending on 
$i,j\in \{1,2 \}$. 
 
  For $i= j =1$, the move   $M_i^{-1}$ on an  $\alpha$-pairing 
$(S,s,b)$ adds 
one element $A$   and then  $M_j$ removes   one  element  $A'\in S\amalg \{A\}$. If 
$A'=A$, 
then   $  M_j\circ M_i^{-1}=\id$. If   $A'\neq A$, then 
$A'\in S$ 
is an annihilating element. The 
transformation $  
M_j\circ M_i^{-1}$   can be achieved by first applying $M_j $ that 
removes 
$A'$ and then applying $M_i$ that adds $A$.

   Let $i=1, j=2$. The move 
 $M_i^{-1}$ on $(S,s,b)$ adds an annihilating element $A$  and  
  $M_j$ removes a pair of twins $A_1,A_2\in S\amalg \{A\}$. 
  If   $A_1,A_2\in S$, then  $  
M_j\circ M_i^{-1}$ 
 can be achieved by first removing $A_1,A_2$ and then adding $A$. If 
$A_1=A$, 
then   $A_2 $ is an annihilating element of $S$ and $  M_j\circ M_i^{-1}$ is the 
move 
$M_1$    removing $A_2$. The case $A_2=A$ is similar.

   Let $i=2, j=1$. The move 
 $M_i^{-1}$  on $(S,s,b)$ adds a pair of twins $A_1,A_2$   and  
$M_j$ removes an annihilating element  $A\in S\amalg \{A_1,A_2\}$. 
  If   $A\in S$, then  $  M_j\circ 
M_i^{-1}$  can 
be achieved by first removing $A$ and then adding $A_1,A_2$.
   If $A=A_1$, then   $ A_2$ is an annihilating element of $S\amalg \{A_2\}$   
and $  
M_j\circ M_i^{-1}=M_1^{-1}$. The case $A=A_2$ is similar.

      Let $i=j=2$. The move 
 $M_i^{-1}$  on $(S,s,b)$ adds a pair of twins $A_1,A_2$   and  
$M_j$ removes a pair of twins  $A'_1,A'_2\in S \amalg \{A_1,A_2\}$. If these 
two pairs  
 are disjoint, then 
  $  M_j\circ M_i^{-1}$  can be achieved by first removing  $A'_1,A'_2\in S$ 
and then 
adding $A_1,A_2$. If these two pairs  
  coincide, then $  M_j\circ M_i^{-1}=\id$. It remains to 
consider 
the case where these pairs have one common element, say $A'_1=A_1$, 
while 
$A'_2\neq A_2$. Then $A'_2\in S$ and for all $C\in S$,
  $  b(A'_2,C)=   b(A'_1,C)=   b(A_1,C)=  
b(A_2,C)$.
 Therefore the move   $  M_j\circ M_i^{-1}$ gives an $\alpha$-pairing isomorphic 
to 
$(S,s,b)$. The isomorphism $S\to (S-\{A'_2\}) \cup \{A_2\}$ is the 
identity on 
$S- \{A'_2\}$ and  sends $A'_2$ into $A_2$.
\end{proof}

  For an  $\alpha$-pairing $p$ we  denote the unique (up to isomorphism)   primitive $\alpha$-pairing homologous to $p$ by $p_\ast$. It follows from the definitions that $(p^-)_\ast =(p_\ast)^-$, 
  ${\overline p}_\ast= \overline {p_\ast}$, and $(p_1\oplus p_2)_\ast = (p_1)_\ast \oplus (p_2)_\ast$.

 \subsection{ Homology invariants of  $\alpha$-pairings}\label{dopa3end1}
 Any isomorphism invariant of primitive $\alpha$-pairings yields a homology invariant of an $\alpha$-pairing $p=(S,s,b)$    by computing   on $p_\ast=(S_\ast, s_\ast, b_\ast)$. For example, 
 $\rho (p)=\card (S_\ast) -1$ is a homology invariant of $p$. Given   $a\in  \alpha, x\in \pi$, the number
 $$\rho_{a,x} (p)=\card \{A\in S_\ast-\{s_\ast\} \,\vert\, \vert A\vert =a, \, b_\ast (A,s_\ast)=x\}$$  is a homology invariant of $p$. Clearly $\rho (p)= \sum_{a\in \alpha, x\in \pi } \rho_{a,x} (p)$. We have $\rho_{a,x} (p^-)= \rho_{a,x^{-1}} (p)$ and $\rho_{a,x} (\overline p )= \rho_{\tau (a),x } (p)$. Both  $\rho $   and $\rho_{a,x}$ are additive with respect to the direct sum of $\alpha$-pairings.

 Imitating the self-linking function of nanowords, we can define a  homology invariant of   $p=(S,s,b)$  as follows.   
  For      $a\in \alpha$,  set $$[a]_p= \sum_{ A \in S-\{s\},  b(A,s) \neq 1, \vert A\vert =a} b(A,s)
\in
\ZZ\pi.$$
Define   a function $u^p=u(p)$  on $\alpha$    by 
 $ u^p(a)= 
[a]_p (\modu 2) \in (\ZZ /2\ZZ) \pi$ if 
$\tau (a)=a$ and $
  u^p(a)=  [a]_p-[\tau(a)]_p \in \ZZ\pi$ if $  \tau (a)\neq a$.   It is obvious that $u^p $  is   invariant under the moves $M_1, M_2$ and is thus a  homology invariant of $p$.
For  $\tau (a)=a$, $$ u^p(a)= u^{p_\ast} (a)= \sum_{x\in \pi-\{1\}} \rho_{a,x} (p) \, x\, (\modu 2)$$  and  for  $  \tau (a)\neq a$, $$
  u^p(a)= u^{p_\ast} (a)=  \sum_{x\in \pi-\{1\}} (\rho_{a,x} (p) - \rho_{\tau (a),x} (p))  \,x.$$

 \subsection{ From $\alpha$-forms to $\alpha$-pairings}\label{dddextend1}
Every  $\alpha$-form $ (\A, n, l)$
determines an $\alpha$-pairing as follows. Set 
$S=\A\cup \{s\}$ where $s$ is not a letter of the alphabet $\A$.
There is a unique skew-symmetric pairing  $b:S\times S\to \pi$ such that for any $A,B\in
\A\subset S $, 
$$ b(A,s)= \prod_{ C \in   \A}   \vert C\vert^{n  (A,C)} \in \pi
\,\,\,\,{\rm {and}}\,\,\,\, 
b(A, B)=(l(A,B))^2 \vert A\vert^{n(A,B)}\vert B\vert^{n(A,B)} \in \pi.$$  
The $\alpha$-pairing $(S,s,b)$    { associated} with $f=(\A, n, l)$ is denoted $p(f)$. The operation $f\mapsto p(f)$ commutes with direct sum and the operations of taking the opposite and the inverse. 

\begin{lemma}\label{exttobased}    $\alpha$-pairings determined by  homologous  $\alpha$-forms are   homologous.
\end{lemma}
                     \begin{proof} It suffices to analyze the behavior of the
associated $\alpha$-pairing under the moves (i)$^*$--(iii)$^*$ on   $\alpha$-forms. 
Consider the move (i)$^*$ replacing an  $\alpha$-form $(\A, n, l)$   with  $(\A'= \A-\{A\},
n\vert_{{\A'}}, l\vert_{{\A'}})$  where  $A\in \A$ is such
that $n(A,\A)=0$ and   $l (A,\A)=1$.  Let $(S=\A\cup \{s\},s,b)$ be the $\alpha$-pairing associated with  $(\A, n, l)$. 
It follows from the definitions that 
$A$ is an annihilating element for $b$.
  The  $\alpha$-pairing associated with $(\A' ,
n\vert_{{\A'}}, l\vert_{{\A'}})$ is $(S'= S-\{A\}, s\in S', 
b\vert_{{S'}} )$.  
It is obtained from $(S ,s,b)$ by the move $ M_1 $ deleting $A$.

Consider the move (ii)$^*$ replacing   $(\A, n, l)$   with  $(\A'= \A-\{A,B\},
n\vert_{{\A'}}, l\vert_{{\A'}})$  where  $A,B\in \A$ are such
that $\vert A\vert =\tau (\vert B\vert) $   and $  n(A,C)=n(B,C), l (A,C)=l(B,C) \vert
B\vert^{n(B,C)} $
for all $C\in \A $.  Let $(S=\A\cup \{s\},s,b)$ be the $\alpha$-pairing associated with  $(\A, n, l)$. For $C\in \A$,
$$ b(A,C)=  (l(A,C))^2 \vert A\vert^{n(A,C)}\vert C\vert^{n(A,C)} =(l(B,C))^2 \vert B\vert^{2 n(B,C)}
 \vert A\vert^{n(B,C)} \vert C\vert^{n(B,C)}  $$ $$
= (l(B,C))^2 \vert B\vert^{  n(B,C)}
   \vert C\vert^{n(B,C)}=b(B,C). $$ Similarly, $  b(A,s)=b(B,s) $. Therefore $A$ and $B$ are twins. The  $\alpha$-pairing associated with $(\A' ,
n\vert_{{\A'}}, l\vert_{{\A'}})$ is $(S'= S-\{A,B\}, s\in S', 
b\vert_{{S'}} )$. 
It  is obtained from $(S ,s,b)$ by the move $ M_2 $ deleting $A,B$.

Finally,  under a move    $(\A, n, l)\mapsto (\A, n', l')$ of type (iii)$^*$, the associated $\alpha$-pairing $(S,s,b)$ is preserved. Indeed, let  $(S,s,b')$ be the $\alpha$-pairing associated with $(\A, n', l')$.
It is clear that  $b'$ coincides with $b$ except possibly on the pairs $(A,B),
(B,A), 
(A,C), (C,A), (B,C), (C,B)$.    We have
$$b'(A, B)=  (l'(A,B))^2 \vert A\vert^{n'(A,B)}\vert B\vert^{n'(A,B)} =  (l'(A,B))^2=   (l(A,B)  \vert C\vert) ^2$$
$$=   (l(A,B))^2 \vert A\vert\,\vert B\vert=  (l(A,B))^2 \vert A\vert^{n(A,B)}\vert B\vert^{n(A,B)}=b(A, B).$$
 Similarly, $b'(B,C)=b(B,C)$ and $b'(A, C)=b(A, C)$. By skew-symmetry, $b'=b$.
\end{proof}

Passing from an $\alpha$-form  to the associated $\alpha$-pairing we   loose certain information contained in the 2-subgroup of $\pi$. It is interesting to find further invariants of 
 $\alpha$-forms retaining this information. 
 
 \section{Linking forms of nanowords}\label{edrbased22}

We associate with each  nanoword $w$ over $\alpha$ an $\alpha$-form and an $\alpha$-pairing  reflecting   the 
 linking  of
letters in $w$.

 \subsection{Notation}\label{nota1}  
 Given a  nanoword
$(\A, w:\hat n
\to
\A)$ and a  letter
$A\in
\A$,  we   denote by $i_A$ (resp.\ $j_A$) the minimal (resp.\ the maximal) element of the 2-element set $\omega^{-1} (A)\subset \hat
n$. Thus  $w=\cdots
 A \cdots  A \cdots $ where the first entry of $A$ appears on the $i_A$-th position and the second entry of $A$ appears on the $j_A$-th position.

 \subsection{Linking form}\label{based1} For  a  nanoword
$(\A, w:\hat n
\to
\A)$, we define an  $\alpha$-form $f^w=(\A,n_w,lk_{w} ) $ called the {\it linking form} of $w$. The function  $n_{w}$ was defined in
Sect.\ \ref{covfunct}.   For   letters $D,E\in
\A$, set $$
D\circ E= \prod_{F\in \A, i_D<i_F<j_D \, {\rm {and}} \, i_E<j_F<j_E} \vert  F \vert \in \pi
 $$
and 
$$lk_w(D,E)= (D\circ E) (E\circ D)^{-1}\in \pi.$$
 In particular, if  
$w=\cdots
 D \cdots  D \cdots  E
\cdots  E
\cdots$, then
$E\circ D=1$ and   
$lk_w(D,E)$ $ =D\circ E$. If $w=\cdots  D \cdots  E \cdots  E \cdots  D \cdots$, then
$$lk_w(D,E)= \prod_{F\in \A, i_D<i_F<  i_E <  j_F<j_E} \vert  F \vert
\prod_{F\in \A, i_E<i_F<  j_E <  j_F<j_D} \vert  F \vert^{-1}.$$
If $w=\cdots  D \cdots  E \cdots  D \cdots  E \cdots$, then
$$lk_w(D,E)= \prod_{F\in \A, i_D<i_F<  i_E <  j_F<j_D\, {\rm {or}} \,  i_E<i_F<  j_D <  j_F<j_E  \, {\rm {or}} \,  i_D<i_F< i_E< j_D <  j_F<j_E} \vert  F \vert.$$
For all other positions of $D, E$ in $w$, we can use these formulas to compute $lk_w(E,D)$ and then   use the  
skew-symmetry of $lk_w$  to compute $lk_w(D,E)$.

\begin{lemma}\label{extmatrad}    Under the homotopy moves (i), (ii), (iii) on nanowords, the associated 
$\alpha$-forms are transformed by the   moves (i)$^*$, (ii)$^*$, (iii)$^*$, respectively.
\end{lemma}
                     \begin{proof}     Consider the first  homotopy move
$w=(\A,  x A A  y) 
\mapsto  (\A- \{A\},  x  y)=v$.  Pick $D,E\in \A- \{A\}$. The   formulas for $lk_w(D,E)$   show that if a letter $F\in \A$  
contributes
non-trivially to  $lk_w(D,E)$, then $i_F<j_F-1$.
Hence  the letter $A$ does not contribute to $lk_w(D,E)$. This implies that     $lk_{v}  $ is the restriction of $lk_w$ to $\A'$. It is also obvious that
$n_{v}  $ is the restriction of $n_w$ to $\A'$.  It  remains to observe  that $n_w(A,\A)=0$ and $lk_w(A,\A)=1$. Thus $f^v$   
is obtained from $f^w$  by   (i)$^*$.

Consider the      move
$w=(\A,  x A B  y BAz) 
\mapsto  (\A- \{A,B\},  x  y z)=v$  where   $\vert A\vert =\tau (\vert B\vert)$. It follows from the definitions 
that $n(A,B)=0$ and $ lk_w (A,B)=1$. A case by case study of all possible positions of   two
entries of a letter 
$C\in \A-\{A,B\}$ in $w$ with respect to $AB, BA$  gives $
 n(A,C)=n(B,C)$ and $ lk_w (A,C)=lk_w(B,C) \vert B\vert^{n(B,C)} $. It is obvious that $n_v$ is the restriction of $n_w$  to $
\A'=\A-\{A,B\}$. Since the contributions of the
letters $A,B$ to $lk_{w}(C,C')$  
cancel each other for all $C,C'\in \A-\{A,B\}$, we obtain that $lk_{v}  $ is the restriction of $lk_w$ to $\A'$.  Thus $f^v$
is obtained from $f^w$  by   (ii)$^*$.

Consider the     move $w=(\A, xAByACzBCt)\mapsto  (\A, xBAyCAzCBt)=v$ 
where $ \vert A\vert =  \vert B\vert
=  \vert C\vert$. 
   It follows from the definitions
 that    $ n_{w}(A,B)=n_{w}(B,C)=1, n_{w}(A,C)=0 $ and $n_v(A,B)=n_{v}(B,C)=0, n_{v}(A,C)=1$.
It is obvious that for any $D, E\in \A-\{A,B, C\}$, we have $lk_{w} (D,E)= lk_{v} (D,E)$.  A case by case study of  
possible positions of the two entries of    $D$ in $w$ 
 give  $$lk_{w} (D,A)= lk_{v} (D,A), \,\,\,\, lk_{w} (D,B)= lk_{v} (D,B), \,\,\,\,  lk_{w} (D,C)= lk_{v} (D,C)$$   for any $D
\in
\A-\{A,B, C\} $.  It remains to observe that 
$$
lk_{v} (A,B)=lk_{w}
(A,B) \,\vert C\vert, \,\,\,\,  lk_{v} (A,C)=lk_{w}  (A,C)\, \vert B\vert^{-1}$$
and $   lk_{v}(B,C)=lk_{w} (B,C) \vert A\vert 
$. 
Thus $f^v$
is obtained from $f^w$  by (iii)$^*$.
\end{proof}   

This lemma implies that the homology class of the   $\alpha$-form $f^w=(\A,n_w,lk_{w} ) $  is a homotopy invariant of $w$. This form behaves in a natural way under the standard operations on nanowords: $$f^{\overline w} =\overline {f^w}, \quad   f^{ w^-} =  (f^w)^-,\quad  f^{w_1w_2} =
f^{w_1 } \oplus f^{ w_2}.$$

 \subsection{Linking pairing}\label{based112121216497df} For  a  nanoword
$w$,   the $\alpha$-pairing   associated with the $\alpha$-form $f^w $  is denoted $p^w$ and called the {\it linking pairing} of $w$. 
It   behaves  in a natural way under the standard operations on nanowords: $$p^{\overline w} =\overline {p^w}, \quad   p^{ w^-} =  (p^w)^-,\quad  p^{w_1w_2} =
p^{w_1 } \oplus p^{ w_2}.$$

As we know, the $\alpha$-pairing $p^w$ can be compressed to a primitive $\alpha$-pairing   $p^w_\ast=(S_\ast,s_\ast, b_\ast: S_\ast \times S_\ast \to \pi)$. The isomorphism class of $p^w_\ast$ is a homotopy invariant of $w$.   If $w$ is contractible,  then $p^w_\ast$ is a  trivial $\alpha$-form. The formulas above imply that 
$$p_\ast^{\overline w} =\overline {p_\ast^w}, \quad   p_\ast^{ w^-} =  (p_\ast^w)^-,\quad  p_\ast^{w_1w_2} =
p_\ast^{w_1 } \oplus p_\ast^{ w_2}.$$
The homology invariants of  $\alpha$-pairings defined in Sect.\ \ref{dopa3end1} yield homotopy invariants of nanowords. In particular, set $\rho (w)=\rho (p^w)=\card (S_\ast)-1$ and for $a\in \alpha , x\in \pi$, set $$\rho_{a,x} (w)=\rho_{a,x} (p^w) =\card \{A\in S_\ast-\{s_\ast\} \,\vert\, \vert A\vert =a, \,\,b_\ast (A,s_\ast)=x\} .$$
We have      $\rho_{a,x} (\overline w) =\rho_{\tau (a),x} (w)$, $\rho_{a,x} (w^-) =\rho_{a, x^{-1}} (w)$ and $\rho_{a,x} (w_1w_2)=\rho_{a,x} (w_1)+\rho_{a,x} (w_2)$. Application:  if $\rho_{a,x} (w)\neq 0$ for some $a,x$, then   $ww'$ is non-contractible for any nanoword $w'$.

  For any nanoword $(\A,w)$, the homotopy invariance of $\rho$ implies that   
 $\vert \vert w\vert \vert\geq \rho (w)$.  Application:  if the $\alpha$-pairing   $p^w$ 
   is primitive, then $w$ is not homotopic to a nanoword of   smaller length and $\vert \vert w\vert \vert=\rho(w)=\card (\A)$. 
   
   Note finally that the self-linking function $u^w$ of a nanoword $w$ is computed from   $p_\ast^w$: it follows from the definitions that 
 $u^w= u (p^w)= u(p^w_\ast)$.

 \subsection{Nanowords of length 4}\label{nanexacacacacota1} It follows from the definitions  that all nanowords of length $\leq 2$ are contractible.  We  give a   homotopy classification of nanowords of length 4. Clearly, each such nanoword is either contractible or isomorphic to the nanoword 
 $(\A=\{A,B \}, w=AB AB )$ with $\vert A\vert =a  \in \alpha,  \vert B\vert =b  \in \alpha$.  Denote the latter nanoword by $w_{a,b}$  (possibly $a=b$). Clearly, $\overline{w_{a,b}}=w_{\tau(a), \tau(b)}$
  and $(w_{a,b})^-\approx w_{b,a}$.
  
  By the second homotopy move,  if $b=\tau(a)$, then $w_{a,b}$ is contractible. The following theorem establishes the converse and gives a homotopy classification of these nanowords.

  \begin{theor}\label{classid}   For $a, b\in   \alpha$ with $\tau(a)\neq b$, the nanoword $w_{a,b}$ is non-contractible. Two nanowords $w_{a,b}$ and $w_{a',b'}$ with  $\tau(a)\neq b, \tau(a')\neq b'$ are homotopic if and only if $a=a'$ and $b=b'$.     \end{theor}
    \begin{proof}  Set $w=w_{a,b}$. The linking $\alpha$-form $(\A, n_w, lk_w)$   is given by the matrices
 $$ n_w  = \left [ \begin{array}{ccccc}   
        0& 1   \\
-1& 0  \\
\end{array} \right ] , \,\,\,\, lk_w= \left [ \begin{array}{ccccc}   
        1& 1   \\
 1& 1  \\
\end{array} \right ].$$
Here the 1-st and 2-nd rows (and columns) correspond to   $A,B $, respectively.   The associated   $\alpha$-pairing $p^w=(S, s, S\times S\to \pi)$  with   $S=\{s, A, B \}$ is 
given by the matrix 
\begin{equation}\label{1ner}  p^w= \left [ \begin{array}{ccccc}   
        1& b^{-1}    & a  \\
b  & 1 & ab   \\
a^{-1}   & a^{-1}b^{-1} & 1  \\
\end{array} \right ] .\end{equation}
Here the   rows  (resp.\  columns)  correspond to   $s, A,B\in S$ counting from top to bottom (resp.\ from left to right).
Since $b\neq \tau(a)  $, we have $ab\neq 1\in \pi $  so that $p^w$ has no annihilating elements and no twins. Hence $p^w$ is   primitive     and  $\vert \vert w\vert \vert= \rho (w)=2$.

Suppose that the  nanowords $w=ABAB$ and $w'=A'B'A'B'$  with  $\vert A\vert =a,  \vert B\vert =b \neq \tau(a), \vert A'\vert =a',  \vert B'\vert =b'\neq \tau(a')$ are homotopic. Then the associated  primitive $\alpha$-pairings are isomorphic.  
If the isomorphism sends $A,B$ to $A', B'$, respectively, then $ a'=a, b'=b$  and we are done. If  the isomorphism sends $A$ to $B'$ and $B$ to $ A'$, then $ a'=b$ and $  b'=a$. 
We claim that in the latter case  $a=b$. If $a\neq b$, then $a,b$ lie in different orbits $\hat a, \hat b\subset \alpha$ of $\tau$. 
    Recall  from Sect.\ \ref{applega}  the functions $\mu_w, \mu_{w'}:\hat \alpha \times \hat \alpha \to \ZZ$  where $\hat \alpha$ is the set of orbits of $\tau$.  Clearly,   $\mu_w (\hat a, \hat b)=1$ and $\mu_{w'}(\hat a, \hat b)=-\mu_{w'}(\hat b, \hat a)=-1$. This contradicts the invariance of  $\mu_w$   under homotopy. Therefore   $a=b=a'=b'$.
 \end{proof} 
 
  \begin{corol}\label{cororlspsid}   A nanoword of length $\leq 4$ is homotopically skew-symmetric if and only if it is contractible. A nanoword of length $\leq 4$ is homotopically  symmetric if and only if it is contractible or symmetric.      \end{corol}
    \begin{proof}   The nanoword  $w=w_{a,b}$ is homotopically skew-symmetric if and only if  $w \simeq w_{\tau(b), \tau (a)} $. By Theorem \ref{classid},  this happens if and only if $ \tau(a)=b$.
  The nanoword  $w=w_{a,b}$ is homotopically  symmetric if and only if   $w \simeq w_{b,a} $. By Theorem \ref{classid},  this happens if and only if  $ \tau(a)=b$ or $a=b$.
 \end{proof} 

    \begin{corol}\label{homsymspsid} The group of homotopy automorphisms  of $w_{a,b}$ is the full group of 
   $\tau$-equivariant bijections $\alpha\to \alpha$ if $a=\tau(b)$ and its  subgroup consisting of bijections preserving $a$ and $b$ if $a\neq \tau (b)$.    \end{corol}

 \subsection{Example}\label{exacacacacota1}  The following example shows that for  the 
 nanowords of length $ 6$, the associated $\alpha$-pairing is a   strong but non-faithful invariant.  
   Pick three letters $a,b,c\in \alpha$ (possibly coinciding)
and consider the  $\alpha$-alphabet   $\A=\{A,B,C\} $ where $\vert A\vert =a,  \vert B\vert =b, \vert C \vert = c$.    Consider the nanoword $(\A , w=ABACBC)$.  The linking $\alpha$-form $(\A, n_w, lk_w)$   is given by  
 $$ n_w  = \left [ \begin{array}{ccccc}   
        0& 1 & 0 \\
-1& 0 & 1\\
0& -1 & 0\\
\end{array} \right ] , \,\,\,\, lk_w= \left [ \begin{array}{ccccc}   
        1& 1 & b \\
 1& 1 & 1\\
b^{-1} &  1 & 1\\
\end{array} \right ].$$
Here the   rows  and columns  correspond to  $A,B,C$, respectively. For example, $lk_w(A,B)=1\in \pi,  lk_w(A,C)=b\in \pi$, etc. The associated   $\alpha$-pairing $p^w=(S, s, S\times S\to \pi)$  with   $S=\{s, A, B, C\}$ is 
given by  
$$ p^w=  \left [ \begin{array}{ccccc}   
        1& b^{-1}   & ac^{-1} & b  \\
b  & 1 & ab  &  b^2  \\
a^{-1} c & a^{-1}b^{-1} & 1 & bc \\
b^{-1}  &  b^{-2}  & b^{-1}c^{-1} & 1\\
\end{array} \right ] .$$
Here the   rows  and columns  correspond to   $s, A,B,C $, respectively.

 Observe
  that all the generators $\{e\}_{e\in \alpha}$  of $\pi$ represent non-trivial elements of $\pi$. More generally,  a product of an odd number of generators is   a non-trivial element  of $\pi$.
This implies that $b^{\pm 1}\neq a^{-1} c$ so that $A, B$ are not twins and  $B, C$ are not twins. It is also clear that $A$ and $C$ are not annihilating elements. The element $B\in S$ is annihilating if and only if $a=c=\tau (b)$.  The letters $A$ and $C$ are twins   if and only if $a=\tau (c)$ and $b=\tau (b)$.  In all other cases $p^w$ is primitive and $\vert \vert w\vert \vert =3$.  If $a=c=\tau (b)$, then deleting $B$ we obtain a homologous $\alpha$-pairing 
$$ \left [ \begin{array}{ccccc}   
        1& a    & a^{-1}  \\
a^{-1}  & 1    &  a^{-2}  \\
a  &  a^{2}    & 1\\
\end{array} \right ] .$$ If  $a\neq \tau (a)$, then  this $\alpha$-pairing is primitive and  $3\geq \vert \vert w\vert \vert \geq \rho(w)=2$. If $a=b=c=\tau(a)$, then $w$ is contractible by $w\simeq BACACB \simeq BB \simeq \emptyset$. If $a=\tau (c)$ and $b=\tau (b)$, then deleting $A,C$ we obtain a homologous $\alpha$-pairing 
$$ \left [ \begin{array}{ccccc}   
        1& a^2  \\
a^{-2}  & 1       \\
\end{array} \right ] .$$
If $a\neq \tau (a)$, then it  is primitive and $3\geq \vert \vert w\vert \vert \geq \rho(w)=1$.  In the  remaining  case     $a=\tau(a)=c \neq b=\tau(b)$, the $\alpha$-pairing      $p^w$ is 
homologous to a trivial $\alpha$-pairing  and gives no information about $w$.
The homomorphisms $\gamma$ and $\tilde \gamma$ do not help    since in this case   $  \tilde \gamma (w)=1$.     However,   $w$  is non-contractible by Theorem \ref{1classid} below. 

The nanoword $ABACBC$ with $\vert A\vert =\vert B\vert =\vert C\vert =c$ is the desingularization of the word $ccc$ in the alphabet $\alpha$. Therefore if $c\neq \tau(c)$, then
\begin{equation}\label{normccc} \vert \vert ccc\vert \vert =3. \end{equation}

\section{Analysis litterae}\label{cl555v}

 We give   a homotopy classification of words of length $\leq 5$.
Since the desingularization of a word deletes all letters of multiplicity 1, we   consider only    words    which contain  no letters of multiplicity 1. We call such words {\it multiplicity-one-free}.  We first summarize the results   concerning the words of length $\leq 4$.

\begin{theor}\label{erf4444sid}  A  multiplicity-one-free word  of length $  \leq 4$ in the alphabet $\alpha$   has one of the following   forms:
 $aa, a^3=aaa, a^4=aaaa, aabb, abba, bbaa, abab$ with distinct $a, b\in \alpha$. The words $aa, aabb, abba, bbaa$ are   contractible. The words $a^3, a^4$ are contractible if and only if $\tau (a)=a$. The word  $abab$ is contractible if and only if $\tau (a)=b$. Non-contractible words of type $a^3, a^4, abab$ are homotopic if and only if they are equal, i.e., coincide letter-wise.
   \end{theor}
  \begin{proof} That all multiplicity-one-free words  of length $  \leq 4$ belong to our list is obvious.
The words $aa$, $aabb$, $abba$, $bbaa$ are nanowords and can be   contracted by the first homotopy move. 
If $\tau (a)=a$, then   $a^3, a^4$ are contractible by Example \ref{exithhhg}.3.
If $\tau (a)\neq a$, then   $a^3, a^4$ are not contractible by   Sect.\  \ref{wor2597unkn1}.
If $\tau (a)=b\neq a$, then   $abab$ is  contractible by the second homotopy move.
If $\tau (a)\neq b\neq a$, then    $abab$ is  not contractible by Theorem \ref{classid}. The same theorem shows that two non-contractible homotopic words of this type are equal. That   non-contractible homotopic   words of type $a^3$ or $a^4$ are equal  follows from Theorem \ref{lipidyd}.
That a word $abab$ with $\tau (a)\neq b \neq a$ is not   homotopic to a monoliteral word follows from
the results of Sect.\ \ref{wor2597unkn1}. 
 \end{proof}

\begin{theor}\label{erfclassid}  A  multiplicity-one-free word  of length $  5$ in the alphabet $\alpha$ 
has one of the following   forms:
\begin{equation}\label{1woki} aaabb,\, aabba,\, abbaa,\, bbaaa,\end{equation}
\begin{equation}\label{woki} a^5,\, abaab,\, baaba,\, aabab,\, ba baa,\, baaab,\, ababa\end{equation}
   with distinct $a, b\in \alpha$.   The words    (\ref{1woki}) are homotopic to $aaa$; they are contractible if and only if $\tau (a)=a$.  The words  $a^5$ and $   baaab $ are contractible if and only if $\tau (a)=a$.  The word  $  ababa$   is contractible if and only if $\tau (a)=b$.     The words  $  abaab$, $baaba$,   
 $aabab$, $babaa$    are  not contractible.  Two non-contractible words of type   (\ref{woki}) are  homotopic   if and only if they coincide letter-wise with the following   exception:   $$aabab\simeq ba baa \simeq baaab$$ for   $a=\tau(b) \neq b$.  Non-contractible words  of type   (\ref{woki})  are not  homotopic to   words of length $\leq 4$.   \end{theor}

Note a subtlety in the last claim: it does not exclude   that  a word of type  (\ref{woki})  is homotopic to a {\it nanoword}  over $\alpha$ of length $\leq 4$ (cf. Example \ref{exithhhg}.5). This nanoword however can not be a word in the alphabet $\alpha$.  
 
 Among the words listed in Theorem \ref{erfclassid} only the words $a^5, ababa, baaab$ are symmetric. Here is a complete list of homotopically symmetric words of length $ 5$.

\begin{corol}\label{c49dpmid}   A  multiplicity-one-free word of length $ 5$ is homotopically  symmetric if and only if it is   contractible or symmetric or belongs to the list (\ref{1woki}) or has   the form  $ aab a b$ or  $ babaa $ with $b=\tau(a)\neq a$.        \end{corol}

 Among the words listed in Theorem \ref{erfclassid} only the words $a^5, ababa, baaab$ may be  skew-symmetric and this happens  when $a=\tau(a)$ and $b=\tau (b)$.  
  
\begin{corol}\label{cd444963pmid}     A  multiplicity-one-free word of length $ 5$ is homotopically   skew-symmetric if and only if it is contractible or  skew-symmetric.      \end{corol}

\begin{corol}\label{1corrrolrdsid}  Let $\tau=\id:\alpha\to \alpha$. A
 multiplicity-one-free word  of length $\leq 5$ in the alphabet $\alpha$   is  non-contractible if and only if it has one of the following six forms:
 $$abab, \,  abaab,  \,  baaba,    \, 
 aabab,   \,   ba baa,   \,  ababa $$
 where $a\neq b$.
Such  words   are  homotopic   if and only if they coincide letter-wise.        \end{corol}

The rest of this section is devoted to a  proof of Theorem \ref{erfclassid}. 

 \subsection{Reduction to a lemma}\label{easyfunct}  That all multiplicity-one-free words  of length $   5$ belong to our list is obvious. The words (\ref{1woki}) are homotopic to $aaa$ by the first homotopy move.
 By Example \ref{exithhhg}.3 and Theorem \ref{lipidyd}, the words $a^3$ and   $a^5$ are  contractible if and only if $\tau(a)=a$. By Sect.\ \ref{wor2597unkn1},  if $a^5$ is non-contractible, then it is not homotopic to a word of length $\leq 4$.   In the rest of the proof we  use the notation $ w^1_{a,b}=abaab$ and similarly 
$$  w^2_{a,b}=baaba,  \,
w^3_{a,b}=aabab, \,w^4_{a,b}=ba baa, \,  w^5_{a,b}=baaab,\,   w^6_{a,b}=ababa.  $$ 
 where $a\neq b$. Using Example \ref{exithhhg}.3, it is easy to show that  the word $w^5_{a,b}$ with $ \tau (a)=a$ is contractible.  In the sequel, writing $w^5_{a,b}$ we always suppose that $ \tau (a)\neq a$. By   Example \ref{exithhhg}.4, the word $w^6_{a,b}$ with $a=\tau (b)$ is contractible. In the sequel, writing $w^6_{a,b}$ we always suppose that $a\neq \tau (b)$. 
 
It remains  to prove the following lemma for  $i=1,..., 6$.

\begin{lemma}\label{roki}    ($\ast$) The  word $w^i_{a,b }$ is not  homotopic to a  word of length $\leq 4$ or to a monoliteral word.  In particular,  $w^i_{a,b }$ is non-contractible.
 
  ($\ast \ast$) The  words    $w^i_{a,b}, w^i_{a',b'}$   are homotopic if and only if $a=a'$ and   $b=b'$.  
  
   ($\ast \ast \ast$) For $j<i$, the  words of type $w^i$   are not homotopic to the  words  of type  $w^j$    with the following   exception: $w^3_{a,b}\simeq w^4_{a,b} \simeq w^5_{a,b}$ for $a=\tau(b) \neq b$.   \end{lemma}
   
   We   prove this lemma case by case.      The cases $i=3$ and  $i=6$  need    additional techniques at least for some $a,b$.  The proof   in this cases    will be accomplished  in Sect.\ \ref{dd822}.

   \subsection{Case $i=1$}\label{i=1} The \'etale word associated with $  w^1_{a,b}=abaab$ is the pair $(\A=\alpha, aabab)$. Its desingularization  yields the nanoword $w=A_3A_2BA_3A_1A_2A_1B$ in the $\alpha$-alphabet 
\begin{equation}\label{bet} \{A_1=a_{2,3}, A_2=a_{1,3}, A_3=a_{1,2}, B=b_{1,2}\}\end{equation} where  $\vert A_1\vert=\vert A_2\vert=\vert A_3\vert =a, \vert B \vert=b$.  The   linking $\alpha$-form  of $w$  is given by  
 $$ n_w  = \left [ \begin{array}{ccccc}   
        0& -1 & 0 & 0 \\
 1& 0 & -1 & 1\\
0& 1 & 0 & 1\\
0& -1 & -1 &  0\\
\end{array} \right ] , \,\,\,\, lk_w= \left [ \begin{array}{ccccc}   
        1& 1 & a^{-1} & 1 \\
        1& 1 & 1 & a \\
 a& 1 & 1 & a\\
1 & a^{-1} &  a^{-1} & 1\\
\end{array} \right ] $$
where the   rows  and columns  correspond to  $A_1,A_2,A_3,B$, respectively.   Consider the associated   $\alpha$-pairing  $p^w= (S,s, S\times S \to \alpha)$. Here   $S=\{s, A_1,A_2,A_3,B\}$ and $S\times S \to \alpha$ is given by the matrix 
$$   \left [ \begin{array}{ccccc}   
        1&  a & b^{-1}  &  a^{-1} b^{-1} & a^2 \\
        a^{-1}& 1 & a^{-2}  &  a^{-2} & 1 \\
        b & a^{2}  & 1  &  a^{-2} & a^3 b \\
        ab & a^{2} & a^2  &  1 & a^3 b \\
a^{-2} & 1 & a^{-3}b^{-1}  &  a^{-3}b^{-1} & 1 \\
\end{array} \right ]  $$
where the   rows  and columns  correspond to   $s, A_1,A_2,A_3,B$, respectively.
An easy inspection (using only that $a\neq b$) shows that  $p^w$ is primitive and therefore $\vert \vert w^1_{a,b} \vert \vert=\vert \vert w\vert \vert =4$.  In particular $w^1_{a,b}$ is  non-contractible. We can recover $a$ and $b$ from $p^w$ as follows. Since $p^w$ is primitive, it is determined by the homotopy class of $w^1_{a,b}$ uniquely up to isomorphism. The projection $S-\{s\}\to \alpha$ takes one value (the letter $a$) with multiplicity 3 and another value (the letter $b$) with multiplicity 1. This \lq\lq multiplicity argument" shows that   if  $w^1_{a,b}\simeq w^1_{a',b'}$, then  $a=a', b=b'$. 
The equality $\vert \vert w^1_{a,b} \vert \vert =4$ implies that $w^1_{a,b}$ is not homotopic to a  nanoword of length $\leq 4$. 
That $w$ is not homotopic to a monoliteral word follows from the fact that the projection $S-\{s\}\to \alpha$ takes two distinct values while the similar projection associated with a monoliteral word takes only  one value.

 \subsection{Case $i=2$}\label{i=2} Since $w^2_{a,b}$ is   opposite to $w^1_{a,b} $, we have $\vert \vert w^2_{a,b} \vert \vert=\vert \vert w^1_{a,b} \vert \vert =4$ and the claims ($\ast  $),  ($\ast \ast$)   for $i=1$ imply similar  claims for $i=2$.   We need to check only    that   $w^2_{a,b}$ is not  homotopic to  $w^1_{a',b'}$. As we know, the $\alpha$-pairing of  
 $w^2_{a,b}=(w^1_{a,b})^-$ is opposite to the one of $w^1_{a,b}$.
If $w^2_{a,b}\simeq w^1_{a',b'}$,  then the multiplicity argument counting   the values of the projection  $S-\{s\}\to \alpha$ shows that $a=a'$ and $b=b'$. The matrix of the   $\alpha$-pairing   of $w^2_{a,b} $ is obtained by transposition from  the one of   $w^1_{a,b} $ (and possibly by  simultaneous permutation of rows and columns keeping the first row and the first column). Comparing the values in the first column  we obtain that  $\{a^{-1},b\}=\{a,b^{-1}\}$. Since $a\neq b$, we   have  $\tau(a)=a, \tau(b)=b $. Let  $H=\{1 ,b\} \subset \pi$. The $H$-coverings of  $w^1_{a,b} \simeq A_3A_2BA_3A_1A_2A_1B$ and $w^2_{a,b}   \simeq B A_3 A_2 A_3 A_1 B A_2 A_1 $ are  the nanowords $ A_2B A_2 B$ and $  B A_2 B A_2 $, respectively. Since $a\neq b$, these nanowords  are homotopic only if $a=\tau(b)$
which is excluded by $\tau(b)=b$. 

 \subsection{Case $i=3$}\label{i=3}  The desingularization of the \'etale word associated with $  w^3_{a,b}=aabab$    yields the nanoword $A_3A_2A_3A_1BA_2A_1B$ in the   $\alpha$-alphabet (\ref{bet}).
Its   linking $\alpha$-form    is given by  
 $$ n   = \left [ \begin{array}{ccccc}   
        0&  -1 & 0 & 1 \\
  1& 0 & -1 & 1\\
0&  1 & 0 & 0\\
-1& -1 & 0 &  0\\
\end{array} \right ] , \,\,\,\, lk = \left [ \begin{array}{ccccc}   
        1& 1 & a^{-1} & 1 \\
        1& 1 & 1 & a \\
 a& 1 & 1 & a\\
1 & a^{-1} &  a^{-1} & 1\\
\end{array} \right ] . $$
The associated   $\alpha$-pairing    is given by the matrix 
$$  p= \left [ \begin{array}{ccccc}   
        1&  ab^{-1} & b^{-1}   &  a^{-1}  & a^{2} \\
        a^{-1}b & 1 & a^{-2}  &  a^{-2} & ab \\
        b  & a^{2}  & 1  &  a^{-2} & a^3 b \\
        a  & a^{2} & a^{2}  &  1 & a^2 \\
a^{-2} & a^{-1}b^{-1} & a^{-3}b^{-1}  &  a^{-2} & 1 \\
\end{array} \right ]  $$
where the   rows  and columns  correspond to   $s, A_1,A_2,A_3,B$, respectively.

We first study the case where  $a\neq \tau(b)$ (and  $a\neq b$). Then $ab\neq 1$ and $p $   is primitive. The arguments as above show that:   $\vert \vert w^3_{a,b} \vert \vert  =4$;  $w^3_{a,b}$ is    not homotopic to a  nanoword of length $\leq 4$ or to a monoliteral word; if $w^3_{a,b} \simeq w^3_{a',b'}$ with $a'\notin \{b', \tau (b')\}$, then     $a=a', b=b'$. 
 If  $w^3_{a,b} \simeq w^1_{a',b'}$ with $a'\neq b' $, then  the multiplicity argument gives    $a=a', b=b'$ and  comparing the values in the $s$-column we obtain that  $\{a^{-1},b\}=\{a,b \}$. Hence   $a=\tau(a)\neq b$.  We must prove the following. 

 \begin{lemma}\label{eril} If $a=\tau(a)\neq b $, then  the nanowords $A_3A_2A_3A_1BA_2A_1B$  and 
 $A_3A_2BA_3A_1A_2A_1B$  in the   $\alpha$-alphabet (\ref{bet}) are not homotopic.\end{lemma} 
    
The   $\alpha$-pairings of  the  nanowords  in question are isomorphic and   coverings do not help.  We prove  this lemma in Sect.\ \ref{erilpotnct} using other techniques. 

  If  $w^3_{a,b} \simeq w^2_{a',b'}$  with $a'\neq b' $, then   as above   $a=a', b=b'$ and   $\{a,b^{-1}\}=\{a,b \}$, $\{a^{-1}b^{-1}, a^2\}=\{a^{-1}b , a^{-2} \}$. Hence $b=\tau(b)\neq a=\tau (a)$. Let  $H=\{1 ,b\} \subset \pi$. The $H$-coverings of  $w^3_{a,b} $ and $w^2_{a,b} $ are  the nanowords $ A_2B A_2 B$ and $  B A_2 B A_2 $. If  they are homotopic, then $a =b$
which contradicts the assumptions.

We now study the case $a=\tau (b) $. By Example \ref{exithhhg}.5,  $w^3_{a,b}$ is  homotopic to   the nanoword $w_{a,a}$ of length $4$. The latter  is non-contractible since $a\neq b= \tau(a)$. Therefore $\vert \vert w^3_{a,b} \vert \vert =2$. The  word $w^3_{a,b}\simeq w_{a,a}$  is not homotopic to a word $cdcd$ with distinct $c,d \in \alpha$ by Theorem \ref{classid}. The word $w^3_{a,b} $ is not homotopic to  a monoliteral word $c^m $ with $c\in \alpha, c\neq \tau (c), m\geq 3$  because  $\vert \vert c^m\vert \vert \geq 3$ by  (\ref{ine})   and  (\ref{normccc}).
The word $w^3_{a,b} $ is not homotopic to $ w^i_{a',b'}$ for $i=1, 2$   since they have different norms. 
If   $w^3_{a,b}\simeq w^3_{a',b'}$, then $\vert \vert w^3_{a',b'} \vert \vert=\vert \vert w^3_{a,b} \vert \vert=2$ and   we must have $a'=\tau(b')$. Then $w_{a,a}\simeq w_{a',a'}$ and therefore
$a=a'$ and $b=\tau(a)=\tau(a')=b'$. 

\subsection{Case $i=4$}\label{i=4}  Since $w^4_{a,b}=(w^3_{a,b})^-$,  the claims ($\ast  $),  ($\ast \ast$)   for $i=4$ follow from the  similar  claims for $i=3$. The claim  ($\ast   \ast \ast$)   for $i=3$ implies that  $w^4_{a,b}$ is not homotopic to  $w^j_{a',b'}$ with $j=1,2$. If  $a=\tau (b) $, then  $w^4_{a,b}\simeq w^3_{a,b}$ by Example \ref{exithhhg}.5.  It remains  to prove that
 $w^4_{a,b}$ with $a\neq \tau(b)$ (and  $a\neq b$) is not homotopic to  $w^3_{a',b'}$. 
If they  are homotopic then the multiplicity argument yields  $a=a'$ and $b=b'$. The matrix of the   $\alpha$-pairing   of $w^4_{a,b}=(w^3_{a,b})^- $ is obtained by transposition from  the one of   $w^3_{a,b} $. Comparing the values in the $s$-column   we obtain that   $\{a^{-1},b^{-1}\}=\{a,b\}$. Since $a\neq \tau(b)=b^{-1}$, we have  $a=a^{-1}=\tau(a)$ and $b= \tau(b) $. Let  $H=\{1,b\}\subset \pi$. The nanowords   $(w^3_{a,b})^H= A_2B A_2 B $ and $(w^4_{a,b})^H= B A_2 B A_2 $     are homotopic only if $a=\tau(b)$ or $a=b$
which is excluded by the assumptions. 

\subsection{Case $i=5$}\label{i=5}  
The desingularization of the \'etale word associated with $  w^5_{a,b}=baaab$    yields the nanoword $BA_3A_2A_3A_1 A_2A_1B$  in the   $\alpha$-alphabet (\ref{bet}).
Its   linking $\alpha$-form   is given by  
 $$ n   = \left [ \begin{array}{ccccc}   
        0&  -1 & 0 & 0 \\
  1& 0 & -1 & 0\\
0&  1 & 0 & 0\\
0& 0 & 0 &  0\\
\end{array} \right ] , \,\,\,\, lk = \left [ \begin{array}{ccccc}   
        1& 1 & a^{-1} & a^{-1} \\
        1& 1 & 1 & 1 \\
 a& 1 & 1 & a\\
a & 1 &  a^{-1} & 1\\
\end{array} \right ] .$$
The associated   $\alpha$-pairing   is given by the matrix 
$$  p= \left [ \begin{array}{ccccc}   
        1&  a & 1  &  a^{-1}  & 1 \\
        a^{-1} & 1 & a^{-2}  &  a^{-2} & a^{-2} \\
       1  & a^{2}  & 1  &  a^{-2} & 1 \\
        a  & a^{2} & a^{2}  &  1 & a^2 \\
1 & a^{2} & 1  &  a^{-2} & 1 \\
\end{array} \right ]  $$
where the   rows  and columns  correspond to   $s, A_1,A_2,A_3,B$, respectively.
By our assumptions $\tau(a)\neq a \neq b$.   
Assume first that   $a\neq \tau(b)$. Then   $p $   is primitive and the arguments as above show that:   $\vert \vert w^5_{a,b} \vert \vert  =4$;  $w^5_{a,b}$ is    not homotopic to a  nanoword of length $\leq 4$ or to a monoliteral word; if $w^5_{a,b} \simeq w^5_{a',b'}$ with $a'\notin \{b', \tau (b')\}$, then     $a=a', b=b'$. 
 If  $w^5_{a,b} \simeq w^j_{a',b'}$  with $j=1,2,3,4$,  then comparing the $s$-columns of  the corresponding $\alpha$-pairings
we always obtain that $b'=a'$, a contradiction.

We now study the case $a=\tau (b)\neq b$. By Example \ref{exithhhg}.5,  $w^5_{a,b}\simeq  w_{a,a}$. As in the case $i=3$, we obtain  that   $\vert \vert w^5_{a,b} \vert \vert =2$ and  $w^5_{a,b} $ verifies the claim ($\ast$) of the lemma. 
If  $w^5_{a,b} \simeq  w^j_{a',b'}$ for $j\leq 5$, then comparing the norms we observe that $j\in \{3,4, 5\}$ and $a'=\tau (b')$. Then $w_{a,a}\simeq w^5_{a,b} \simeq  w^j_{a',b'} \simeq w_{a',a'}$ so that  $a=a'$ and $ b=\tau (a)=\tau(a')=b'$. Thus we recover either one of the exceptional homotopies  from   ($\ast \ast \ast$) or the tautological homotopy 
   $w^5_{a,b}\simeq w^5_{a,b}$.

 \subsection{Case $i=6$}\label{i=6}  The desingularization of the \'etale word associated with $  w^6_{a,b}=ababa$    yields the nanoword $w= A_3A_2BA_3A_1 BA_2A_1 $  in the   $\alpha$-alphabet (\ref{bet}).  The   linking $\alpha$-form  of $w$  is given by  
 $$ n_w  = \left [ \begin{array}{ccccc}   
        0&  -1 & 0 & -1 \\
  1& 0 & -1 & 0\\
0&  1 & 0 & 1\\
1& 0 & -1 &  0\\
\end{array} \right ] , \,\,\,\, lk_w= \left [ \begin{array}{ccccc}   
        1& b^{-1}  & a^{-1}  b^{-1}  & 1 \\
        b& 1 & b^{-1}  & 1 \\
 ab& b & 1 & 1\\
1 & 1 &  1 & 1\\
\end{array} \right ] .$$
The associated   $\alpha$-pairing      is given by the matrix 
\begin{equation}\label{scol} p^w=  \left [ \begin{array}{ccccc}   
        1&  ab & 1   &  a^{-1}   b^{-1}  & 1 \\
        a^{-1}  b^{-1}  & 1 & a^{-2} b^{-2}   &  a^{-2} b^{-2}   & a^{-1} b^{-1}   \\
        1  & a^{2} b^{2}    & 1  &  a^{-2} b^{-2}   & 1\\
        a b & a^{ 2}b^2 & a^{2}b^2  &  1 & a  b \\
1 & ab & 1  &  a^{-1}b^{-1} & 1 \\
\end{array} \right ]  \end{equation}
where the   rows  and columns  correspond to   $s, A_1,A_2,A_3,B$, respectively. This matrix depends only on   $ab\in \pi$. The assumption $a\neq \tau (b)$ guarantees that $ab\neq 1$.  It is clear that $A_1,A_3,B$ are not annihilating and $A_2$ is annihilating if and only if $a^2b^2=1$.  
The latter is equivalent to  $a=\tau(a)$ and $ b=\tau(b)$.  The  only possible   twins are $A_1,A_3$ and again this happens     if and only if  $a=\tau(a), b=\tau(b)$.
Thus, if $a=\tau(a),  b=\tau(b)$, then  $p^w$ is homologous to a trivial $\alpha$-pairing.  Since this is not the case for   the   words considered above in this lemma (i.e.,   the words $w^j , j=1,...,5$,  the  non-contractible monoliteral words, and  the  non-contractible nanowords of length 4),   $w$ is not homotopic to either of them. It remains to prove the following.
 
  \begin{lemma}\label{6eril} Let  $a=\tau(a)\neq b=\tau(b)$. Then  $w^6_{a,b}$ is non-contractible. If 
 $w^6_{a,b} \simeq w^6_{a',b'}$ with   $a'=\tau(a')\neq b'=\tau(b')$, then $a=a'$ and $b=b'$. \end{lemma}

 We shall prove this lemma in Sect.\ \ref{e666lpotnct}.

Consider now the case where $a\neq \tau(a)$ or $b\neq \tau(b)$. Then $p^w$ is primitive and  $\vert\vert w\vert \vert  =4$. The   arguments given in the case $i=1$ apply here and yield the claims ($  \ast$), ($\ast \ast$) for $i=6$. That $w^6_{a,b}$ is not homotopic  to $w^j_{a',b'}$ with $j\leq 5$ follows from the fact that  their $\alpha$-pairings are not   isomorphic:  the $s$-column of  (\ref{scol}) does not contain generators $\{e\}_{e\in \alpha} \subset \pi$ while the $s$-column  of the $\alpha$-pairings associated with $w^j$    contains such generators.

\section{Colorings of nanowords}\label{col56822}

 Throughout this section we fix a set $\beta\subset \alpha$ such that $\tau(\beta)=\beta$. (The set $\beta$ may be empty or equal to $\alpha$.) 
We shall use notation introduced in Sect.\ \ref{nota1}.
 
 \subsection{Tricolorings}\label{tricolgroupga}  A {\it
tricoloring}  (with respect to $\beta\subset \alpha$) of  a nanoword
$(\A, w:\hat n
\to
\A)$ over
$\alpha$ 
   is a function $f:\{0,1,...,n\}\to \ZZ/3\ZZ$ such that   
for  any  $A\in \A$ with $\vert A\vert \in \beta$,    $$f(i_A)=f(i_A-1) ,\,\,\,\, f(j_A-1)+   f(j_A) +  f(i_A)=0$$
and  for   any $A\in \A$ with $\vert A\vert \in \alpha-\beta$,    $$f(i_A-1) + f(i_A)   +   f(j_A)=0 ,\,\,\,\, f(j_A)=  f(j_A-1).$$  Note that a sum of three residues
$(\modu 3)$ is $0$ if  and only if   these residues  are  either  equal to each other  or form the triple $\{0,1,2\}$.   It is convenient   to insert dashes
between consecutive letters of $w$ and also at the beginning and the end of $w$.  Thus instead of   $w(1) w(2)\cdots w(n)$ we  
write
$-w(1) -w(2)-\cdots  - w(n)-$.  The  residue $f(i)$ can be interpreted as the label of the dash between $w(i)$ and $w(i+1)$. The residues $f(0)$ and $f(n)$ are
 the labels of the leftmost and the rightmost dashes, respectively.  They are called the {\it input} and the {\it output} of    $f$, respectively.

The sum of two tricolorings of  a nanoword $w$ is a tricoloring of $w$. Thus  the tricolorings of $w$ form a   vector space over $\ZZ/3\ZZ$.  It is non-zero since  $w$ has     constant  tricolorings   $f=const$. Let $c_{k,l}=c (\beta,w)_{k,l}$ be the number   of tricolorings of $w$ with given input  and output $k,  l\in \ZZ/3\ZZ$. The difference between two tricolorings with the same input and output has 
input and output zero. Therefore  either  $c_{k,l}=0$ or $c_{k,l}=c_{0,0}$. Adding
constant   tricolorings and negating,  we   observe that     $c_{k,l}$ depends
only on
$k-l$ and  $c_{k,l}=c_{-k,-l}$.  We conclude that  the matrix $(c_{k,l})$ has the form
\begin{equation}\label{mattt}  \left [ \begin{array}{ccccc}   
        c & c  & c  \\
c & c  & c \\
c & c  & c \\
\end{array} \right ] \,\,\,\, {\rm {or}} \,\,\,\, \left [ \begin{array}{ccccc}   
        c & 0 & 0 \\
0& c  & 0\\
0& 0 & c \\
\end{array} \right ]\end{equation}
  where   $c =c_{0,0}\geq 1$ is a   power of $3$.

   An easy check shows that  the number  $c (\beta,w)_{k,l}$  is
preserved under     homotopy moves on
$w$ for all $k,l$, cf. Sect.\ \ref{colgroupga}.  This   
  can be   used to distinguish homotopy types of nanowords (see examples below). In particular,  if a
nanoword admits   non-constant tricolorings, then it is non-contractible.  
 
  The formula
$w\mapsto   (c (\beta,w)_{k,l})_{k,l\in \ZZ/3\ZZ}$ defines a representation of the monoid of homotopy classes of nanowords $\N_\bullet (\alpha)$ in the monoid of integral  $(3\times 3)$-matrices of type 
(\ref{mattt}).  For $\beta=\emptyset$ or
$\beta=\alpha$,   any     nanoword  has only constant tricolorings.

\subsection{Examples}\label{rexixixinct}   1. Consider  the nanoword $w=ABAB$ such that 
$\vert A\vert  ,  \vert B\vert  $   lie in different orbits of $\alpha$.  Take as $\beta$   the orbit of $\vert A\vert$. 
Then $w$ admits a  non-constant tricoloring 
$$ \stackrel{0}{-}A\stackrel{0}{-}B\stackrel{2}{-}A\stackrel{1}{-}B\stackrel{1}{-}.$$
 This gives another proof of the fact that this nanoword is not contractible.
   
   2. Consider the nanoword 
$w=A_1A_2BA_3A_1BA_2A_3$ where $\vert A_1\vert =\vert A_2\vert=\vert A_3\vert$ and $\vert B\vert$  lie in different orbits of $\tau$. 
 Take as $\beta$   the orbit of $\vert A_1\vert$. 
An easy check shows that $c (\beta,w)_{k,l}=3$ for all $k,l$.  Here is an example of a coloring of $w$ with distinct input and output:
$$\stackrel{0}{-}A_1\stackrel{0}{-}A_2\stackrel{0}{-}B\stackrel{1}{-}A_3\stackrel{1}{-}A_1\stackrel{2}{-}B\stackrel{2}{-}A_2\stackrel{1}{-}A_3\stackrel{1}{-}.$$
Hence $w$ is non-contractible. Note that  $w\in {\rm {Ker}} (\tilde
\gamma)$.

\subsection{Ring  $\Lambda$}\label{ringsfunct} To define more general colorings, we introduce 
  a group $
\Psi_\alpha$  generated by  
$\{{}a, a_{\bullet}  
\}_{a\in
\alpha}$ with    defining relations
$ {}a a_{\bullet}=a_{\bullet} {}a$ and $ a\, \tau (a)=a_{\bullet}\, \tau(a)_{\bullet}=1   $ for   $a\in \alpha$.   
Let  $\Lambda=\ZZ
\Psi_\alpha$   be  the integral group ring of  
$
\Psi_\alpha$.  We   view elements of $\Lambda$ as
non-commutative   polynomials in the variables  $\{{}a, a_{\bullet}  
\}_{a\in
\alpha}$ subject  to the relations $ {}a a_{\bullet}=a_{\bullet} {}a$ and $ a\,  \tau (a)=a_{\bullet} \, \tau(a)_{\bullet}  =1 $ for  all  $a\in \alpha$. 

\subsection{Generalized colorings}\label{colgroupga}    Let $X$ be a left $\Lambda$-module.   An {\it
$X$-coloring} (with respect to $\beta\subset \alpha$) of a nanoword
$(\A, w:\hat n
\to
\A)$    is a function $f:\{0,1,...,n\}\to X$ such that   
for  any $A\in \A$ with $\vert A\vert =a\in \beta$,    $$f(i_A)=a f(i_A-1) ,\,\,\,\, f(j_A)=a_\bullet    f(j_A-1) + (1- a a_\bullet) f(i_A-1)$$
and  for  any $A\in \A$ with $\vert A\vert=a \in \alpha-\beta$,    $$f(i_A)=a_\bullet  f(i_A-1)   + (1- a a_\bullet  ) f(j_A-1) ,\,\,\,\, f(j_A)= a f(j_A-1).$$  
For $X= \ZZ/3\ZZ$, $ax=x, a_\bullet x=-x$ for all $a\in \alpha, x\in X$, this gives
tricolorings of $w$.

   The sum or the difference  of two $X$-colorings of   $w$ is always an $X$-coloring of $w$.
 The zero function $f=0$ is an $X$-coloring of $w$.    Generally speaking, non-zero  constant functions are not $X$-colorings unless we make additional assumptions like for instance   $ax=x$ for all $a\in \alpha, x\in X$.  
If $X$ is finite as a set, then  there is only a finite number  of
$X$-colorings $f$  of
$w$ with given input 
$k=f(0)\in X$ and output  $  l=f(n)\in X$. Denote this number $X(\beta,w)_{k,l}$. As in Sect.\ \ref{tricolgroupga},  either $X(\beta,w)_{k,l}=0$ or $X(\beta,w)_{k,l}=X(\beta,w)_{0,0}$.  An
easy check shows that  $X(\beta,w)_{k,l}$   is preserved under   homotopy moves on
$w$ (we shall verify this    in Sect.\ \ref{3matrinvaricovfunct}).  The formula
$w\mapsto   (X(\beta,w)_{k,l})_{k,l\in X}$ defines a representation of   $\N_\bullet (\alpha)$  in the monoid of integral  $(\card (X) \times \card (X))$-matrices 
  with   non-negative entries.
    
A vast family of $\Lambda$-modules can be obtained by the following construction.
Pick a unital associative ring $R$   and   two functions
$p, p_\bullet:\alpha\to R $ such  that $p(a) \, p(\tau (a))=p_\bullet(a) \, p_\bullet(\tau (a))=1  $ for all $a\in \alpha$.  Any left $R$-module $X$ becomes a
$\Lambda$-module by $ax= p(a) x$ and $a_\bullet x =p_\bullet (a) x$ for $a\in \alpha, x\in X$.    

As an application, consider the nanoword $ABAB$ over $\alpha=\{a, \tau (a), b, \tau (b)\}$ with $\vert A\vert =a, \vert B\vert =b$. We check that this nanoword is
not homotopic to its inverse. Set $X=R=\ZZ/5\ZZ$, $\beta=\{a, \tau (a)\}$, $p=1$ (the constant function) and let $  p_\bullet:\alpha\to \ZZ/5\ZZ $ be given by   
$p_\bullet(a)=p_\bullet(\tau(b))=2$ and
$  p_\bullet(\tau (a))=p_\bullet(b)=3$. The nanoword
$ABAB$ admits a non-constant $X$-coloring
$$\stackrel{0}{-}A\stackrel{0}{-}B\stackrel{3}{-}A\stackrel{1}{-}B\stackrel{1}{-} $$
while the inverse nanoword $A'B'A'B'$   with $\vert A'\vert =\tau(a), \vert B'\vert =\tau(b)$ admits only constant $X$-colorings.

\subsection{Exercise}\label{excicicnct} Prove the identity $X(\beta,{\overline w}^-)_{k,l}=X(\alpha-\beta,w)_{l,k}$ for any finite $\Lambda$-module $X$ and $k,l\in X$.

\section{Matrices, modules, and polynomials}\label{polyorticcdfd}

Throughout this section we fix a set $\beta\subset \alpha$ such that $\tau(\beta)=\beta$.  We  associate with each   nanoword $w$ a certain module $K_\beta(w)$.  We begin with     so-called marked modules and weighted matrices.

\subsection{Marked modules and weighted matrices}\label{matrinvaricovfunct}  By a {\it marked module} over a unital associative ring $R$, we   mean 
a  left $R$-module $K$ endowed with an ordered pair of vectors $v_-,v_+\in K$ called
 the {\it input} and the {\it output}, respectively. By {\it isomorphisms} of marked modules we mean isomorphisms of modules preserving the input and the output. 

  Let  $R^{ab}$ be the  commutative ring obtained by quotienting $R$ by the 2-sided ideal
generated by the set $\{rs-sr\,\vert\, r,s \in R\}$. Denote the projection $R\to R^{ab}$ by $q$. 

  A {\it weighted matrix} over
$R$ is   a   finite  matrix over
$R$  endowed with an ordered pair of elements of $  R^{ab}$  called   the  {\it weights}.  
  Consider the following   transformations on a weighted   matrix $M$  with $n\geq 1$ rows and $m\geq 1 $ columns.

(i)$'$ Permutation of two consecutive rows and simultaneous multiplication of both weights  by $-1$.

(ii)$'$ Addition  to the $k$-th column of  the $j$-th column multiplied on the right by an  element of $R$ where   $1\leq k\leq m $, $2\leq  j\leq m-1$,  and $k\neq
j $.  The weights are preserved. 

(iii)$'$  Striking out the $i$-th row and the $j$-th column provided $n\geq 2$, $ 1\leq i \leq n$,  $2\leq j
\leq m-1$ and  the $(i,j)$-th entry  $M_{i,j}$ of $M$ is invertible in  $R$ and  is the only non-zero entry in
the
$i$-th row. The first weight is multiplied by
$(-1)^{i+j+1} q(M_{i,j})$ and the second weight is multiplied by
$(-1)^{i+j} q(M_{i,j})$.

We say that two weighted matrices (possibly of different sizes) are {\it equivalent} if  they can be obtained from each other by a finite sequence of transformations (i)$'$
-- (iii)$'$ and the inverse transformations.   The number $m-n$, called the {\it deficiency} of
$M$, is preserved under the transformations (i)$'$ --
(iii)$'$.   For a weighted matrix $M$,  denote by
$\div M$ the    weighted matrix obtained by multiplying both weights of $M$ by
$  -1$.  

Each weighted $(n\times m)$-matrix $M$ presents an $R$-module denoted $K(M)$. More precisely, $M$ determines an $R$-homomorphism $R^n\to R^{m}$
where the canonical basis vectors  of $R^{m}$ (resp.\ of
$R^n$) correspond to the columns (resp.\ the rows) of $M$.   The cokernel of this homomorphism is   $K(M)$. The first and the last basis
vectors  of $R^{m}$ project to certain  $v_-, v_+ \in K(M)$.  This gives a marked module  $(K(M), v_-,v_+)$. It does not depend on the weights of $M$
and  is preserved  (up to isomorphism) under the transformations (i)$'$ -- (iii)$'$ on $M$ and the involution $M\mapsto \div M$.

 A weighted matrix
$ (M,
r_-\in
R^{ab}, r_+\in R^{ab})$ of deficiency 1    gives  rise to two elements of   $R^{ab}$ as
follows.     Consider  the square 
matrix  $M^{ab}_-$ (resp.\ $M^{ab}_+$)  obtained from
$M$     by deleting  the first (resp.\ the last) column and projecting all entries  to $R^{ab}$. For $\varepsilon\in \{+,-\}$, set  
 $\nabla^\varepsilon (M  )=  r_\varepsilon  \det (M^{ab}_\varepsilon)\in R^{ab}$.
It is obvious that $\nabla^\varepsilon(M)$ is  preserved under the transformations (i)$'$ -- (iii)$'$ and   $ \nabla^\varepsilon( \div M)=-\nabla^\varepsilon(M)$.   
 Note  that   $M_{ \varepsilon }^{ab}$ is a presentation matrix of the $R^{ab}$-module  $R^{ab}\otimes_R K(M)/Rv_{\varepsilon }$.   The principal    ideal  of   $R^{ab}$ generated by $\nabla^\varepsilon(M)$   coincides with 
the order ideal of this module.  (This  $R^{ab}$-module can be presented by a square matrix over $R^{ab}$ and the  determinant of any such matrix generates the order ideal.)

\subsection{Matrices of nanowords}\label{2matrinvaricovfunct}   Recall the ring $\Lambda=\ZZ \Psi_\alpha$ introduced in Sect.\ \ref{ringsfunct}.  With each non-empty nanoword $(\A, w:\hat
n\to
\A)$ over
$\alpha$ we associate a   $(2n\times (2n+1))
$-matrix 
$M  $ over   $\Lambda $ as follows.    Let us numerate all letters of $\A$ by the numbers $1,2,..., n$. The 
$k$-th letter
$ A_k$ determines the $(2k-1)$-st and $2k$-th  rows  of $M $. To define these rows, we   shall use notation of  Sect.\ \ref{nota1}. If $a=\vert A_k\vert \in \beta$, then
$$M_{2k-1, i_A-1}=a,  M_{2k-1, i_A} = -1,  M_{2k, i_A-1}=1-a a_\bullet,  M_{2k, j_A-1} = a_\bullet, M_{2k, j_A} = -1$$
while all other entries of these two rows are   $0$. If $a=\vert A\vert \in \alpha- \beta$, then one   exchanges $i_A, j_A$ in these formulas which gives
$$ M_{2k-1, j_A-1}=a,  M_{2k-1, j_A} = -1, M_{2k, i_A-1} = a_\bullet, M_{2k, i_A} = -1,  M_{2k, j_A-1}=1-a a_\bullet$$
while all other entries of these two  rows are  $0$.   We provide $M$ with weights   as follows.  The      ring $\Lambda^{ab}$ is the
ring of      polynomials in the commuting variables  $\{{}a, a_{\bullet}  
\}_{a\in
\alpha}$ subject   to the  relations $ a \tau (a)=a_{\bullet} \tau(a)_{\bullet}=1   $ for   $a\in \alpha$. 
 For   $a\in \alpha$, set  $\langle a,
w\rangle  =\card \{A\in \A\,\vert \,  \vert A\vert=a\}$ and 
$$r^w =\prod_{a\in \alpha- 
\beta} 
(-aa_\bullet)^{-\langle a,
w \rangle }  = \prod_{a\in
\alpha- \beta}
(-\tau(a) \tau(a)_\bullet)^{ \langle a,
w \rangle } \in \Lambda^{ab} .$$  The triple $ (M,  r^w, r^w)$ is a weighted matrix over $\Lambda$ denoted $M(w)$ or $M_\beta (w)$. Its 
equivalence class   does not depend on the choice of numeration of letters of $\A$ since the   permutation of the rows of $M$ correspondinding to two
different choices  is even. 
The next lemma describes the behavior of $M(w)$ under homotopy moves on $w$.

\begin{lemma}\label{nezzd}  The equivalence class of $M(w)$ is preserved under the first and second homotopy moves on $w$. If a nanoword $v$ is obtained from
$w$ by the third homotopy move, then $M(v) $ is equivalent to $\div {M(w)}$.
\end{lemma}
    \begin{proof}    It is obvious that isomorphic nanowords give
 the same weighted matrices.  
 Consider the first  homotopy move
$w=(\A,  x A A  y) 
\mapsto  (\A- \{A\},  x  y)=v$ with $v\neq \emptyset$.   We numerate the letters of $\A$ starting with $A$ and set $a=\vert A\vert \in \alpha $.  If
$a\in
\beta$, then 
$$M(w)=     
\Bigl ( \left [ \begin{array}{ccccc}  0& a&-1& 0& 0 \\
        0& 1- a a_\bullet & a_\bullet& -1& 0 \\
m_1& m_2& 0& m_3&m_4
\end{array} \right ], r^w , r^w \Bigr )$$
where $m_1, m_4$ are matrices and $m_2, m_3$ are columns over $\Lambda$, all with the same number of rows equal to twice the length of $v$. 
It is clear that $r^v=r^w $ and $M(v)=     
 ( [  
m_1\, \vert \,  m_2+ m_3\, \vert \,  m_4
  ], r^v, r^v)$ where the vertical bar    separates  submatrices. Multiplying the
middle column of $M(w)$  by $a\in \Lambda$ and adding it to the previous column we obtain a matrix in which the first row contains only one non-zero entry    $-1$. Applying
the transformation (iii)$'$ and using that $a a_\bullet=  a_\bullet  a $, we obtain the weighted matrix
$$  \Bigl (\left  [ \begin{array}{ccccc}   
        0& 1 & -1& 0 \\
m_1& m_2 & m_3&m_4
\end{array} \right ], (-1)^{k+1 } r^w, (-1)^{k } r^w \Bigr )$$
where $k$ is the number of columns in $m_1$, that is the length of  the word $x$.  
If $x\neq \emptyset$,  then $k\neq 0$ and we can apply the transformation (ii)$'$    adding  the  column containing $m_2$ to the next one. Applying  then  (iii)$'$ to
the first row   we obtain
$M(v) $. If $x= \emptyset$,  then  the assumption $v\neq \emptyset$ implies that $y\neq \emptyset$ so that the  column   containing $m_3$ is not the last one. We can apply the transformation (ii)$'$    adding  this  column   to the previous one.
Then   applying   (iii)$'$ to
the first row we again obtain $M(v) $. If $a\in \alpha- 
\beta$, then 
$$M(w)=     
\Bigl (  \left [ \begin{array}{ccccc}  0& 0& a&-1& 0 \\
        0&   a_\bullet & -a a_\bullet& 0& 0 \\
m_1& m_2& 0& m_3&m_4
\end{array} \right ], r^w , r^w\Bigr ).$$
It is clear that  $M(v)=     
 ( [  
m_1\, \vert \,  m_2+ m_3\, \vert \,  m_4
  ], r^v, r^v)$ where $r^v=  - aa_\bullet r^w $. Multiplying the
middle column of $M(w)$  by $a^{-1}=\tau(a)\in \Lambda$ and adding it to the next  column  we obtain a matrix in which the first row contains only one non-zero entry   
$a$. Applying  (iii)$'$ to the first row we obtain the weighted matrix
$$ \Bigl (\left  [ \begin{array}{ccccc}   
        0& a_\bullet  & -  a_\bullet  & 0 \\
m_1& m_2 & m_3&m_4
\end{array} \right ], (-1)^{k  } a r^w, (-1)^{k+1 }a  r^w \Bigr )$$
where $k$ is the   length of   $x$.  If $x\neq \emptyset$,  then we can    add the   column  containing $m_2$  to the next one and then apply (iii)$'$ to  the first row to  obtain
$  {M(v)}  $. If $x=\emptyset$, then $y\neq \emptyset$ and we proceed as above.

Consider the       move $w=(\A,  xA 
B  y BAz) \mapsto (\A-\{A, B\},   xyz)  =v  $ where  $\vert A\vert =\tau (\vert B\vert) $.  Applying if necessary     several          inverse first homotopy moves
to
$w$, we can assume that the words $x,y$, and $z$ are non-empty.
   We numerate the letters of
$\A$ starting with
$A, B$. Set
$a=\vert A\vert$, $b=\vert B\vert =\tau(a)$, and $r=r^w\in \Lambda^{ab}$. By the definition of $\Lambda$, we have $ab= a_\bullet b_\bullet=1\in \Lambda$. Note that $r^v=r^w=r$. This is obvious if $a,b\in \beta$ and follows from the
equality $a a_\bullet b b_\bullet =a b a_\bullet   b_\bullet =1$ in $\Lambda^{ab}$ if $a,b\in \alpha - \beta$. Suppose that
$a,b\in
\beta$.   Then 
$$M(w)=     
 \Bigl ( \left [ \begin{array}{ccccccccc}  0& a&-1& 0& 0 & 0& 0& 0& 0\\
        0& 1- a a_\bullet & 0& 0& 0& 0 & a_\bullet& -1& 0 \\
  0& 0  & b &  -1 & 0& 0& 0& 0& 0\\
 0& 0  & 1-bb_\bullet  &  0 & 0& b_\bullet & -1& 0& 0\\
m_1& m_2& 0& m_3&m_4 & m_5&0 & m_6&m_7
\end{array} \right ], r,r \Bigr )$$
where $m_1, m_4, m_7$ are non-void matrices and $m_2, m_3, m_5, m_6$ are columns over $\Lambda$, all with the same number of rows equal to twice the length of
$v$. It is clear that $$ M(v)=     
 ( [  
m_1\, \vert \,  m_2+ m_3\, \vert \,  m_4 \, \vert \,   m_5+ m_6 \, \vert \,  m_7
  ], r,r).$$ Multiplying by
$b\in
\Lambda$ the   column of
$M(w)$ containing
$m_3$ and adding it to the previous column    we obtain a matrix in which the third row contains only one non-zero entry   
$-1$. Applying   (iii)$'$ to this row we obtain the weighted matrix
$$ \Bigl (\left [ \begin{array}{cccccccc}  0& a&-1 & 0 & 0& 0& 0& 0\\
        0& 1- a a_\bullet &   0& 0& 0 & a_\bullet& -1& 0 \\
 0& 0  & 1-bb_\bullet    & 0& b_\bullet & -1& 0& 0\\
m_1& m_2 &   m_3 b&m_4 & m_5&0 & m_6&m_7
\end{array} \right ] ,  (-1)^{k} r, (-1)^{k +1} r \Bigr )$$
where $k$ is the number of columns in $m_1$, that is the length of  the word $x$.  Similarly using the entry $-1$ in the second row and applying  
(ii)$'$, (iii)$'$,  we can transform  the latter weighted matrix into 
$$ \left [ \begin{array}{cccccccc}  0& a&-1 & 0 & 0& 0 & 0\\
 0& 0  & 1-bb_\bullet    & 0& b_\bullet &   -1& 0\\
m_1& m_2+ m_6 (1- a a_\bullet) &   m_3 b &m_4 & m_5&  m_6 a_\bullet&m_7
\end{array} \right ]$$ with both weights equal to $(-1)^{l+1} r $
where $l$ is the number of columns in $m_4$. 
The latter matrix is similarly transformed into
$$  \left [ \begin{array}{cccccccc}  
 0&  ( 1-bb_\bullet )a   & 0& b_\bullet &   -1& 0\\
m_1& m_2+ m_6 (1- a a_\bullet) +  m_3 b a &m_4 & m_5& m_6 a_\bullet &m_7
\end{array} \right ] $$ with weights $   (-1)^{k+l  } r, (-1)^{k+l +1} r   $
and finally into 
$$   ([    
m_1\, \vert \, m_2+ m_6 (1- a a_\bullet) +  m_3 ba +m_6 a_\bullet ( 1-bb_\bullet )    a \, \vert \,m_4 \, \vert \, m_5+   m_6  a_\bullet b_\bullet \, \vert \,m_7
  ] , r, r   ).$$
Since $ab= b_\bullet a_\bullet=1$,  the latter weighted matrix coincides with $M(v)$.    

The case where $a,b\in \alpha - \beta$ is treated similarly. For completeness, we specify the matrix $M(w)$ in this case  leaving to the reader to verify that it is
equivalent to
$M(v)$:
 $$M(w)=     
 \Bigl ( \left [ \begin{array}{ccccccccc}  0&  0& 0 & 0& 0& 0&a&-1&  0\\
        0& a_\bullet& -1&    0& 0& 0&1- a a_\bullet & 0 &  0 \\
  0& 0  &  0& 0& 0&b &  -1 &  0& 0\\
 0& 0  &b_\bullet & -1&       0 & 1-bb_\bullet  & 0&  0& 0\\
m_1& m_2& 0& m_3&m_4 & m_5&0 & m_6&m_7
\end{array} \right ], r,r \Bigr ).$$

Consider the    move $w=(\A, xAByACzBCt)\mapsto  (\A, xBAyCAzCBt)=v$ where $\vert A\vert=\vert B\vert=\vert C\vert$.  
Applying if necessary     several          inverse first homotopy moves
to
$w$, we can assume that the words $x,y, z,t$ are non-empty. Set $a=\vert A\vert\in \alpha $,  $q=1- a a_\bullet \in \Lambda$,  and $r=r^w\in \Lambda^{ab}$.  Clearly,  $r^v=r^w=r $.  If
$a\in
\beta$, then 
$$M(w)\!=\!     
 \left [ \begin{array}{ccccccccccccc}  0& a&-1& 0& 0 & 0& 0& 0& 0& 0& 0 & 0& 0\\
        0& q & 0& 0& 0&   a_\bullet& -1& 0 & 0& 0 & 0& 0 & 0\\
  0& 0  & a &  -1 & 0& 0& 0& 0& 0& 0& 0& 0& 0\\
 0& 0  & q  &  0 & 0&  0& 0& 0& 0&  a_\bullet & -1& 0& 0\\
  0& 0& 0& 0& 0 & 0  & a&-1& 0& 0 & 0 & 0 & 0\\
0& 0& 0& 0& 0 & 0  & q & 0& 0& 0 &  a_\bullet & -1& 0\\
m_1& m_2& 0& m_3&m_4 & m_5&0 & m_6&m_7 &m_8 &0 &m_9 &m_{10}
\end{array} \right ] $$
with weights $r ,r $ and
$$M(v)\!=\!    \left [ \begin{array}{ccccccccccccc}  0& a&-1& 0& 0 & 0& 0& 0& 0& 0& 0 & 0& 0\\
        0& q & 0& 0& 0&    0 & 0& 0 & 0& 0 & a_\bullet& -1& 0\\
  0& 0  & a &  -1 & 0& 0& 0& 0& 0& 0& 0& 0& 0\\
 0& 0  & q  &  0 & 0&  0& a_\bullet & -1& 0& 0& 0&   0& 0\\
  0& 0& 0& 0& 0   & a&-1& 0& 0 & 0 & 0 & 0 & 0\\
0& 0& 0& 0& 0 & q & 0& 0  &  0 &  a_\bullet & -1& 0 & 0\\
m_1& m_2& 0& m_3&m_4 & m_5&0 & m_6&m_7 &m_8 &0 &m_9 &m_{10}
\end{array} \right ] $$
with weights $r ,r $.
It is an exercise in linear algebra to verify that $ M(w) $ is equivalent to the   matrix 
$$   [    
m_1\, \vert \, m_2+  m_3 a^2 +  m_6 aq +  m_9 q\, \vert \,  m_4 \, \vert \,  m_5 + m_6 aa_\bullet  +  m_9  a_\bullet q \, \vert \, m_7
\, \vert \, m_8+ m_9 a_\bullet^2 \, \vert \, m_{10}  ]  $$
with weights $r,r$. The matrix   $ M(v) $ is equivalent to the same   matrix with weights $-r,-r$.
Thus, $M(v) $ is equivalent to $\div {M(w)}$.  Similarly, if  $ a  \in \alpha-\beta$, then 
$ M(w) $ is equivalent to the   matrix 
$$   [    
m_1\, \vert \, m_2+  m_3 a_\bullet^2   \, \vert \,  m_4 \, \vert \,  m_5  +  m_3  a_\bullet q
+  m_6 aa_\bullet \, \vert \, m_7
\, \vert \, m_8+ m_3 q + m_6   aq + m_9 a^2 \, \vert \, m_{10}  ]  $$
with weights $-r,-r$  and  $ M(v) $ is equivalent to the same   matrix with weights $r,r$.
\end{proof}

To describe the behavior of the weighted matrix   under operations on nanowords,  we need   two involutions on the set of weighted
matrices over $\Lambda$.  Define first a ring involution  $\lambda\mapsto \overline \lambda: \Lambda\to \Lambda$   by 
$ \overline a=a^{-1} = \tau(a)$ and $ \overline{a_\bullet}=a^{-1}_\bullet =\tau(a)_\bullet$  for   $a\in \alpha$.  Given a weighted matrix $M$, denote by $\overline M$ the matrix obtained by applying this ring involution to all  entries of $M$  and applying the induced involution of $\Lambda^{ab}$ to the weights. For an $(n\times m)$-matrix $ M=(M_{i,j})_{i,j}$,
 set $M^-=(M_{n+1-i, m+1-j})_{i,j}$.  The second involution on the set of weighted
matrices over $\Lambda$  transforms   $(M ,
r_-, r_+)$   into   $ (M^-,   r_+, 
  r_-)$ if the deficiency of $M$ is odd and into $ (M^-,   r_-, 
  r_+)$ if the deficiency of $M$ is even. Both   involutions are compatible with the equivalence, that is transform equivalent weighted matrices into equivalent weighted matrices. 

It follows from the definitions that
$M(\overline w)=
\overline {M(w)}$. The next lemma computes  $M_{\alpha-\beta}(w^-)$ from $M_\beta (w)$.

\begin{lemma}\label{modulesnezzd}  For any nanoword $(\A,w )$, the weighted matrix 
$M_{\alpha-\beta}(w^-)$ is obtained from
${\overline {M_\beta (w)}\,}^-$ by transformations (i)$'$ -- (iii)$'$ and  the  following        transformations on  weighted   matrices: (iv)$'$   addition  to a row of another row multiplied on the left by an
  element of
$\Lambda$, the weights being preserved  and (v)$'$ multiplication of a row by an invertible element $\lambda\in \Lambda$, the weights being multiplied by
$\lambda^{-1}$. 
\end{lemma}
    \begin{proof}    Let us numerate the letters of $\A$. 
Recall that the $k$-th letter $A_k$ contributes  2 rows to $M_\beta (w)$.    Set  $a=\vert A_k\vert , i=i_{A_k},
j=j_{A_k}$. Assume first that $i \leq  j -2$. If $a\in \beta$, then the  rows  contributed by $A_k$ are
$$\left  [ \begin{array}{ccccccc}    
        { {0}}_{i-1} & a  & -  1  & { {0}}_{j-i-2}  &    0  & 0 & { {0}}_{n-j} \\
{ {0}}_{i-1} & 1  -aa_\bullet & 0 & { {0}}_{j-i-2}  &  a_\bullet & -1 & { {0}}_{n-j}
\end{array} \right ]$$
where   ${ {0}}_s$ stands for a row of $s$ zeros. The corresponding rows in  ${\overline {M_\beta (w)}\,}^-$
are
$$ \left  [ \begin{array}{ccccccc}    
 { {0}}_{n-j} &  -1  & a_\bullet^{-1}  & { {0}}_{j-i-2}  &   0 & 1  -a^{-1} a_\bullet^{-1} & { {0}}_{i-1}\\
{ {0}}_{n-j} & 0  & 0  & { {0}}_{j-i-2}  &   -1 & a^{-1} & { {0}}_{i-1}
\end{array} \right ].$$
Permuting these two rows we obtain
$$ \left  [ \begin{array}{ccccccc}    
        { {0}}_{n-j} & 0  & 0  & { {0}}_{j-i-2}  &   -1 & a^{-1} & { {0}}_{i-1} \\
{ {0}}_{n-j} &  -1  & a_\bullet^{-1}  & { {0}}_{j-i-2}  &   0 & 1  -a^{-1} a_\bullet^{-1} & { {0}}_{i-1}
\end{array} \right ].$$
Multiply the first row by $a_\bullet^{-1}-a$ and add it to the  second
row. This gives $$ \left  [ \begin{array}{ccccccc}    
        { {0}}_{n-j} & 0  & 0  & { {0}}_{j-i-2}  &   -1 & a^{-1} & { {0}}_{i-1} \\
{ {0}}_{n-j} &  -1  & a_\bullet^{-1}  & { {0}}_{j-i-2}  &   a-a_\bullet^{-1}  & 0   & { {0}}_{i-1}
\end{array} \right ].$$
Multiplying the first row by  $-a$ and the second row by $-a_\bullet$ we obtain 
$$ \left  [ \begin{array}{ccccccc}    
        { {0}}_{n-j} & 0  & 0  & { {0}}_{j-i-2}  &   a & -1 & { {0}}_{i-1} \\
{ {0}}_{n-j} &  a_\bullet  &  {-1}  & { {0}}_{j-i-2}  &   1 -a a_\bullet   & 0   & { {0}}_{i-1}
\end{array} \right ].$$
These are exactly the rows contributed by   $A_k$ to $M_{\alpha-\beta} (w^-)$.  
 The cases where $i=j-1$ or  $a\in \alpha - \beta$ are  similar. That the weights behave appropriately is left to the reader as an exercise.
\end{proof}

Given    two nanowords $w, w'$,  we can  compute  the weighted matrix  $M(ww')$ from $M(w)= (M, r_-, r_+)$ and  $M(w')= (M', r'_-, r'_+)$.  Let  us write
$M =[  M_+\,\, m_+ ]$ where
$m_+$ is the last column of $M$ and $M_+$ is obtained from $M$ by deleting $m_+$. Similarly,   $M' =[  m'_-\,\,    M'_-]$ where $m'_-$ is
the first column of $M' $. Then 
$$M(ww')=\Bigl (\left  [ \begin{array}{ccc}    
       {M_+}  & m_+  & 0   \\
  0 & m'_- & {M'_- }
\end{array} \right ], r_-r'_-, r_+r'_+\Bigr ).$$

\subsection{Module $K_\beta (w)$}\label{3matrinvaricovfunct} By   Sect.\ \ref{matrinvaricovfunct}, the  weighted matrix $M_\beta(w)$ of a nanoword
$w$ gives rise to a  marked $\Lambda$-module   
$ K_\beta (w) = K(M_\beta(w))$. Lemma \ref{nezzd}
implies that  
$ K_\beta(w) $  is preserved under homotopy moves on
$w$. It is clear that   $K_\beta(\overline w)=\overline { K_\beta(w)}$ where for a  $\Lambda$-module $K$, the $\Lambda$-module $\overline K$   is obtained
from
$K $ by pulling back along the  involution 
$\lambda\mapsto
\overline
\lambda: \Lambda\to \Lambda$.     Lemma \ref{modulesnezzd} implies that  $K_{\alpha-\beta}(w^-)$ is obtained from
$\overline {K_\beta (w)}$ by permuting the input with the output.

For a  $\Lambda$-module $X$, the
$X$-colorings of $w$ bijectively correspond to $\Lambda$-linear homomorphisms   $K_\beta(w)\to X$. Hence  for
finite (as a set) 
 $X$, the number of $X$-colorings of $w$ with given input and output is a homotopy invariant of $w$. 

Further homotopy invariants of $w$ are provided   by    the elementary
ideals and other invariants  of the   $\Lambda^{ab}$-module $\Lambda^{ab}\otimes_\Lambda K_\beta(w)$ and  of its quotients by the input  and/or the output.
 We   discuss  two  such invariants in the next subsection. 
 
\subsection{Polynomials $\nabla^\pm_\beta (w)$}\label{epqmovfunct}  Since the deficiency of the 
weighted matrix $M=M_\beta (w)$   is equal to 1,   Sect.\ \ref{matrinvaricovfunct}  yields  two polynomials  $\nabla^- (M  ), \nabla^+ (M)\in \Lambda^{ab}$. Lemma \ref{nezzd} implies
that they are homotopy invariants of $w$ at least up to multiplication by $-1$.  We  eliminate the sign indeterminacy as follows.  Consider the homomorphism
${\rm {aug}}:\Lambda^{ab}\to \ZZ$ sending all $a, a_\bullet$ to $1$.  In other words,  the value of ${\rm {aug}}$ on $\lambda\in \Lambda^{ab}$ is the sum of coefficients of $\lambda$.  Applying ${\rm {aug}}$  to all entries of $M$ and permuting the rows, we obtain the matrix 
$$  \left  [ \begin{array}{ccccccc}    
        1 & -1  & 0  & 0   & \cdots   \\
0 &  1  & -1  & 0  & \cdots   \\
& & \cdots\\
 \cdots   & 0  &   0 & 1 & -1   \\
\end{array} \right ].$$
 Therefore ${\rm {aug}}(\nabla^\varepsilon (M  ))=\pm 1$ for   $\varepsilon=\pm 1$. The polynomial $$\nabla^\varepsilon_\beta(w)= {\rm {aug}}(\nabla^\varepsilon (M  )) \, \nabla^\varepsilon (M  ) \in \Lambda^{ab} $$ is a homotopy invariant  of $w$  without any indeterminacy.  Clearly, ${\rm {aug}} (\nabla^\varepsilon_\beta (w ))= 1$. 
  The results of Sect.\ \ref{2matrinvaricovfunct} imply  that $$\nabla^\varepsilon_\beta(\overline w)= \overline{\nabla^\varepsilon_\beta(w)},\,\,\,\, \nabla^\varepsilon_{\alpha-\beta} (w^-) =\overline {  \nabla^{-\varepsilon}_\beta(w) },\,\,\,\, \nabla^\varepsilon_\beta(ww')=\nabla^\varepsilon_\beta(w) \nabla^\varepsilon_\beta(w').$$

\section{Invariant $\lambda$}\label{dd822}

 We  use  the modules of nanowords to define a $\Lambda$-valued homotopy invariant of nanowords $\lambda$. As an application we give a homotopy classification of nanowords of length 6. 

\subsection{Invariants $\lambda'$ and  $\lambda$}\label{spebetnct} For a nanoword $(\A,w:\hat n\to \A)$,  consider the $\Lambda$-module $ K_\alpha (w)$ corresponding to $\beta=\alpha$. This module is free of rank 1. Indeed, let $v_-=x_0,
x_1,..., x_{n}=v_+
 $ be the generators of    $ K_\alpha (w)$ determined by  the first, second, etc.\ columns of the matrix $M=M_\alpha(w)$. The two rows of
$M$  derived from a letter
$A \in
\A$ express  $x_{i_A}$ as a linear combination of $x_{s}$ with $s< i_A$ and  express  $x_{j_A}$ as a linear combination of $x_{s}$ with $s< j_A$.
Therefore  $K_\alpha (w)=\Lambda v_-$. Then $v_+= \lambda'(w) v_-$ for a unique
$ \lambda' (w) \in \Lambda$ which  is a homotopy invariant of $w$.  It admits an equivalent but  more convenient version $\lambda(w)$ defined as follows.  Consider the involutive anti-automorphism $\iota: \Lambda\to \Lambda$   equal to  the identity on the   generators $\{a, a_\bullet\}_{a\in \alpha}\subset \Lambda$. Thus  $\iota$ acts on   monomials   by reading them  from right to left.  For instance, $\iota( a b a_\bullet)= a_\bullet b a $. Set $$\lambda(w)= \iota( \lambda'(w))\in \Lambda.$$

We describe  a simple  method allowing to compute $\lambda(w)$ (and extending the Silver-Williams method \cite{sw} in the setting of virtual  strings).     For each $i=1,...,n$ there is a unique row in $M$ expressing the generator $x_i\in K_\alpha (w)$ via the generators with smaller indices. This relation has either the form $x_i= k_i x_{i-1}$ with $k_i\in \Lambda$ or the form
 $x_i= k_i x_{i-1} +l_i x_{i'}$ with $k_i, l_i\in \Lambda$, $l_i\neq 0$,  and $i'<i-1$. We form a graph  $\Gamma_w$ with $n+1$ vertices $\{0,1,2,...,n\}$  and $3n/2$ edges as follows. A relation of the form $x_i= k_i x_{i-1}$ gives rise to one edge - the segment $[i-1,i] \subset \mathbf R  $ with endpoints $i-1, i$. We place on this edge the weight $k_i$. A relation of the form $x_i= k_i x_{i-1} +l_i x_{i'}$ gives rise to two edges:  the segment $[i-1,i] \subset \mathbf R$   and an  arc connecting $i$ with $i'$ (and disjoint from the rest of $\Gamma_w$). We place  on these two edges the weights $k_i$ and $l_i$, respectively. Consider  a    path   leading from $n$ to $0$ and following along the edges of $\Gamma_w$, always from the bigger endpoint to the smaller one.   Compute the product of weights   associated with the edges lying on the path writing the    weights from right  to left in the order of  appearance of the edges on the path. Summing these products over all such paths from  $n$ to $0$ we obtain   $\lambda(w)$.
 
 \begin{figure}
\centerline{\includegraphics[width=5.8cm]{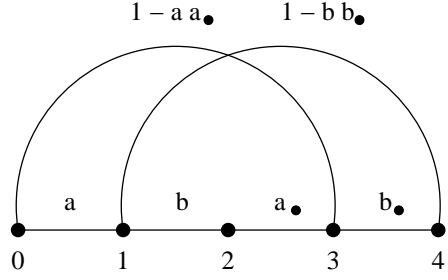}}
\caption{Graph  $\Gamma_{w}$ for $w=ABAB$}\label{figure3}
\end{figure}

 Example:  for    $w=ABAB$ with $\vert A\vert =a\in \alpha,  \vert B\vert =b\in \alpha $, there are three   paths in $\Gamma_w$ contributing to $\lambda(w)$, cf.\  Figure \ref{figure3}. This gives  
  \begin{equation}\label{xyz} \lambda  (w)=ab a_\bullet b_\bullet+ (1-aa_\bullet) b_\bullet+ a(1- bb_\bullet).\end{equation}
 We   describe  a more general class of nanowords for which $\lambda$ can be  computed explicitly. Pick a finite sequence of letters (possibly with repetitions)  $a_1,..., a_m\in  \alpha$  and a permutation $\sigma: \hat m \to \hat m$.
   Consider the nanoword  $(\A=\{A_1,..., A_m\},w:\widehat{2m} \to \A)$ where $\vert A_i \vert =a_i$,  $w(i)=A_i$  and $w(i+m)= A_{\sigma (i)}$ for $i= 1,..., m$. Thus $w=A_1A_2\cdots A_m A_{\sigma (1)} A_{\sigma (2)} \cdots A_{\sigma (m)}$ with  $\vert A_i \vert =a_i$.  The corresponding graph $\Gamma$ contains precisely $m+1$ paths from $2m$ to $0$ and \begin{equation}\label{fir}   \lambda(w)=
 a_1 a_2 \cdots a_m (a_{\sigma(1 )})_\bullet \,  (a_{\sigma(2 )})_\bullet \cdots  (a_{\sigma(m )})_\bullet \end{equation}
   $$+\sum_{i=1}^m  \, a_{ 1} a_{ 2} \cdots a_{\sigma(i)-1}  (1-  a_{\sigma(i )} (a_{\sigma(i )})_\bullet)  (a_{\sigma(i+1)})_\bullet \, (a_{\sigma(i+2)})_\bullet \cdots (a_{\sigma(m )})_\bullet  
     .$$

\subsection{Properties of  $\lambda$}\label{dmmmetnct}    It follows from the definitions  that $\lambda(\overline w)= \overline {\lambda(w)}$ and $\lambda(ww')=\lambda(w) \lambda(w')$ for any nanowords $w, w'$.  
The graph  $\Gamma_{w^-}$ is obtained  from $\Gamma_w$ by relabeling the vertices ($i$ becomes $n-i$) and exchanging   $a\leftrightarrow a_\bullet$ for all $a\in \alpha$. Therefore 
$ \lambda(w^-)=\kappa  (\lambda (w))$ where $\kappa$ is the ring anti-automorphism of $\Lambda$ defined on the generators     $a\in \alpha$  by $\kappa (a)=a_\bullet, \kappa (a_\bullet)=a$. Thus if $w$ is homotopically symmetric then $\kappa   (\lambda (w))=\lambda(w)$. If  $w$ is homotopically skew-symmetric then $\overline {\kappa  (\lambda (w))}=\lambda(w)$.  If $w$ is contractible, then $\lambda (w)=1$. 

To state further properties of $\lambda(w)$, recall the group   $\Pi$  defined in Sect.\ \ref{groupga} and consider the   ring homomorphisms $ p, r, r_\bullet:\Lambda\to \ZZ \Pi$    defined on the generators of  $\Lambda$ by   $ p (a)=z_a, \, p(a_\bullet)= z_{\tau(a) }$, $ r(a)=z_a,\,  r(a_\bullet)=1$ and $ r_\bullet (a)=1,\,  r_\bullet (a_\bullet)=z_a$.   
 
\begin{lemma}\label{42garalaezzd}  For any nanoword $w$, 
$p(\lambda(w))=  \gamma(w) \in \Pi  $   is the invariant  defined in Sect.\ \ref{groupga} and  $r(\lambda(w))=r_\bullet (\lambda(w))=  1\in  \ZZ \Pi$.\end{lemma}
    \begin{proof}     
    Substituting $a_\bullet = a^{-1}$ in   the matrix $M_\alpha (w)$    we obtain    $x_{i_A+1}= \vert A\vert x_{i_A}$ and $x_{j_A+1}= \vert A\vert^{-1} x_{j_A}$. Comparing   with the definition of $\gamma(w)$,
    we obtain  $p(\lambda(w))=  \gamma(w)$. Substituting $a_\bullet = 1$ in   the matrix $M_\alpha (w)$   and adding all rows we obtain a row with first entry 1, last entry $-1$, and all other entries $0$. Hence $r(\lambda(w))=  1$. Substituting $ \lambda(w^-)=\kappa  (\lambda (w))$  in $r(\lambda(w^-))=  1$ we obtain a  formula equivalent to $r_\bullet (\lambda(w))=  1$. 
\end{proof}

The polynomials $\nabla^-_\alpha (w), \nabla^+_\alpha(w)$ defined in Sect.\ \ref{epqmovfunct} (for $\beta=\alpha$) can be computed   as follows:   $\nabla^-_\alpha(w)=1$ and $\nabla^+_\alpha (w)=  q(\lambda'(w))=  q(\lambda(w))$  where  $q$ is the projection $\Lambda \to \Lambda^{ab}$. Thus $\lambda(w)$ is   a  common  extension of    $\nabla^+_\alpha (w)$ and $\gamma(w)$.

     \subsection{Derived invariants}\label{dfmmmexcebetnct} 
     For a  monomial $m\in \Lambda$ in the generators $\{a,a_\bullet\}_{a\in \alpha}$, denote    $\deg(m)$ the  number of generators without $\bullet$   appearing in $m$ and denote  $\deg_\bullet (m)$
  the   number of generators with $\bullet$ in $m$.  For instance $\deg(aba_\bullet )=2$ and $\deg_\bullet (aba_\bullet )=1$. We have 
$$\Lambda=\Lambda_{0,0}\oplus \Lambda_{0,1}\oplus  \Lambda_{1,0}\oplus  \Lambda_{1,1}$$
where $\Lambda_{i,j} \subset \Lambda$ is additively generated by the monomials $m$   such that   $\deg(m)=i\, ({\rm {mod}} \,2)$  
and   $\deg_\bullet (m)=j\, ({\rm {mod}}\, 2)$. As a corollary,  the invariant $\lambda(w)$ of a nanoword $w$ splits as a sum of 4 homotopy invariants 
$$\lambda (w) =\lambda_{0,0} (w) + \lambda_{0,1}(w) + \lambda_{1,0}(w) +  \lambda_{1,1}(w) $$
such that $\lambda_{i,j} (w)\in \Lambda_{i,j} $ for all $i,j$.

Consider the ring homomorphism  $\psi: \Lambda \to \ZZ \Pi \otimes_{\ZZ}  \ZZ\Pi$
  defined on the generators $a\in \alpha$  of  $\Lambda$ by $  \psi(a)=z_a\otimes 1$ and $ \psi( a_\bullet) =1\otimes z_a$. For instance,  $\psi (ab a_\bullet)=z_a   z_b \otimes z_a$.  The set $\{x \otimes y\}_{x,y\in \Pi}$ is a basis of the   additive group   $\ZZ \Pi \otimes_{\ZZ}  \ZZ\Pi$. For      $\lambda\in \Lambda$, we have a unique (finite)   expansion
 $\psi (\lambda )=\sum_{x,y\in \Pi}  \lambda^{x,y}   \, x \otimes y$ with $\lambda^{x,y}\in \ZZ$. For a  nanoword $w$,  this gives   a set of   homotopy invariants  $\{(\lambda  (w))^{x,y} \in \ZZ \}_{x,y\in \Pi}$.

 \subsection{Proof of Lemma \ref{eril}}\label{erilpotnct}   
\begin{figure}
\centerline{\includegraphics[width=9cm]{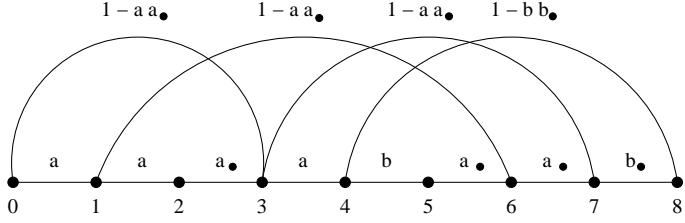}}
\caption{Graph  $\Gamma_{w}$ for $w=A_3A_2A_3A_1BA_2A_1B$}\label{figure4}
\end{figure}

For $w=A_3A_2A_3A_1BA_2A_1B$, there are 7 paths in the graph $\Gamma_w$ contributing to $\lambda(w)$, see Figure \ref{figure4}.  We have
  $$\lambda (A_3A_2A_3A_1BA_2A_1B)=
  (1-aa_\bullet)^2 b_\bullet + (1- a a_\bullet) a (1-b b_\bullet)+ (1-a a_\bullet) ab {a^2_\bullet}  b_\bullet$$
  $$+  a (1-a a_\bullet) a_\bullet b_\bullet + a^2 a_\bullet (1-aa_\bullet) b_\bullet +a^3 a_\bullet (1-b b_\bullet) +a^3 a_\bullet b a_\bullet^2 b_\bullet.$$
  This implies that $\lambda_{1,1} (A_3A_2A_3A_1BA_2A_1B)=a^3a_\bullet -a^3 a^2_\bullet b_\bullet$. In particular if $a=\tau(a)  $, then we have $\lambda_{1,1} (A_3A_2A_3A_1BA_2A_1B)=a a_\bullet -a  b_\bullet$. Similarly,
  $$\lambda (A_3A_2BA_3A_1 A_2A_1B)=
  (1-aa_\bullet)^2 b_\bullet + (1- a a_\bullet) a a_\bullet^2 b_\bullet+ a (1-a a_\bullet)  {a_\bullet} b_\bullet$$
  $$+  a^2 (1-b b_\bullet)   + a^2 b a_\bullet (1-aa_\bullet) b_\bullet +a^2b a a_\bullet^3   b_\bullet .$$
  This implies that $\lambda_{1,1} (A_3A_2BA_3A_1A_2A_1B)=a a^2_\bullet b_\bullet-b b_\bullet$. 
  For  $a=\tau(a)  $,   we have $\lambda_{1,1} (
  A_3A_2BA_3A_1A_2A_1B)=a   b_\bullet-b b_\bullet$. The equality 
  $a a_\bullet -a  b_\bullet =a   b_\bullet-b b_\bullet$ implies that    $a a_\bullet =a   b_\bullet $ so that $a=b$ which contradicts the assumptions of the lemma. 
  
  \subsection{Proof of Lemma \ref{6eril}}\label{e666lpotnct}    A desingularization of  $w^6_{a,b}$ gives the nanoword $A_3A_2BA_3A_1BA_2A_1 $ in the $\alpha$-alphabet (\ref{bet}). We have
  $$\lambda (A_3A_2BA_3A_1BA_2A_1 )=
  (1-aa_\bullet)^2  +  (1-aa_\bullet) a b_\bullet a^2_\bullet  + a (1-a a_\bullet) a_\bullet$$
  $$+  a^2 (1-b b_\bullet) a_\bullet^2  + a^2 b a_\bullet(1- a a_\bullet)   +a^2 b a a_\bullet  b_\bullet  a_\bullet^2.$$
  Substituting here $a^2=a_\bullet^2=b^2=b_\bullet^2=1$ (due to the assumptions $a=\tau (a), b=\tau (b)$, we easily obtain  that 
  $$\lambda_{1,1}=\lambda_{1,1} (A_3A_2BA_3A_1BA_2A_1)=
  2-ba-a_\bullet b_\bullet+ba a_\bullet b_\bullet.$$
  Note that $ba\neq 1$ and $a_\bullet b_\bullet\neq 1$ since $\tau (a)= a\neq b$.
  Therefore $\lambda_{1,1}  \neq 0$. Moreover, knowing 
  $\lambda_{1,1}  $ we can recover $ba$: This is the only non-trivial element   $x\in \Pi$ such that  $ (\lambda_{1,1} )^{x,1}=-1$. 
  The equality $ba=b'a'$ implies that $b=b', a=a'$ which concludes the proof of the lemma.
  
    \subsection{Remark}\label{k50000986jaga} For a nanoword $w=A_1\cdots A_n$, set     $w_*=\vert A_1\vert_* \cdots \vert A_n\vert_*\in \Lambda$ where $\vert A_i \vert_*= \vert A_i \vert$ if $i$ numerates the first occurence of the letter $A_i$ in $w$ and  $\vert A_i \vert_*= \vert A_i \vert_\bullet$ if $i$ numerates the second occurence of $A_i$ in $w$.
 The computation of $\lambda$ in Sect.\ \ref{spebetnct}  implies that $\lambda(w)-w_*$ is an algebraic  sum of monomials $m$ in the generators $a, a_\bullet$ such that
 $\max (\deg(m) , \deg_\bullet (m)) \leq n/2$ and  $\deg(m)+\deg_\bullet (m) < n$.
 
  \subsection{Exercise}\label{exerc55ga} Deduce from Lemma \ref{42garalaezzd} that if $\alpha$ is a one-element set, then $\lambda(w)=1$ for any nanoword $w$ over $\alpha$.

  \section{Analysis litterae  II}\label{dqmqmqmqmqm2}

  In this section we give a homotopy classification of nanowords of length $\leq 6$. We begin with 
   nanowords of length 4. We shall use the homomorphism $\psi$ introduced in Sect.\ \ref{dfmmmexcebetnct}.
  
 \subsection{Nanowords of length 4 re-examined}\label{excrifgbetnct}    Let $w=ABAB$ with $\vert A\vert =a\in \alpha,  \vert B\vert =b\in \alpha $. Formula \ref{xyz} implies that 
$$\lambda_{0,0}  (w)= ab a_\bullet b_\bullet,\,  \lambda_{1,0}  (w)= 
a -aa_\bullet b_\bullet,\,  \lambda_{0, 1}  (w)= b_\bullet -abb_\bullet,\,  \lambda_{1,1}  (w)= 0.$$  
  Then  $\psi(\lambda_{0, 1}  (w))= (1-z_a z_b)\otimes z_b$. Note that  $z_a z_b=1\in \Pi$ if and only if $a=\tau(b)$. Therefore  $\psi(\lambda_{0, 1}  (w))=0$ implies that $ a=\tau(b)$ and then
  $w$ is contractible. Conversely, if $w$ is contractible then $\lambda_{0, 1}  (w)=0$. We conclude that $w$ is contractible  if and only if $\psi(\lambda_{0, 1}  (w))=0$. If $\psi(\lambda_{0, 1}  (w)) \neq 0$, then we can recover $a$ and $b$ from $\psi(\lambda_{0, 1}  (w))$.  This gives an alternative proof of  Theorem \ref{classid} and shows that $\lambda$ is a faithful homotopy invariant of nanowords of length 4.

\subsection{Nanowords of length 6}\label{fgsszeiikaga}  Pick three letters $a,b,c\in \alpha$ (possibly coinciding). 
Let $\A$ be the  $\alpha$-alphabet  consisting of 3 letters $A,B,C$ with $\vert A\vert =a,  \vert B\vert =b, \vert C \vert = c$. Consider the nanowords  $ w^1_{a,b,c}=ABCABC$ and $$ w^2_{a,b,c}=ABCACB,  
w^3_{a,b,c}=ABCBAC, w^4_{a,b,c}=ABCBCA, w^5_{a,b,c}=ABACBC$$ in this $\alpha$-alphabet.   It is easy to check that any nanoword of length 6 is either  homotopic to a nanoword of length $\leq 4$ or is isomorphic to   $w^i_{a,b,c}$ with $i\in \{1,...,5\}$.   We now point out   obvious sufficient conditions for   $w^i_{a,b,c}$ to be contractible. 
 
If $a=\tau(b)$ or $c=\tau(b)$, then $w^1_{a,b, c}=ABCABC$ is contractible. 
 We say that an (ordered) triple $a,b,c \in \alpha$ is {\it 1-regular} if $ a \neq \tau(b)\neq c$.
 
If   $c=\tau(b)$, then $w^2_{a,b, c}=ABCACB$ is   contractible. 
A  triple $a,b,c \in \alpha$ is {\it 2-regular} if $ c \neq \tau(b)$.
  
If   $a=\tau(b)$, then $w^3_{a,b, c}=ABCBAC$ is contractible. 
A  triple $a,b,c \in \alpha$ is {\it 3-regular} if $ a \neq \tau(b)$.
  
If   $c=\tau(b)$, then $w^4_{a,b, c}=ABCBCA$ is contractible. 
A  triple $a,b,c \in \alpha$ is {\it 4-regular} if $ c \neq \tau(b)$. (This coincides with the 2-regularity).

  If   $ a=b=c=\tau(a)$, then $w^5_{a,b, c}=ABACBC$ is contractible by the moves $ABACBC \mapsto BACACB \mapsto BB \mapsto \emptyset$.
  We say that a triple $a,b,c \in \alpha$ is    {\it singular} if $  a=b=c=\tau(a)$ and {\it 5-regular} otherwise.
  
  The following theorem gives a   homotopy classification of nanowords of length $6$. Its proof shows that $\lambda$ is a faithful homotopy invariant of nanowords of length 6 corresponding to regular triples $a,b,c$. However,   one of the claims of the theorem - that $  w^5_{a,b,c} $ is not homotopic to a nanoword of length 4  -  needs    subtler techniques in the case $a=c=\tau(b)\neq b$. This is due to the fact that in this case  the invariant $\lambda$ and all the  other invariants of nanowords introduced above do not distinguish $w^5_{a,b,c} =ABACBC$ from  
 $w_{a,c}=ACAC$. That these nanowords are not homotopic    will be proven in Sect.\ \ref{ed56quaga}.
 
  \begin{theor}\label{1classid}   For $i=1,..., 5$ and any $i$-regular triple $a, b,c\in   \alpha$, the nanoword $w^i_{a,b, c}$ is neither contractible nor  homotopic to a nanoword of length 4. The nanowords $w^i$ corresponding to $i$-regular triples $a,b,c$ and $a',b',c'$ are homotopic if and only if $a=a'$,  $b=b'$, and $c=c'$.  For $i\neq j$, the nanowords $w^i$ corresponding to $i$-regular triples are not homotopic to nanowords   $w^j$  corresponding to $j$-regular triples  with one exception: $w^4_{a,b,c}\simeq w^5_{a,b,c} $ for  $a=b=c\neq \tau(a)$. \end{theor}
    \begin{proof}     Set $w=w^1_{a,b,c} $.  Formula (\ref{fir}) with $\sigma  =\id$  yields
$$\lambda(w)= abc a_\bullet b_\bullet c_\bullet +(1-aa_\bullet) b_\bullet c_\bullet + a(1-bb_\bullet ) c_\bullet  +ab (1-cc_\bullet ).$$
Therefore 
$\lambda_{0,1} (w)= \lambda_{1,0} (w) =0$ and
$$ \lambda_{0,0} (w)=ab+b_\bullet c_\bullet -ab b_\bullet c_\bullet  ,\,\,\,\,\,  
\lambda_{1, 1} (w)=a c_\bullet  -a a_\bullet  b_\bullet c_\bullet  -abc c_\bullet+ abc a_\bullet b_\bullet c_\bullet
 .$$
Clearly,   $\psi(1-\lambda_{0,0} (w))=(1-z_a z_b )\otimes (1-z_b z_c)$. If $w$ is contractible, then $\lambda_{0,0} (w)=1$ and we must have  $z_a z_b=1$ or  $z_b z_c=1$. Then  $a=\tau(b)$ or $c=\tau(b)$ which is excluded by the   1-regularity.

This argument shows, moreover,  that    we can recover   $z_a z_b$ and $z_b z_c$ from $\lambda_{0,0}(w)$. Therefore if $w^1_{a,b,c}\simeq w^1_{a',b',c'}$, then $z_a z_b= z_{a'} z_{b'}$ and $z_b z_c= z_{b'} z_{c'}$.
The equality $z_a z_b= z_{a'} z_{b'}\in \Pi$ holds   in only two cases: when $a=a', b=b'$ and when $a=\tau(b)$ and $a'=\tau(b')$. The second option is excluded by the 1-regularity. Therefore  $a=a'$,  $b=b'$, and $c=c'$.

That $w$ is  not homotopic to a nanoword of length 4   follows from the formula $\lambda_{0,1} (w)=0$ and the fact that   nanowords of length 4 with 
 $\lambda_{0,1} =0$ are contractible.

    Set $w=w^2_{a,b,c}$. Formula (\ref{fir}) with $\sigma (1)=1, \sigma (2)=3, \sigma (3)=2$ yields
$$\lambda(w)= abc a_\bullet c_\bullet b_\bullet +(1-aa_\bullet) c_\bullet b_\bullet + ab (1-cc_\bullet ) b_\bullet  +a  (1-bb_\bullet ).$$
Therefore $\lambda_{0,0} (w)=c_\bullet b_\bullet$, 
$\lambda_{0,1} (w)= 0$,  and
$$   \lambda_{1,0} (w) =a-abc c_\bullet b_\bullet,\,\,\,\,\,  
\lambda_{1, 1} (w)= abc a_\bullet c_\bullet b_\bullet
-a a_\bullet  c_\bullet  b_\bullet   .$$
If $w$ is contractible, then $\lambda_{1,0}(w)=0$. This   implies  $c=\tau(b)$, excluded by the 2-regularity.
For  a 2-regular triple $a, b,c$,    we can recover   $b,c$ from $ \lambda_{0,0}(w) $ and    $a$ from $ \lambda_{1,0}(w) $.  (Indeed, $x=z_a$ is the only element   of $  \Pi$ such that $(\lambda_{1,0}(w))^{x, 1}  \neq 0$.)  That 
  $w$ is not homotopic to a nanoword of length 4  follows from the formula $\lambda_{0,1} (w)=0$ and the fact that a nanoword of length 4 with 
 $\lambda_{0,1} =0$ must be contractible. That 
  $w$ is not homotopic to a nanoword of type $w^1$  follows from the formulas $\lambda_{1,0} (w^1)=0$ and  $\lambda_{1,0}(w)\neq 0$.

   Set $w=w^3_{a,b,c}$. Formula (\ref{fir}) with $\sigma (1)=2, \sigma (2)=1, \sigma (3)=3$ yields
$$\lambda(w)= abc b_\bullet a_\bullet c_\bullet + a(1-bb_\bullet) a_\bullet c_\bullet +    (1-aa_\bullet ) c_\bullet  +ab  (1-cc_\bullet ).$$
Therefore 
$\lambda_{0,0} (w)= ab$, $ \lambda_{1,0} (w) =0$, and
$$ \lambda_{0,1} (w)=c_\bullet -ab b_\bullet  a_\bullet c_\bullet,\,\,\,\,\,  
\lambda_{1, 1} (w)= abc b_\bullet a_\bullet c_\bullet
-a bc c_\bullet   .$$
If $w$ is contractible, then $\lambda_{0,1}(w)=0$. This   implies  $a=\tau(b)$, excluded by the 3-regularity.
For  a 3-regular triple $a, b,c$,    we can recover   $a,b$ from $ \lambda_{0,0}(w) $ and    $c$ from $ \lambda_{0,1}(w) $.  (Indeed, $y=z_c$ is the only element   of $  \Pi$ such that $(\lambda_{0,1}(w) )^{1, y} \neq 0$.)   That 
  $w$ is not homotopic to a nanoword of length 4  follows from the formula $\lambda_{1,0} (w)=0$ and the fact that a nanoword of length 4 with 
 $\lambda_{1,0} =0$ must be contractible. That 
  $w$ is not homotopic to a nanoword of type $w^1$ or $w^2$  follows from the formulas $\lambda_{0,1} (w^1)=\lambda_{0,1} (w^2)=0$ and  $\lambda_{0, 1}(w)\neq 0$.
    
    Set $w=w^4_{a,b,c}$. Formula (\ref{fir}) with $\sigma (1)=2, \sigma (2)=3, \sigma (3)=1$ yields
$$\lambda(w)= abc b_\bullet c_\bullet a_\bullet + a(1-bb_\bullet) c_\bullet a_\bullet +   ab  (1-cc_\bullet ) a_\bullet  + 1-aa_\bullet .$$
Therefore 
$\lambda_{0,0} (w)= 1$,  and
$$ \lambda_{1,0} (w)= a c_\bullet a_\bullet-abc c_\bullet a_\bullet,  \,\,\,\,\,   \lambda_{0,1} (w)=aba_\bullet -ab b_\bullet  c_\bullet a_\bullet,\,\,\,\,\,  
\lambda_{1, 1} (w)= abc b_\bullet c_\bullet a_\bullet
-a a_\bullet   .$$
If $w$ is contractible, then $\lambda_{1,0}(w)=0$. This   implies  $c=\tau(b)$, excluded by the 4-regularity.
For  a 4-regular triple $a, b,c$,    we can recover   $a,b, c$ from $ \lambda_{1,0}(w) $. Indeed, $x=z_a$ is the only element   of $  \Pi$ such that $(\lambda_{1,0} (w) )^{x, y}   =1$ for some $y\in \Pi$.
This recovers $a$. Similarly, $x'=z_a z_b z_c$ is the only element   of $  \Pi$ such that $(\lambda_{1,0}(w) )^{x', y}  =-1$ for some $y\in \Pi$.  Knowing $x, x'$, we  recover $x^{-1} x'=z_b z_c$ and the letters $b,c$.
That 
  $w$ is not homotopic to a nanoword of length 4  follows from the formula $\lambda_{0,0} (w)=1$ and the fact that a nanoword of length 4 with 
 $\lambda_{0,0} =1$ must be contractible. That 
  $w$ is not homotopic to a nanoword of type $w^1$ or $w^2$  follows from the formulas $\lambda_{0,1} (w^1)=\lambda_{0,1} (w^2)=0$ and  $\lambda_{0, 1}(w)\neq 0$ (in the 4-regular case).  That 
  $w$ is not homotopic to a nanoword of type $w^3$  follows from the formulas $\lambda_{1,0} (w^3)=0$ and  $\lambda_{1,0}(w)\neq 0$.

     Set $w=w^5_{a,b,c} =ABACBC$. There are 5 paths in the graph $\Gamma_w$ contributing to $\lambda(w)$.  This gives 
$$ \lambda(w)= (1-aa_\bullet)(1-c c_\bullet) +aba_\bullet (1-cc_\bullet) 
 + a(1-bb_\bullet)c_\bullet+(1-aa_\bullet) c b_\bullet c_\bullet +aba_\bullet c b_\bullet c_\bullet.$$ 
 Therefore  $$\lambda_{0,0} (w)= 1-abb_\bullet c_\bullet +a a_\bullet cc_\bullet, \,\,\,\,  \lambda_{1,0} (w)= c b_\bullet c_\bullet  -ab a_\bullet  c c_\bullet,    $$
$$  \lambda_{0,1} (w)=aba_\bullet -aa_\bullet  c b_\bullet c_\bullet, \,\,\,\,\lambda_{1, 1} (w)= a c_\bullet -cc_\bullet
-a a_\bullet  + aba_\bullet c b_\bullet c_\bullet .$$
Then $\psi( \lambda_{1,0} (w))= z_c \otimes z_b z_c -z_a z_b z_c \otimes z_a z_c$. 
Observe that $z_c = z_a z_b z_c \Leftrightarrow a=\tau (b)$ and  $z_b z_c = z_a z_c  \Leftrightarrow  b=a$.  Thus  if $ \lambda_{1,0} (w)=0$, then   $\tau(a)=a=b$. A similar analysis shows that if $ \lambda_{0,1} (w)=0$, then $\tau(b)=b=c$. Therefore if $w$ is contractible, then the triple $a,b,c$ is singular.  

 Now we show how to recover   $a,b,c$ from $w=w^5_{a,b, c}$ provided   $w$ is non-contractible. Suppose first that  $\psi( \lambda_{1,0} (w)) \neq 0$.   Then $x=z_c$ is the only element   of $  \Pi$ such that $(\lambda_{1,0} (w) )^{x, y}   =1$ for some $y\in \Pi$. This $y$ is unique and equal to $z_b z_c$. This   recovers $c$ and $b$. Then it is easy to recover $a$ from 
 $\psi( \lambda_{1,0} (w))$.  If  $\psi( \lambda_{1,0} (w)) = 0$, then  $\tau(a)=a=b$ and
 $\psi( \lambda_{1,1} (w))=a  \otimes (c_\bullet -   a_\bullet)$. Since  $w$ is non-contractible, $a \neq c$.
Thus we can  recover  $a$ and $c$
 from $\lambda_{1, 1} (w)$. 
 
 If $w$ is homotopic to a nanoword $w^1$, then $\lambda_{0,1} (w)=\lambda_{1,0} (w)=0$.
As shown above,  this  implies that   the triple $a,b,c$ is singular. 
 
  If $w $ is homotopic to a nanoword $ w^2$, then $\lambda_{0,1} (w)=0$ and $\lambda_{0,0} (w) $ is a monomial.
 The latter implies that   $ abb_\bullet c_\bullet =1 $ or $ abb_\bullet c_\bullet =a a_\bullet cc_\bullet$. Applying $\psi$, we easily deduce  that  in both cases $a=c$. The equality $\lambda_{0,1} (w)=0$ implies that $\tau(b)=b=c $. Hence  the triple $a,b,c$ is singular. 
 
 If $w $ is homotopic to a nanoword $ w^3$, then $\lambda_{1,0} (w)=0$ and $\lambda_{0,0} (w) $ is a monomial.
 The former implies that $a=b=\tau(b)$ and  the latter implies that $a=c $.   Hence  the triple $a,b,c$ is singular. 

If $w $ is homotopic to a nanoword $ w^4$, then  $\lambda_{0,0} (w) =1$.
 Hence $ abb_\bullet c_\bullet =a a_\bullet cc_\bullet$. Multiplying on the left by $\tau(a)$ and on the right by $(\tau(c))_\bullet$ we obtain  $bb_\bullet=a_\bullet c$.  Applying $\psi$, we easily deduce  that $a=b=c$. For $a=b=c$, the nanowords $w=ABACBC$ and $w^4_{a,b,c}=ABCBCA$
 are related by the   homotopy
 $$ ABACBC\mapsto  \underline {AD} E B \underline {AC} B E \underline {DC} 
 \mapsto DA\underline {EB} CA \underline {BE} CD \mapsto DACACD $$
 where $\vert D\vert =a$, $\vert E\vert =\tau (a)$ (the nanoword $DACACD$ is isomorphic to $ABCBCA$).
 
 It remains to show that  $w$ is not homotopic to a nanoword of length 4. Suppose first  that   $a\neq c$ or $c \neq \tau(b)$. If $w$ is   homotopic to a nanoword of length 4, then
 $\lambda_{1,1}(w)=0$ and hence $$\psi(\lambda_{1,1}(w))= z_a\otimes z_c -z_c \otimes z_c -z_a \otimes z_a+ z_a z_b z_c \otimes z_a z_b z_c=0.$$
 The term $z_a\otimes z_c$ must cancel with either $z_c \otimes z_c$ or $z_a \otimes z_a$. In both cases   $a=c$. The term $z_a z_b z_c \otimes z_a z_b z_c$ must then cancel with 
 $z_a \otimes z_a$. This is possible only if   $c=\tau (b)$ which contradicts our assumptions. 
 
Consider   the remaining  case $a=c=\tau(b)$. By the 5-regularity, $a\neq \tau(a)$.  We prove here the following weaker claim: if $w$ is homotopic to a nanoword $KLKL$, then $\vert K\vert =\vert L\vert =a$. That this is also impossible will be proven in Sect.\ \ref{ed56quaga}. To prove the weaker claim, 
set $k=\vert K\vert, l=\vert L\vert$. Then  $  \lambda_{1,0}(KLKL)=k-kk_\bullet l_\bullet$. Since $a=c=\tau(b)$, we have
$\lambda_{1,0}(w)=c b_\bullet c_\bullet  -ab a_\bullet  c c_\bullet=a-a a_\bullet a_\bullet$.
The equality $k-kk_\bullet l_\bullet = a-a a_\bullet a_\bullet$ is possible in only two cases: when $k=l=a$ and when $k=\tau(l), a=\tau(a)$. The latter is excluded by the 5-regularity. Therefore $k=l=a$. 
  \end{proof} 
   
    \begin{corol}\label{1edrcororlspsid}   For $i=1,..., 5$ and any $i$-regular triple $a, b,c\in   \alpha$, we have  $\vert \vert w^i_{a,b, c} \vert \vert =3$.      \end{corol}

 \begin{corol}\label{12cororlspsid}   A nanoword of length $\leq 6$ is homotopically (skew-) symmetric    if and only
 if   it is contractible  or (skew-) symmetric.      \end{corol}
 
One needs to prove only that a homotopically (skew-) symmetric  nanoword of length $\leq 6$     is contractible  or (skew-) symmetric.  This follows from the homotopy classification above.  Note that $w^i_{a,b,c}$ is   symmetric if and only
 if  $i\in \{1,5\}, a=c$  or $i=4, b=c$. The nanoword $w^i_{a,b,c}$ is   skew-symmetric if and only if   $i\in \{1,5\},  \tau(a)=c,  \tau(b)=b$  or 
 $i=4, \tau(a)=a, \tau(b)=c$.

\section{$\alpha$-quandles and  $\alpha$-keis}\label{4dd822}

 We  introduce $\alpha$-quandles and  $\alpha$-keis  generalizing the classical   quandles and keis.
 Their connections with words will be discussed in the next section.

\subsection{Keis and quandles}\label{kzeiikaga} Keis were introduced by M. Takasaki  in 1942 as   abstractions of symmetric
transformations. A more general notion of a quandle  was introduced in the 1980's by D. Joyce, S. Matveev, and E. Brieskorn independently,  see
\cite{ka}  for  a survey. A {\it quandle} is a non-empty set $X$ with   binary operation $ (x,y)\mapsto x\ast y  $ 
such that
   $x\ast x=x$,      $(x\ast y) \ast z= (x\ast z) \ast  (y\ast z)$,  for all $x, y, z\in X $, and  the mapping  $x\mapsto x\ast z:X\to X$ is a bijection  for any $z\in X$.
The identity  $(x\ast y) \ast z= (x\ast z) \ast  (y\ast z)$ means that  the bijection  $x\mapsto x\ast z$ preserves $\ast$. A quandle    is a {\it kei} if  this bijection 
   is   involutive for all $z$  that is  $(x\ast z)
\ast  z=x$ for all
$x, z\in X
$.   For example, $X=\ZZ/m\ZZ$ with $x\ast y=2y-x$ is a kei for any $m\geq 1$. 

\subsection{$\alpha$-quandles}\label{twisquaga} For a   set $\alpha$, we introduce   a notion of an $\alpha$-quandle. Let $X$ be a non-empty  set. Suppose
that  each
$a\in
\alpha$   gives rise to a bijection  $x\mapsto ax: X\to
X$ and  to a binary operation $ (x,y)\mapsto x\ast_a y  $  on $X$.  These operations form   an   {\it 
$\alpha$-quandle} if the following axioms are satisfied:

(i) $ax\ast_a x=x$ for all $a\in \alpha, x\in X$;

(ii) $a (x\ast_a y)= ax \ast_a ay$ for all $a\in \alpha, x, y\in X$;

(iii) $ (x\ast_a y)\ast_a z= (x \ast_a az) \ast_a (y\ast_a z)$ for all $a\in \alpha, x, y, z\in X$;

(iv)    the mapping  $x\mapsto x\ast_a z:X\to X$ is a bijection for any $a\in \alpha,  z\in X$.

  The operations $\{x\mapsto ax: X\to
X\}_{a\in \alpha}$ and     $ \{(x,y)\mapsto x\ast_a y\}_{a\in \alpha}  $ are called the {\it quandle operations}.  A morphism $X\to X'$  of $\alpha$-quandles is a
set-theoretic   map  $X\to X'$  commuting with these operations.

  To give   examples of $\alpha$-quandles, consider the  semi-group  $\tilde  \Psi $ with generators $\{{}a, a_{\bullet}   \}_{a\in
\alpha}$ subject to the relations
$ {}a a_{\bullet}=a_{\bullet} {}a $ for  all $a\in \alpha$.  Suppose that $\tilde  \Psi $ acts on a group $X$ by group automorphisms $x\mapsto ax$ and $x\mapsto
a_{\bullet} x$ where $x\in X$. Then the   automorphisms 
$x\mapsto {}a x
$   together with the binary operations
\begin{equation}\label{eeee}x\ast_a y= y (a_{\bullet}x) (a_{\bullet} {}ay)^{-1} \in X\end{equation} for
$x,y\in X$ form  an
$\alpha$-quandle.  We check the axioms.  Axiom  (i):
${}a x\ast_a x= x  (a_{\bullet} {}a x)  (a_{\bullet} {}ax)^{-1} =x$.  Axiom  (ii): $$ {}a(x\ast_a y)= a ( y (a_{\bullet}x) (a_{\bullet} {}ay)^{-1})=(ay) (aa_{\bullet} x) (a a_{\bullet} {} a y)^{-1}=
{}a x\ast_a  {}a y.$$  Axiom (iii):
$$  (x \ast_a {}a z) \ast_a (y\ast_a z)= (az) (a_{\bullet}x) (a_{\bullet} {}aaz)^{-1}\ast_a  z (a_{\bullet}y) (a_{\bullet} {}az)^{-1} $$
$$=
z (a_{\bullet}y) (a_{\bullet} {}az)^{-1}               (a_{\bullet} az) (a_{\bullet} a_{\bullet}x) (a_{\bullet} a_{\bullet} {}aaz)^{-1}   ((a_{\bullet} a z) (a_{\bullet} a a_{\bullet}y) (a_{\bullet} a a_{\bullet} {} az)^{-1})^{-1}      $$
$$= z (a_{\bullet}y) (a_{\bullet} {}az)^{-1}               (a_{\bullet} az) (a_{\bullet} a_{\bullet}x) (a_{\bullet} a_{\bullet} {}aaz)^{-1}  (a_{\bullet} a_{\bullet} {}aaz)  (a_{\bullet} a_{\bullet} ay)^{-1}     (a_{\bullet} a z)^{-1}  
$$
$$ =z (a_{\bullet}y)   (a_{\bullet} a_{\bullet}x)   (a_{\bullet} a_{\bullet} ay)^{-1}     (a_{\bullet} a z)^{-1}  
= z  \, a_{\bullet} (x\ast_a y) ({} a_{\bullet} az)^{-1}=  (x\ast_a y)\ast_a z.$$
Axiom (iv) follows from the assumption that $x\mapsto a_{\bullet}x:X \to X$ is a bijection.

In particular, any  left module over the semi-group ring $\ZZ\tilde \Psi $   is an $\alpha$-quandle with quandle  operations 
$x\mapsto
{}a x
$   and $x\ast_a y=  a_{\bullet}x+ (1-a_{\bullet} {}a)y $.

Note that if $\alpha=\{a\}$ is a 1-element set and $ax=x$ for all $x$, then an $\alpha$-quandle is simply a quandle.

\subsection{$\alpha$-keis}\label{keisga}  Consider a set $\alpha$ with involution $\tau$.  An {\it
$\alpha$-kei}   is an $\alpha$-quandle  $X$   such that 

(v) $a\tau (a) x=x$ for all  $x \in X, a\in \alpha$ and 

(vi) 
  $(x\ast_a y) \ast_{\tau(a)} ay=x$  for all  $x, y\in X, a\in \alpha$. 

These two axioms strengthen   Axiom (iv) above: they imply  that  for any $y\in X$, the mappings $x\mapsto 
x\ast_{a} y: X \to X $ and   $x\mapsto x
\ast_{\tau(a)} ay : X \to X $ are inverse to each other.  

The quandle operations $\{x\mapsto ax: X\to
X\}_{a\in \alpha}$ and     $ \{(x,y)\mapsto x\ast_a y\}_{a\in \alpha}  $ in a kei $X$ are called the {\it kei operations}.  Morphisms of $\alpha$-keis are their
morphisms as
$\alpha$-quandles.  Isomorphism of $\alpha$-keis is denoted by $\approx$.

Recall the ring $\Lambda=\ZZ\Psi$ defined in Sect.\ \ref{ringsfunct}. Any left $\Lambda$-module $X$ becomes  an $\alpha$-kei  with kei operations  
$x\mapsto {}a x
$   and $x\ast_a y= a_{\bullet}x + (1-a_{\bullet} {}a)y $ for 
$x,y\in X$. Axioms (i) -- (iv) were checked above, Axiom (v) follows from the definition of $\Lambda$ and Axiom (vi)
follows from the formulas 
$$(x\ast_a y) \ast_{\tau(a)} ay =\tau(a)_{\bullet} (x\ast_a y) + (1-\tau(a)_{\bullet} \tau(a)) ay$$
$$
=\tau(a)_{\bullet} a_{\bullet}x + \tau(a)_{\bullet} y - \tau(a)_{\bullet} a_{\bullet} {}a  y +ay - \tau(a)_{\bullet} \tau(a)  ay=x.$$
The $\alpha$-keis obtained by this construction from $\Lambda$-modules are said to be {\it abelian}.

   When
$\alpha$ is a 1-element set and $ax=x$ for all $x$,  the notion of an
$\alpha$-kei is equivalent to the one of a kei.

\subsection{Presentations by generators and relations}\label{keisalquaga}  For an $\alpha$-kei  $X$, a    set   $S\subset X$ {\it generates} $X$ if  all elements of $X$
can be obtained from  elements of $S$   using the  kei  operations.  Any  set $S$ generates a unique {\it  free $\alpha$-kei} $X(S)$ characterized by the
condition that  every set-theoretic map from $S$ to an
$\alpha$-kei $X$ extends  to a morphism of $\alpha$-keis $X(S)\to X$.  The elements of $X(S)$ are obtained from  elements of $S$   using the kei
operations modulo the identities imposed by Axioms (i) -- (iii), (v), (vi).  

As in the theory of groups,  we can  present an
$\alpha$-kei   by a set of generators  $S$  and a set  of   relations $R$.  A {\it relation} is a pair  $p,q$ of elements of $X(S)$ which we  
  write as an equality $p=q$.  The $\alpha$-kei presented by    $S$ and   $R$ is obtained  by
quotienting $X(S)$ by all the relations   from $R$ and all their corollaries.  For example, pick    
$a,b\in
\alpha$ and consider the $\alpha$-kei $X$ with generators  $x,y$ subject to one defining  relation $x\ast_b
b\tau (a) y=(bx\ast_{a} y)\ast_{\tau (b)} aax$. This means that every element of $X$ can be obtained from $x,y$   using the kei
operations, that  the relation in question is satisfied, and that all other relations between $x$ and $y$ are
corollaries of this one and the axioms of an
$\alpha$-kei.  Following a similar train of ideas, one can define   free products of  $\alpha$-keis.  We leave the details to the reader, cf. \cite{fr} for presentations
of quandles.

 \subsection{Abelianization}\label{abelquaga} Each $\alpha$-kei $X$ gives rise to a  $\Lambda$-module   with generators   $\{[x]\}_{x\in X}$â  subject to the
\lq\lq commutation relations" 
$[ax]=a [x]$ and $ [x\ast_a y]= a_{\bullet} [x]+  (1- a_{\bullet} {}a) [y]$  for all $ x, y\in X, a\in \alpha$.  This module  is called the {\it abelianization} of $X$
and denoted $X^{ab}$.  It can be described by the following universal property:  any kei morphism from $X$ to  an  abelian $\alpha$-kei $Y$ factors uniquely as a
composition of the mapping $x\mapsto [x]: X\to X^{ab}$ and a $\Lambda$-homomorphism $X^{ab}\to Y$.

 Given a presentation of $X$ by generators $S$ and
relations
$R$ we can   compute
$X^{ab}$ in terms of generators and relations. This module is generated by the set $\{[s]\}_{s\in  S}$ and each relation $p=q$ from   $R$ gives rise to a relation
obtained by expanding both $p$ and
$q$ as linear combinations of the vectors $\{[s]\}_{s\in  S}$  via the commutation relations.  In particular, for the  free $\alpha$-kei $X(S)$ generated by $S$,
we have  $X(S)^{ab}=\oplus_{s\in S} \Lambda [s]$. For each $v\in X(S)$, we have a unique expansion $[v]=\sum_{s\in S}  \lambda_s(v) [s]$ with $\lambda_s(v)\in
\Lambda$. 

\subsection{Marked $\alpha$-keis}\label{frararaelquaga} An $\alpha$-kei  $X$ endowed with an ordered pair of distinguished elements   $v_-, v_+\in X$ is  {\it
marked}.  The  elements $v_-, v_+$ are called   the  {\it input} and the  {\it output}, respectively. By  {\it  (iso)morphisms}  of marked $\alpha$-keis, we mean kei
(iso)morphisms preserving the input and the output.
 The abelianization of a marked $\alpha$-kei $(X, v_-  ,v_+  )$ is the marked  $\Lambda$-module $(X^{ab}, [v_-]\in X^{ab}, [v_+]\in X^{ab})$.

Given a marked $\alpha$-kei  $X $ we define a marked $\alpha$-kei $\overline X$ to be the same  set $X$ with the same input and output and new kei operations  
$ax: = \tau (a) x$, $x\ast_a y : =x \ast_{\tau  (a)} y$ for $x,y\in X, a\in \alpha$.   Clearly, $\overline  {\overline  X}=X$.  A presentation of $X$ by generators $S$ and
relations
$R$ yields a presentation of $\overline X$ by generators $S$ and
relations
$\overline R$ where the relations in $\overline R$ are obtained from those in $R$ by replacing each letter $a\in \alpha$ appearing in these relations by $\tau(a)$.
 We   define a marked $\alpha$-kei $X^-$ to be $\overline X$ with permuted input and output. Clearly, 
 $(X^-)^-=X$.

 We point out  a simple 
 numerical invariant of a finitely generated marked $\alpha$-kei $ \K$.  Fix  a marked
$\alpha$-kei $X$ that is finite as a set. Then the set of  marked  morphisms   $ \K   \to X$ 
  is finite.  The number  of  such  morphisms    is an isomorphism invariant of $\K$.

\subsection{Reconstruction}\label{reconsssaga}   Given an $\alpha$-quandle  $X$ and a  set $\alpha_0\subset \alpha$, we can restrict the quandle operations  
 in $X$ to only those $a$ which belong to $\alpha_0$. This gives an $\alpha_0$-quandle $Y$ coinciding with $X$ as a set and called the {\it
restriction} of $X$ to $\alpha_0$.  If $X$ is an $\alpha$-kei, then for any 
$a\in
\alpha_0 \cap \tau (\alpha_0)$ and any $x, y\in Y$, we have
$a\tau (a) x=x$ and 
  $(x\ast_a y) \ast_{\tau(a)} ay=x$.  An $\alpha_0$-quandle $Y$ satisfying the latter conditions is  said to be {\it compatible} with $\tau$.

\begin{lemma}\label{recoccsovl} Let $\alpha_0$ be a subset of
$\alpha$ such that
$\alpha_0 \cup \tau (\alpha_0)=\alpha$.   Let $Y$ be an 
$\alpha_0$-quandle compatible with $\tau$. Then 
 there is a unique $\alpha$-kei $X$ whose restriction to $\alpha_0$ is   $Y$. 
\end{lemma} 
                     \begin{proof}  Clearly, $X=Y$ as a set and an $\alpha_0$-quandle.  We need only to define the operations $x\mapsto ax$ and $(x,y)\mapsto x\ast_a
y$ for  
$a\in
\alpha-\alpha_0$.  By   assumption,  
  $\tau(a)\in \alpha_0$. By the definition of an $\alpha_0$-quandle,  the mappings $x\mapsto \tau (a)x:X\to X$ and $x\mapsto x \ast_{\tau(a)} ay:X\to X$ are bijective. The axioms of an $\alpha$-kei
show that we must take their inverses as  the mappings $x\mapsto ax$ and $x\mapsto x\ast_a y$, respectively.  Thus
$z=ax$ is the only element   of
$X$ such that
$\tau(a)z =x$ and $t=x\ast_a y$ is the only element   of $X$ such that $t \ast_{\tau(a)} ay=x$. We  show  that this  makes $X$ into an $\alpha$-kei. 

We first  check the axioms of  an $\alpha$-quandle. It suffices to check them for  $a\in \alpha-\alpha_0$. Axiom (i): to verify that  $ax\ast_a x=x$ it suffices to check that 
$x
\ast_{\tau(a)} ax= ax$.  This holds since $$x \ast_{\tau(a)} ax= a \tau (a) (x \ast_{\tau(a)} ax)=a (\tau (a) x  \ast_{\tau(a)} \tau(a) ax) = a (\tau (a) x 
\ast_{\tau(a)} x) =ax.$$ To check (ii), set $t= x\ast_{a} y$. The  required equality  $a t= ax \ast_{a} a y$
is equivalent to $a t \ast_{\tau(a)}  aay= ax$. We have
$$a t \ast_{\tau(a)}  aay= a\tau (a) (a t \ast_{\tau(a)}  aay)=a( \tau (a)  a t \ast_{\tau(a)}  \tau (a) aay)= a(t \ast_{\tau(a)} ay)=ax$$
by the choice of $t$.  To check (iii), set $t=(x\ast_a y)\ast_a z$.  The  required equality  $t=(x \ast_a az) \ast_a (y\ast_a z)$ is equivalent to
$t \ast_{\tau(a)} a (y\ast_a z)=x \ast_a az$. The latter is equivalent to
$$(t \ast_{\tau(a)} a (y\ast_a z)) \ast_{\tau(a)} aaz=x.$$  By Axiom (iii) applied to $\tau (a)\in \alpha_0$, 
$$(t \ast_{\tau(a)} a (y\ast_a z)) \ast_{\tau(a)} aaz=(t \ast_{\tau(a)}  \tau(a)aaz)  \ast_{\tau(a)} (a (y\ast_a z)  \ast_{\tau(a)} aaz)$$
$$=(t \ast_{\tau(a)}   az)  \ast_{\tau(a)} (  (ay\ast_a az)  \ast_{\tau(a)} aaz)= (x\ast_a y)  \ast_{\tau(a)}  ay=x.$$
Axiom (iv) follows directly  from the definitions.  

To see that $X$ is an $\alpha$-kei, we must check that 
$(x\ast_a y) \ast_{\tau(a)} ay=x$  for all  $x, y\in X, a\in \alpha$.  This follows from the definition of $\ast_a$ for $a\in \alpha-\alpha_0$
and from the compatibility assumption for $a\in \alpha_0 \cap \tau (\alpha_0)$. It remains to consider the case where $a\in \alpha_0$ and $\tau (a)\notin
\alpha_0$.  By the definition of $\ast_{\tau(a)}$, we know that $t=(x\ast_a y) \ast_{\tau(a)} ay$ is the only element of $X$ such that $t \ast_a \tau(a) ay=x\ast_a
y$. Clearly,
$t=x$ is such an element. Therefore $(x\ast_a y) \ast_{\tau(a)} ay=x$.
 \end{proof} 
 
This lemma  establishes   a bijective correspondence between (the isomorphisms classes of)  $\alpha$-keis and $\alpha_0$-quandles compatible with $\tau$.  If $\alpha_0\cap \tau(\alpha_0)=\emptyset$, then the compatibility condition   is empty and we obtain a bijective correspondence between  $\alpha$-keis
and $\alpha_0$-quandles.

\section{Keis of nanowords}\label{56822}

 With each nanoword $w$ over $\alpha$ and a $\tau$-invariant  set $\beta\subset \alpha$ we associate a marked  $\alpha$-kei   $ \K_\beta (w)$. It is  preserved  under homotopy moves on 
$w$.

\subsection{The $\alpha$-kei  $ \K_\beta (w)$}\label{gikl12keisalquaga} Consider a nanoword $(\A, w:\hat n \to \A)$ over
$\alpha$.  As in Sect.\ \ref{nota1}, for  a letter
$A\in
\A$, we denote by $i_A$ (resp.\ $j_A$) the minimal (resp.\ the maximal) element of the 2-element set $\omega^{-1} (A)\subset \hat n$.
The $\alpha$-kei   $\K_\beta (w) $ is  generated 
by $n+1$ symbols $x_0, x_1,..., x_{n}$ subject to the following $n$ defining relations:  for any   $A\in \A$ such that   $a=\vert A\vert \in \beta$,
$$x_{i_A}=a \, x_{i_A-1}, \,\,\,\,\,x_{j_A}=x_{j_A-1} \ast_{a} x_{i_A-1},$$
and  
for any     $A\in \A$ such that   $a=\vert A\vert \in  \alpha- \beta$,
$$x_{i_A }=x_{i_A-1}\ast_{a} x_{j_A-1} ,  \,\,\,\,\,x_{j_A}=a  \,x_{j_A-1} .$$ 
  We take  $x_0$ as the input and  $x_n$ as the output in $\K_\beta (w)$.  Comparing this definition with the one of the marked module $K_\beta (w)$ we obtain that   $$K_\beta (w)=(\K_\beta(w))^{ab}.$$

It is clear that if $w=w_1w_2$ is a product of two nanowords $w_1,w_2$, then $\K(w)$ is the free product of $\K(w_1)$ and $\K(w_2)$ quotiented by the
relation   (the output of $\K(w_1))$ $=$ (the input of $\K(w_2)$). 

It follows   from the definitions that    $\K_\beta (\overline w)\approx \overline {\K_{ \beta} (w)}$.  The   isomorphism is   the identity  $x_i\mapsto x_{ i}$ on the generators.  Similarly,  
$\K_\beta (w^-)\approx   (\K_{\alpha-\beta} (w))^-$.  The latter isomorphism is induced by the permutation of the generators $x_i\mapsto x_{n-i}$.  
 The main property of $ \K_\beta (w)$ is
contained in the next lemma.

\begin{lemma}\label{efrtcovl}  The isomorphism type of the marked $\alpha$-kei $ \K_\beta (w) $ is a homotopy invariant of $w$.
\end{lemma} 
                     \begin{proof}    If $w$ is replaced by an isomorphic nanoword, then    $ \K_\beta (w) $ is replaced by an isomorphic marked $\alpha$-kei. Consider the first homotopy move $  w=xAAy \mapsto xy=v$. Set $a=\vert A\vert\in \alpha$ and let $i=i_A, i+1=j_A$ be  the
 indices numerating the   entries of
$A$  in
$w$.  If $a \in \beta$, then   the generator $x_i$ of $\K_\beta (w)$  is involved in    two relations    $x_i=a x_{i-1}$ and $x_{i+1}=x_i \ast_a x_{i-1}$. We can
exclude
  $x_i$ from the set of generators using the first relation and replace the  second relation with $x_{i+1}=a x_{i-1} \ast_a x_{i-1}$. By Axiom (i), this   may be rewritten  as 
  $x_{i+1}=x_{i-1}$. Thus we can further remove  $x_{i+1}$  from the set of generators and   replace it   in  the remaining relations by  $x_{i-1}$. The resulting set of generators and relations is a
presentation of $\K_\beta (v) $.

If $a \in \alpha-  \beta$, then  
\begin{equation}\label{darz} \K_\beta (w)\approx (\K_{\alpha-\beta} (w^-))^-\approx  (\K_{\alpha-\beta} (v^-))^- \approx \K_\beta (v)\end{equation}  where  the second isomorphism
follows from the inclusion  $a \in \alpha-  \beta$ and the result of the previous paragraph.   For  completeness,  we give a  direct proof of  the
isomorphism   $\K_\beta (w)\approx  \K_\beta (v)$. The generator $x_i$ of $\K_\beta (w)$  is involved in   
two relations   
$  x_{i+1}=ax_i$ and
$x_{i}=x_{i-1} \ast_a x_{i}$.  The first formula is equivalent to $x_i=
\tau (a) x_{i+1}$.  We can exclude
$x_i$ via this relation and rewrite the second relation as  $\tau (a) x_{i+1}=x_{i-1} \ast_a  \tau (a) x_{i+1}$. Note that
for any  elements
$x,y$ of an
$\alpha$-kei, $$\tau (a) x= y \ast_a \tau (a) x \Longleftrightarrow  x= a(y \ast_a \tau (a) x)\Longleftrightarrow   x= ay \ast_a x\Longleftrightarrow  x=y.$$ Therefore we can rewrite the relation $\tau (a) x_{i+1}=x_{i-1} \ast_a  \tau (a) x_{i+1}$ as  $ x_{i+1} =x_{i-1}$. 
 Thus we can further remove  $x_{i+1}$  from the set of generators and   replace it   in  the remaining relations by  $x_{i-1}$. The resulting set of generators and relations is a
presentation of $\K_\beta (v) $.

Consider the second homotopy move $  w=xAByBAz \mapsto xyz=v$. Set $a=\vert A\vert\in \alpha$ and let $i=i_A,j=j_A\geq i+3$ be  the 
indices numerating the   entries of
$A$  in
$w$.  By assumption, $\vert B\vert=\tau (a)$.   If $a \in \beta$, then  $\tau (a)\in \beta$ and  the letters $A,B$ give rise to four relations    $x_i=a x_{i-1},
x_{i+1}={\tau (a)} x_i$  and
$x_{j-1}=x_{j-2} \ast_{\tau (a)} x_{i},  x_{j}=x_{j-1} \ast_a x_{i-1}$. The first two relations imply that $ x_{i+1}=   x_{i-1}$. We exclude
$x_i, x_{i+1}$ from the set of generators replacing them  in all other relations with $a x_{i-1}, x_{i-1}$ respectively. The resulting relation  $x_{j-1}=x_{j-2} \ast_{\tau (a)}
a x_{i-1}$ allows us to exclude $x_{j-1}$ from the set of generators and to rewrite the relation $  x_{j}=x_{j-1} \ast_a x_{i-1}$ as  $x_{j}=(x_{j-2}
\ast_{\tau (a)} a x_{i-1})
\ast_a    x_{i-1}$. The latter is equivalent to  $x_{j}= x_{j-2}$. Thus we can further remove  $x_{j}$  from the set of generators  and   replace it   in the remaining  relations by  $x_{j-2}$. The resulting
set of generators and relations is a presentation of $\K_\beta (v)$. The case $a \in \alpha-\beta$ follows from the   case $a \in  \beta$
by (\ref{darz}).

Consider the third homotopy move $  w=xAByACzBCt \mapsto xBAyCAzCBt =v$. Set $a=\vert A\vert=\vert B\vert=\vert C\vert\in \alpha$. Let
$i=i_A,j=j_A, k=j_B$ be  the  indices numerating the   entries of
$A$  in
$w$ and the second entry of $B$ in $w$.     If $a \in \beta$, then    the letters $A,B,C$ give rise to six
relations   in $\K_\beta(w)$, namely, $x_i=a x_{i-1}$ and 
$$  x_{i+1}=a x_i, 
 x_{j}=x_{j-1} \ast_{a} x_{i-1},  x_{j+1}=a x_{j}, x_k=x_{k-1} \ast_a x_i, x_{k+1} =x_{k} \ast_a x_j.$$  We   exclude  $x_i,x_j,x_k$ using the first, third, and
fifth relations rewriting the other  3 relations as
\begin{equation}\label{dfdfdfd} x_{i+1}=aa x_{i-1}, \, x_{j+1}=a (x_{j-1} \ast_{a} x_{i-1}), \,\end{equation}
$$ x_{k+1} =(x_{k-1} \ast_a ax_{i-1})\ast_a (x_{j-1} \ast_{a}
x_{i-1}).$$
Similarly,  
  $A,B,C$ give rise to 6
relations   in $\K_\beta(v)$, namely, $x_i=a x_{i-1}$ and 
$$ x_{i+1}=a x_i, 
 x_{j}=a x_{j-1} ,  x_{j+1}=  x_{j}\ast_a  x_{i}, x_k=x_{k-1} \ast_a x_{j-1}, x_{k+1} =x_k \ast_a x_{i-1}.$$  We   exclude  $x_i,x_j,x_k$
using the first, third, and fifth relations rewriting the other  3 relations as
$$ x_{i+1}=aa x_{i-1}, \, x_{j+1}=a  x_{j-1} \ast_{a} ax_{i-1} , \, x_{k+1} =(x_{k-1} \ast_a  x_{j-1})\ast_a  x_{i-1}.$$
These relations are equivalent to the relations (\ref{dfdfdfd}) modulo the axioms of  an $\alpha$-kei. 
The case $a \in \alpha-\beta$ follows from the   case $a \in  \beta$
by (\ref{darz}). \end{proof}

\subsection{The invariant $v_+(w)$}\label{betaemvakn1}  For  $\beta=\alpha$, the   defining relations  of the $\alpha$-kei $\K_\beta (w)=\K_\alpha (w)$ consecutively express the generators via the previous ones. Therefore $\K_\alpha (w)$ is a free $\alpha$-kei generated by the input
$v_-$.  We identify  $\K_\alpha (w)$  with the    free $\alpha$-kei $X(s)$ with one generator $s$ via $v_-=s$.  The output $ v_+(w) =v_+ \in \K_\alpha (w)=X(s)$
is   a homotopy invariant of
$w$. It  includes the polynomials $\lambda'(w), \lambda(w)\in \Lambda$ studied  in Sect.\ \ref{spebetnct} via 
  $[v_+(w)]=\lambda' (w) [s] \in (X(s))^{ab}= \Lambda [s]$.  
  
Since $X(s)$ is a free $\alpha$-kei,  the element $v_+(w) \in X(s)$ determines a     unique kei morphism 
$\phi (w):X(s) \to  X(s) $ such that $\phi (w)  (s)=v_+(w) $.  This  morphism is a  homotopy invariant  of $w$.  If $w=w_1w_2$ is a product of two
nanowords, then $\phi(w) =\phi(w_1)\,  \phi(w_2)$. Indeed,  both sides are kei endomorphisms of $X(s)$ transforming   $s$ into   $v_+(w) $. 

We finish this section by showing how to compute $v_+(\overline w)$ from $v_+(w)$. Given two $\alpha$-keis $X, X'$, we call a mapping $f:X\to X'$    a  {\it kei anti-morphism}  if $f(ax)=\tau (a ) x$ and $f(x\ast_a y)= f(x) \ast_{\tau(a)} f(y)$ for all $x,y\in X, a \in \alpha$. In other words, a kei  anti-morphism $X\to X'$ is a kei morphism $X\to \overline {X'}$. It is obvious that $\overline {X(s)}$ is a free $\alpha$-kei  generated by $s$ and therefore there is a unique kei anti-morphism $\iota:X(s) \to X(s)$ preserving $s$.   Clearly, $\iota$ is an involution. It follows   from the definitions that
$v_+(\overline w)= \iota (v_+(w))$.

 \section{Case of free $\tau$}\label{dfdfdfdfd}

Throughout this section  we suppose that    $\tau$ is
fixed-point-free and fix a set   $\alpha_0\subset \alpha$ meeting each orbit of $\tau$ in one element.
We  derive    from  $\alpha$-keis   of   nanowords   a  simpler invariant,  called the characteristic sequence.   It is used to accomplish the homotopy classification of nanowords of length $\leq 6$.

\subsection{The $\alpha$-kei $F$.}\label{keisfreetauquaga}   We  first   construct an  $\alpha$-kei $F$ needed for the sequel.       Recall from Sect.\ \ref{ringsfunct} the group 
$\Psi=
\Psi_\alpha$  with generators 
$\{{}a, a_{\bullet}  
\}_{a\in
\alpha}$ and     defining relations
$ {}a a_{\bullet}=a_{\bullet} {}a,  a \tau (a)=a_{\bullet} \tau(a)_{\bullet}=1   $ for   $a\in \alpha$.      Let $F$ be  the free group (of countable rank)
  freely generated by the set $ \Psi$.  The generator of $F$ corresponding to $\psi\in \Psi$ is denoted $\underline \psi$.    Note that 
$\underline \psi$ and $  \underline {\psi^{-1}}$ are two independent  generators of $F$ for $\psi\neq 1$ and $\underline 1$ is a  non-trivial generator of $F$
where $1$ is the unit of $\Psi$.  A typical element of $F$ has the form
$(\underline {\psi_1})^{\varepsilon_1}\cdots  (\underline {\psi_m})^{\varepsilon_m}$ where $m\geq 0$, ${\psi_1}, ..., {\psi_m} \in \Psi$, and $\varepsilon_1,...,\varepsilon_m\in
\{\pm  1\}$.  

The left action of
$\Psi$ on itself extends   to a group action of $\Psi$ on $F$.  In particular, the generators  $a, {}a_{\bullet}\in \Psi$ act on $F$
by the group automorphisms 
$$ a  \left ( (\underline {\psi_1})^{\varepsilon_1}\cdots  (\underline {\psi_m})^{\varepsilon_m}  \right )= (\underline { a {\psi_1}})^{\varepsilon_1}\cdots  (\underline
{a{\psi_m}})^{\varepsilon_m},$$
$$ {}a_{\bullet}  \left ( (\underline {\psi_1})^{\varepsilon_1}\cdots  (\underline {\psi_m})^{\varepsilon_m}\right )= (\underline {{}a_{\bullet} {\psi_1}})^{\varepsilon_1}\cdots 
(\underline {{}a_{\bullet}{\psi_m}})^{\varepsilon_m}.$$  By   Sect.\  \ref{keisga}, this  induces a structure of an $\alpha_0$-quandle on $F$  with quandle  operations 
$x\mapsto {}ax$ and
$x\ast_a y= y (a_{\bullet} x) (a_{\bullet} {}ay)^{-1}  $ for $a\in \alpha_0$. By Lemma \ref{recoccsovl}, this   extends uniquely to an $\alpha$-kei structure on $F$. The compatibility condition  in Lemma \ref{recoccsovl}   is empty since $\alpha_0  \cap \tau (\alpha_0)=\emptyset$. Note that the resulting $\alpha$-kei structure on $F$, generally speaking, depends on the choice of $\alpha_0$.  The next lemma gives explicit  formulas for the kei operations in $F$ determined by      elements of  $  \alpha-\alpha_0$.

\begin{lemma}\label{ekkirzzd}   
For $a\in \alpha_0$, ${\psi_1}, ..., {\psi_m} \in \Psi$, and $\varepsilon_1,...,\varepsilon_m\in
\{\pm  1\}$, 
$$\tau(a)  \left ( (\underline {\psi_1})^{\varepsilon_1}\cdots  (\underline {\psi_m})^{\varepsilon_m}\right )=  (\underline {a^{-1} {\psi_1}})^{\varepsilon_1}\cdots 
(\underline {a^{-1} {\psi_m}})^{\varepsilon_m}=  (\underline {\tau(a) {\psi_1}})^{\varepsilon_1}\cdots 
(\underline {\tau(a) {\psi_m}})^{\varepsilon_m}.$$
For $a\in \alpha_0$, $x,y\in F$, 
\begin{equation}\label{609} x\ast_{\tau(a)} y= (a_{\bullet}^{-1} a^{-1} y)^{-1} \, ( a_{\bullet}^{-1}  x)\, y \in F.\end{equation}
\end{lemma}
    \begin{proof} Recall from the proof of Lemma \ref{recoccsovl} that    $z=\tau(a) x\in F$    and
$t=x\ast_{\tau(a)} y$ are uniquely determined from the equations
$a z=x$    and  $t\ast_{a} \tau(a) y=x$, respectively. Observe that
$$a \left ((\underline {\tau(a) {\psi_1}})^{\varepsilon_1}\cdots 
(\underline {\tau(a) {\psi_m}})^{\varepsilon_m}\right )=(\underline {a\tau(a) {\psi_1}})^{\varepsilon_1}\cdots 
(\underline {a\tau(a) {\psi_m}})^{\varepsilon_m}= (\underline {  {\psi_1}})^{\varepsilon_1}\cdots 
(\underline {  {\psi_m}})^{\varepsilon_m}.$$
This implies  the first claim of the lemma. We have 
$$ (a_{\bullet}^{-1} a^{-1} y)^{-1} \, ( a_{\bullet}^{-1}  x)\, y \ast_{a} \tau(a) y
= \tau(a) y \, a_\bullet ((a_{\bullet}^{-1} a^{-1} y)^{-1} \, ( a_{\bullet}^{-1}  x)\, y)\, (a_\bullet a \tau(a) y)^{-1}$$
$$= \tau(a) y \,   (a_\bullet a_{\bullet}^{-1} a^{-1} y)^{-1} \, ( a_\bullet a_{\bullet}^{-1}  x)\, a_\bullet y \, (a_\bullet   y)^{-1} 
= \tau(a) y \,   ( a^{-1} y)^{-1} \, x  =x$$
where we use that $a\tau(a) y=y$ and $a^{-1} y = \tau(a) y$. This implies   (\ref{609}).
\end{proof}

\subsection{Characteristic sequences.}\label{charistquaga}  We  can use the $\alpha$-kei  $F$ constructed above  to  study the free $\alpha$-kei  $X(S)$
generated by a set
$S$.   By the definition of $X(S)$, there is a unique  kei morphism $f:X(S)\to F$ sending  $S$ to $\underline 1\in F$.  
For  $v\in X(s)$,    we have $f(v)=(\underline {\psi_1})^{\varepsilon_1}\cdots  (\underline
{\psi_m})^{\varepsilon_m}$  for some 
$\varepsilon_1,..., \varepsilon_m \in \{\pm 1\}$ and   ${\psi_1},..., {{\psi_m}} \in \Psi$.  The sequence $ (\varepsilon_1 \psi_1 ,...,
  \varepsilon_m  {{\psi_m}}) $ is called
a
{\it characteristic sequence} of
$v$. It  is well defined up to insertion/deletion   of  pairs of consecutive terms  $+\psi, -\psi$ or $-\psi, +\psi$   
with     $\psi \in \Psi$.    Deleting  all such pairs     we obtain a unique {\it reduced  characteristic sequence} of $v$.   It   can be efficiently computed. For example, let us compute the reduced characteristic sequence of $v =(ba s\ast_a s)
\ast_{b} a s $ where $a,b\in \alpha, s\in S$.    By definition  
$f(s)=\underline 1$,   
$$f(bas)=ba \underline {1}=\underline {ba},\,\,\,  f(bas\ast_a s) = f(bas) \ast_a f(s)= \underline {ba} \ast_a \underline 1=   \underline 1 \, \underline  {a_{\bullet} ba}\,  (\underline {a_{\bullet} a})^{-1}  ,$$
$$ f((bas\ast_a s) \ast_{b} as)= f(bas\ast_a s) \ast_b f(as)$$
$$= (\underline 1 \, \underline  {a_{\bullet} ba}\,  (\underline {a_{\bullet} a})^{-1}  ) \ast_b \underline a= \underline a
\,  \,
 \underline {b_{\bullet}} \, \, \underline  {b_{\bullet}a_{\bullet} ba}\,  (\underline {b_{\bullet} a_{\bullet} a})^{-1}    (\underline {b_{\bullet} ba})^{-1}  .$$
 Hence $(  a ,
 {b_{\bullet}}  ,    {b_{\bullet}a_{\bullet} ba} ,  -   {b_{\bullet} a_{\bullet} a} , -   {b_{\bullet} ba}  )$ is a  characteristic sequence of $v$. It  is reduced since $b\neq 1\in \Psi$.
 
 The following lemma relates  the  characteristic sequences to the polynomials
$\lambda_s(v)\in \Lambda=\ZZ \Psi$ defined in Sect.\ \ref{abelquaga}.

\begin{lemma}\label{nezzoccsovl}   Let $ (\varepsilon_1 \psi_1 ,...,
  \varepsilon_m  {{\psi_m}}) $ be 
a
 characteristic sequence  of
$v \in X(S)$. Then $$ \varepsilon_1 \psi_1 +...+
  \varepsilon_m  {{\psi_m}}=\sum_{s\in S} \lambda_s(v).$$
\end{lemma} 
                     \begin{proof}  The canonical   inclusion    $\Psi\hookrightarrow \ZZ\Psi=\Lambda$  extends uniquely to a
group homomorphism $i$ from $F$ to the underlying additive group of $ \Lambda$. In particular,
$$i(f(v))=i ((\underline {\psi_1})^{\varepsilon_1}\cdots  (\underline {\psi_m})^{\varepsilon_m})=\varepsilon_1 \psi_1 +...+
  \varepsilon_m  {{\psi_m}}.$$
The ring $\Lambda$ being a left $\Lambda$-module acquires the structure of  an abelian $\alpha$-kei. Comparing the definitions of the kei operations in $F$ and
$\Lambda$ we observe that $i$ is a kei morphism.  The mapping  $if:X(S)\to \Lambda$ splits therefore as a composition of the projection $x\mapsto [x]: X(S)\to
X(S)^{ab}$ and a $\Lambda$-homomorphism $f': X(S)^{ab} \to \Lambda$. For any $s\in S$, we have $f'([s])=if(s)=i(\underline 1)=1$. Therefore 
$$\varepsilon_1 \psi_1 +...+
  \varepsilon_m  {{\psi_m}}= i(f(v))=f'([v])= f' (\sum_{s\in S} \lambda_s(v) [s])=\sum_{s\in S} \lambda_s(v).$$
\end{proof}

\begin{lemma}\label{g5g5nezzoccsovl} Let $\tau_\#:\Psi \to \Psi$ be the involutive  group automorphism  sending the generators $a, a_\bullet$ to $\tau(a)=a^{-1}, \tau(a)_\bullet=a_\bullet^{-1}$, respectively, for   $a\in \alpha$.    If $ (\varepsilon_1 \psi_1 ,...,
  \varepsilon_m  {{\psi_m}}) $ is 
 a
 characteristic sequence  of   $v\in X(S)$, then 
$ (\varepsilon_m \tau_\# (\psi_m)  ,$...,
  $\varepsilon_1   \tau_\#( {\psi_1})) $ is 
 a
 characteristic sequence  of $\iota(v) $ where $\iota$ is    the unique kei anti-automorphism of $X(S) $ preserving $S$ element-wise. \end{lemma} 
                     \begin{proof}  Denote by    $\tilde \tau $ the  unique group anti-automorphism  of $F$ extending $\tau_\#:\Psi \to \Psi$. Observe  that  $\tilde \tau(ax) =a^{-1} \tilde \tau (x)$     and $\tilde \tau(a_\bullet x)=   a_\bullet^{-1} \tilde \tau (x)$  for   $a\in \alpha, x\in F$. Moreover,  $\tilde \tau$ is  a kei anti-morphism. Indeed,  for $x, y\in F$, $a\in \alpha_0$,  
                     $$ \tilde \tau (x\ast_a y)= \tilde \tau ( y (a_\bullet x) (a_\bullet a y)^{-1})=
                     ( \tilde \tau (   a_\bullet a y ) )^{-1}    \,   \tilde \tau (  a_\bullet x)\, \tilde \tau (y)$$
                      $$=
                      (a_{\bullet}^{-1} a^{-1} \tilde \tau (y))^{-1}  \,   a_\bullet^{-1}  \tilde \tau (    x) \, \tilde \tau (y)
                      =\tilde \tau (x) \ast_{\tau(a) } \tilde \tau(y) .$$
                      The case $a\in \alpha-\alpha_0$ is similar.  The compositions $f\iota:X(S)\to F$ and
                      $\tilde \tau f: X(S) \to F$ are kei  anti-morphisms sending $S$ to $\underline 1$.   Therefore $f\iota= \tilde \tau f$. The equality  $f \iota( v)= \tilde \tau  f(v) $  implies the claim of the lemma.          \end{proof}

\subsection{Characteristic sequences of   nanowords.}\label{charnananaquaga}    By a {\it characteristic sequence} of a nanoword $w$,  we mean
a characteristic sequence of the element 
$v_+(w)\in X(s)$ defined in   Sect.\ \ref{betaemvakn1}. This sequence, say   $ (\varepsilon_1 \psi_1 ,...,
  \varepsilon_m  {{\psi_m}}) $,  is well defined up to insertion   of  consecutive terms  $+\psi, -\psi$ or $-\psi, +\psi$   
with     $\psi \in \Psi$. Deleting  all such terms we obtain the  {\it reduced  characteristic sequence} of $w$.  Lemma \ref{nezzoccsovl}  implies that 
$\lambda'(w)= \varepsilon_1 \psi_1 +...+
  \varepsilon_m  {{\psi_m}}$.  Lemma \ref{g5g5nezzoccsovl} implies that  $ (\varepsilon_m \tau_\# (\psi_m) ,...,
  \varepsilon_1   \tau_\#( {\psi_1})) $ is 
 a
 characteristic sequence  of $\overline w$.  By Lemma \ref{efrtcovl}, the reduced  characteristic sequence of a nanoword  is a homotopy invariant. Note that the $\alpha$-kei structure on  $F$ and hence the characteristic
sequences of nanowords   depend on the choice of $\alpha_0$.

\subsection{Examples.}\label{cexcccnaquaga} 1. Let $w=ABAB$ with $\vert A\vert =a\in \alpha_0,  \vert B\vert =b\in \alpha_0$. The $\alpha$-kei $\K_\alpha(w)$ is generated by $v_-=x_0,x_1, ..., x_4=v_+$ subject 
to the relations $x_1=a x_0, x_2=b x_1, x_3=x_2 \ast_a x_0, x_4= x_3 \ast_b x_1$.  
Thus $v_+  =(ba v_-\ast_a v_-)
\ast_{b} a v_-$ and   $v_+(w)=(ba s\ast_a s)
\ast_{b} a s\in X(s)$.  By the computations above,  the reduced characteristic sequence of $w$ is
$(  a ,
 {b_{\bullet}}  ,    {b_{\bullet}a_{\bullet} ba} ,  -   {b_{\bullet} a_{\bullet} a} , -   {b_{\bullet} ba}  )$.  
 
 2. Pick four  letters $a,b,c,d\in \alpha$ (possibly coinciding) and consider  the nanoword $w=ABCDCDAB$ with $\vert A\vert =  a , \vert B\vert =b , \vert C\vert =  c , \vert D\vert =d$.   
 An  inspection    shows that  if $a\neq \tau(b)$ and  $c\neq \tau(d)$, then  the  $\alpha$-pairing of $w$  is primitive. Then  $\vert \vert w\vert \vert =4$ and  $w$ is  non-contractible.  However, for any $a, b \in \alpha_0$ and $c=\tau (b) , d=\tau (a)$ we have $f(v_+(w))=\underline 1$ (a direct computation).
In this case,   the   characteristic sequence does not distinguish $w$ from a contractible nanoword.
This   shows that the   characteristic sequence is not a faithful invariant of nanowords and does not allow to recover the associated $\alpha$-pairing.

\subsection{End of the proof of Theorem \ref{1classid}.}\label{ed56quaga}   We must show that the nanowords $  ABACBC$
and    $ACAC$ with  $\vert A\vert =\vert C\vert=a,  \vert B\vert =\tau(a)\neq a$ are not homotopic.  By Lemma \ref{1f559317olbrad}, it suffices to consider the case where
the alphabet $\alpha$ consists of only two letters $a$ and $\tau(a)\neq a$. Then $\tau$ is fixed-point-free  and  we can consider  the characteristic sequences of nanowords determined by  $\alpha_0=\{a\}\subset \alpha$. The computation above shows that the reduced  characteristic sequence of 
$ACAC$ is  $(  a ,
 {a_{\bullet}}  ,     a^2  a_{\bullet}^2 ,  -   {a a_{\bullet}^2 } , -   {a^2 a_{\bullet} }  )$.
It follows from the definitions that
$$v_+ (ABACBC)= (a(s\ast_a s) \ast_{\tau(a)} as)\ast_a (s\ast_a s).$$
The value of   $f:X(s) \to F $  on
$v_+(ABACBC)$ can be  easily computed:
$$f (s\ast_a s)=\underline 1 \ast_a \underline 1= \underline 1 \,\underline {a_\bullet}  \left (\underline {a_\bullet a}  \right )^{-1},\,\,\,  f (a (s\ast_a s))= \underline a\, \underline {aa_\bullet} \left (\underline {aa_\bullet a}  \right )^{-1},$$
$$f (a (s\ast_a s)\ast_{\tau(a)} as)=\underline a\, \underline {aa_\bullet} \left (\underline {aa_\bullet a}  \right )^{-1}  \ast_{\tau(a)} \underline a=\left (\underline {a_\bullet^{-1}}\right )^{-1}
 \underline { a_\bullet^{-1} a} \,\, \underline {a } \left (\underline {a a}  \right)^{-1}   \underline a ,$$
 $$f (v_+ (ABACBC))=
 \left (\underline {a_\bullet^{-1}}\right )^{-1}
 \underline { a_\bullet^{-1} a} \,\, \underline {a } \left (\underline {a a}  \right)^{-1}   \underline a \ast_a 
 \underline 1\, \underline {a_\bullet} \left (\underline {a_\bullet a}  \right )^{-1}$$
 $$=\underline 1\, \underline {a_\bullet} \left ( \underline {a a_\bullet} \right )^{-1} 
 \left ( \underline {1} \right )^{-1}  \underline a \, \, \underline {aa_\bullet}   \left ( \underline {a^2 a_\bullet} \right )^{-1}  \underline {aa_\bullet}\,
 \underline {a^2 a^2_\bullet}  \left ( \underline {a  a_\bullet^2} \right )^{-1}
  \left ( \underline {a  a_\bullet} \right )^{-1}.$$
Therefore the reduced  characteristic sequence of  $ABACBC$ is 
$$  1,  {a_\bullet}, -  {a a_\bullet},
-  {1},   a,  {aa_\bullet} ,-  {a^2 a_\bullet} ,   {aa_\bullet},  {a^2 a^2_\bullet}, -  {a  a_\bullet^2}, -
   {a  a_\bullet} .$$
  It differs from the one of $ACAC$ already in the first term. Therefore $  ABACBC$
is not homotopic to    $ACAC$. As a check of our computations, note that the sum of all terms in the characteristic sequences of these two nanowords is the same, as it should be  because these nanowords are indistinguishable by the invariant $\lambda$.

                     \end{document}